\documentclass[11pt]{article}
\usepackage{latexsym}
\usepackage[all]{xy}
\usepackage{amsmath}
\usepackage{amssymb}
\usepackage{amsfonts}

\newtheorem{thm}{Theorem}[section]
\newtheorem{prop}[thm]{Proposition}
\newtheorem{lem}[thm]{Lemma}
\newtheorem{sublem}[thm]{Sub-lemma}
\newtheorem{df}[thm]{Definition}
\newtheorem{cor}[thm]{Corollary}

\newtheorem{rmk}[thm]{Remark}
\newtheorem{ex}[thm]{Example}

\begin{document}

\title{\textbf{Moduli of objects in dg-categories}}
\bigskip
\bigskip

\author{\bigskip\\
Bertrand To\"en and Michel Vaqui\'e \\
\small{Laboratoire Emile Picard}\\
\small{Universit\'e Paul Sabatier, Bat 1R2} \\
\small{31062 Toulouse Cedex 9, France}}

\date{April 2007}
\maketitle

\begin{abstract}
The purpose of this work is to prove the existence of 
an algebraic moduli classifying objects in a given 
triangulated category. 

To any dg-category $T$ (over some base ring $k$), 
we define a $D^{-}$-stack $\mathcal{M}_{T}$ in the
sense of \cite{hagII}, classifying certain 
$T^{op}$-dg-modules. When $T$ is saturated, 
$\mathcal{M}_{T}$ classifies compact objects in the
triangulated category $[T]$ associated to $T$. The 
main result of this work states that under certain
finiteness conditions on $T$ (e.g. if it is saturated)
the $D^{-}$-stack 
$\mathcal{M}_{T}$ is locally geometric (i.e. 
union of open and geometric sub-stacks). 
As a consequence we prove the algebraicity of the group 
of auto-equivalences of saturated dg-categories. We also obtain the existence of reasonable
moduli for perfect complexes
on a smooth and proper scheme, as well as complexes
of representations of a finite quiver. 

\vskip .5cm
\begin{center} \textbf{R\'esum\'e} \end{center}

L'objectif de ce travail est de d\'emontrer l'existence 
d'un espace des modules alg\'ebrique classifiant les objets 
dans une cat\'egorie triangul\'ee donn\'ee. 

Pour toute dg-cat\'egorie $T$ (au-dessus d'un anneau de base $k$), 
nous d\'efinisons un $D^{-}$-champ $\mathcal{M}_{T}$ dans le sens de \cite{hagII}, 
qui classifie certains $T^{op}$-dg-modules. 
Quand $T$ est satur\'ee, $\mathcal{M}_{T}$ classifie les objets compacts 
dans la cat\'egorie triangul\'ee $[T]$ associ\'ee \`a $T$. 
Le r\'esultat principal de ce travail \'enonce que sous certaines conditions 
de finitude sur $T$ (par exemple si $T$ est satur\'ee) le $D^{-}$-champ 
$\mathcal{M}_{T}$ est localement g\'eom\'etrique (i.e. r\'eunion de 
sous-champs ouverts g\'eom\'etriques). 
Comme cons\'equence nous d\'emontrons l'alg\'ebricit\'e du groupe des 
auto-\'equivalences des dg-cat\'egories satur\'ees. Nous obtenons 
aussi l'existence d'un espace des modules raisonnable pour les 
complexes parfaits sur un sch\'ema propre et lisse, ainsi que pour les complexes 
de repr\'esentations de carquois finis.   
\end{abstract}

\medskip

\tableofcontents

\bigskip

\section*{Introduction}

This is the first part of a research project whose purpose is
the study of moduli spaces (or rather stacks) of objects in a triangulated category of
geometric or algebraic origin (for instance the bounded derived category of 
coherent sheaves on an algebraic variety, or the bounded derived category of 
finite dimensional modules over an associative algebra).
The purpose of the present work is to prove the existence of a reasonable
moduli of objects in a given triangulated category satisfying some finiteness assumptions
and possessing a dg-enhancement in the sense of \cite{bk}. 

\bigskip

\begin{center} \textit{The notion of moduli of objects in a category} \end{center} 

Before describing the main definitions and results of this work we will 
briefly review the notion of the moduli of objects in a (non triangulated) linear 
category $C$.  
For such a category the problem consists
of finding an algebraic variety $m_{C}$, or more generally an algebraic space or Artin stack,  
whose points are in a natural bijection with the
set of isomorphism classes of objects in $C$. A first problem in order to 
construct $m_{C}$ is to define the associated moduli functor, or in other words to 
define the notion of a family of objects in $C$ parametrized by a given scheme $X$. 
When the moduli functor is defined a second problem consists of 
proving its representability, by a scheme, an algebraic space or an Artin stack, 
which is only expected to hold when the category $C$ satisfies some
futher finiteness assumptions. When it exists, such a geometric object $m_{C}$
will be called a \emph{moduli space of objects in $C$}. 

A well understood example is when $C=A-mod$ is the category of finite dimensional modules
over a finitely generated associative algebra $A$ (over some
algebraically closed field $k$). The corresponding moduli functor is then the
functor sending an affine scheme $Spec\, B$ to the groupoid of 
$A\otimes B$-modules which are projective of finite rank over $B$.
A natural description of the moduli  is then given 
by the following quotient
$$\coprod_{n\in \mathbb{N}}Hom_{k-alg}(A,M_{n}(k))/Gl_{n}(k),$$
where $M_{n}(k)$ is the algebra of $n$ by $n$ $k$-valued matrices.
One important observation here is that the action of $Gl_{n}(k)$ on 
$Hom_{k-alg}(A,M_{n}(k))$ is not free, as its stabilizers correspond to
automorphism groups of objects in $C$. In particular, it is well known 
that the above quotient does not exist as an algebraic variety, but only
as an algebraic (1-)stack in the sense of Artin. In section \S 1 we will review 
this example in more details for the reader's convenience.

The main question adressed in this work is to generalize the previous
construction of the moduli $m_{T}$ to the case where $T$ is now a
triangulated category. 
Again, the problems of defining a reasonable moduli functor, as well as proving the existence of the 
moduli $m_{T}$ itself (as a higher stack, see below)  are only expected to have a solution
when $T$ does satisfy some further finiteness conditions. By analogy with the non-triangulated case, 
these finiteness assumptions state that $T$ is equivalent (as a triangulated category) to
$D_{parf}(B)$, the derived category of finite dimensional dg-modules over some dg-algebra
$B$ of \emph{finite type} (see \ref{d3}). Moreover, we will see below that the moduli $m_{T}$ 
can reasonably exists only in the setting of higher algebraic stacks. 
>From a intuitive point of view our main result can be 
expressed as follows. 

\begin{thm}\label{ti-1}
Let $T$ be a triangulated category which is triangulated equivalent to 
$D_{parf}(B)$, the derived category of finite dimensional (dg-)modules over a
(dg-)algebra $B$ of finite type. Then, there exists an algebraic $\infty$-stack
(see below) $m_{T}$ whose points are in one-to-one correspondence with
isomorphism classes of objects in $T$.
\end{thm}

An important remark is that 
the derived category $T$ of an abelian category $C$ might satisfy these assumptions even though 
$C$ itself does not (a typical example is $C=Coh(X)$ the category of coherent sheaves on 
a smooth proper variety $X$). Therefore, our main result not only 
implies the existence of a moduli for objects in a triangulated category, but also 
provides a new method to construct moduli for objects in a given abelian category
(as a moduli subspace in some $m_{T}$).

\bigskip

\begin{center} \textit{Dg-categories and higher stacks} \end{center} 

Working with triangulated categories brings several technical problems, and 
it turns out that the theory of triangulated categories in  
too coarse for our purpose (already the problem of defining a reasonable moduli
functor does not seem to have any solution in the setting of triangulated categories). 
A much more flexible theory 
is obtained by replacing triangulated categories by \emph{dg-categories} whose theory 
behaves much smoothly (see \cite{ke,to}) and which will be used as a ground setting in this work.

A dg-category is a generalization of a linear
category in the sense that it consists of a set of objects together with
complexes of morphisms between two objects, and composition maps 
preserving the linear and differential structures. To any such dg-category $T$, 
we define its homotopy category $[T]$, which is the category
having the same objects as $T$ and with
$H^{0}(T(x,y))$ as set of morphisms between $x$ and $y$. An important remark is
that all triangulated categories of geometric or algebraic origin are of the
form $[T]$ for some natural dg-category $T$. Moreover, the triangulated 
structure on $[T]$ is then completely determined by $T$ alone (see e.g. \cite{bk,ke}). 
The structure of a dg-category can then be considered as a refinement of the
structure of a triangulated category. 

Let now $T$ be a fixed dg-category.
An object $x\in T$ possesses an automorphism group $Aut(x)$ which is by definition
the group of automorphisms of $x$ as an object in the genuine category $[T]$.
It is also possible to define \emph{higher automorphisms} groups by 
defining the \emph{$i$-automorphism
groups of $x$} to be $Ext^{1-i}(x,x):=H^{1-i}(T(x,x))$ (for all $i>1$).
This definition can be justified as follows: 
the topological intuition suggests that a dg-category is some
sort of higher categorical structure, whose objects are the objects of $T$, 
whose 1-morphisms from $x$ to $y$ are zero cycles in $T(x,y)$, whose
2-morphisms from $u : x\rightarrow y$ to $v : x\rightarrow y$ are 
elements $h\in T(x,y)_{-1}$ such that $dh=v-u$, and more generally whose
$i$-morphisms are elements in $T(x,y)_{1-i}$. Following this intuition, 
the group of \emph{$2$-automorphisms}
of an object $x$ is the group of homotopy classes of self-homotopies of the 
identity  $id\in T(x,x)_{0}$, which is nothing else than 
$Ext^{-1}(x,x)$. More generally, the group of $i$-automorphisms of 
$x$ is the group of homotopy classes of self-homotopies of 
$0\in T(x,x)_{1-i}$, which is $Ext^{1-i}(x,x)$.

Because of the presence of the higher automorphisms groups, 
the classifying space of objects in $T$ is
no longer a groupoid (objects plus isomorphisms), but should rather be 
a higher groupoid. A topological model for this groupoid is  
a homotopy type $|T|$ whose set of connected components 
is the set of equivalence classes of objects in $T$, 
and such that 
$$\pi_{1}(|T|,x)\simeq Aut_{[T]}(x,x) \qquad
\pi_{i}(|T|,x)\simeq Ext^{1-i}(x,x) \; \forall \; i>1.$$
In other words, the moduli of objects in $T$ is not 
expected to be representable by an algebraic $1$-stack, but
only by an algebraic $n$-stack, or even $\infty$-stack, 
taking into account all these higher automorphisms groups.
This explains our choice to work within the theory of 
$\infty$-stacks, simply called \emph{stacks} in the sequel, for which a good
model is the homotopy theory of simplicial presheaves
over the big site of affine schemes endowed with
the \'etale topology (i.e.
functors from the category of commutative rings to 
simplicial sets satisfying a certain \'etale descent 
condition, see 
\cite[\S 2.1]{hagII}).

\bigskip

\begin{center} \textit{$D^{-}$-stacks} \end{center} 

We have just seen that for an object $x$ in a dg-category 
$T$, the negative ext groups $Ext^{1-i}(x,x)$, for $i>1$,
encode the higher automorphisms groups of the object $x$.
The positive ext groups $Ext^{i}(x,x)$, for $i>0$, also have 
an interpretation in terms of \emph{derived deformation 
theory}, and thus in terms of \emph{extended} moduli 
spaces of objects in $T$ (in the sense
of the DDT philosophy, see for example \cite{hin,ka,ko,hagdag}
for references on the subject). Indeed, it is 
very reasonable to expect that the tangent space of
the moduli of objects in $T$ at an object $x$ is 
$Ext^{1}(x,x)$. In the same way, it is expected that 
the first obstruction space at the object $x$ is
$Ext^{2}(x,x)$. Following the derived deformation theory
philosophy, $Ext^{i}(x,x)$ should be interpreted as
a higher obstruction space, and more generally the
whole complex $T(x,x)[1]$ should be thought
as the \emph{virtual tangent space} at the point $x$. 
The existence of this virtual tangent space suggests that 
the moduli of objects in $T$ is not only a 
stack but should also comes equipped with 
a virtual, or derived, structure. A good formalism 
for these \emph{derived stacks} is the theory of
$D^{-}$-stacks\footnote{The notation $D^{-}$ is to remind
the negative bounded derived category. As explained
in \cite{hagdag,hagII} $D^{-}$-stacks are 
\emph{stacks relative to the category of negatively 
graded cochain complexes}.} in the sense of
\cite{hagII} (see also \cite{hagdag} for an introduction).
In terms of moduli functors this means that our moduli 
of objects in $T$ is not only a functor from 
commutative rings to simplicial sets, but rather 
a functor from 
simplicial commutative rings to simplicial sets.
As shown in \cite{hagI,hagII} there exists a reasonable
homotopy theory of such functors, as well 
as a notion of \emph{geometric $D^{-}$-stack},
and \emph{$n$-geometric $D^{-}$-stacks}, obtained
by successive gluing of affine objects with respect to some
natural extension of the smooth topology from
commutative rings to simplicial commutative rings. 
These $n$-geometric $D^{-}$-stacks are derived analogs, and
generalizations,  
of algebraic higher stacks of C. Simpson (see \cite{si} and
\cite[\S 2.1]{hagII}),
and are global counter-part of the
formal extended moduli spaces appearing in 
the derived deformation theory program (e.g. as presented in \cite{hin}). 

As the category of commutative rings sits naturally
in the category of
simplicial commutative rings (they are the
constant simplicial objects), the category of $D^{-}$-stacks contains as a full sub-category 
the category of stacks. Furthermore the two notions
of being geometric are compatible with respect to this 
inclusion, and thus geometric $D^{-}$-stacks are 
truly generalizations of geometric stacks and in particular
of Artin's algebraic stacks, algebraic spaces and schemes.  
Moreover, to any
$D^{-}$-stack is associated its truncation, which is a stack
representing its \emph{classical}, or \emph{un-derived}
part, and obtained by killing the non-trivial derived
structure (in terms of moduli functors
by restricting functors on simplicial rings to rings). 
This truncation functor, which 
is right adjoint to the inclusion, preserves
geometric objects, and thus proving that 
a $D^{-}$-stack is geometric also implies that 
its truncation (or its classical part) is itself
geometric. Moreover, these truncations come naturally 
equipped with a \emph{virtual structure sheaf}, remembering 
the existence of the $D^{-}$-stack itself, which is
an important additional structure from which 
virtual constructions (e.g. virtual fundamental classes)
can be obtained.

All along this work we will place ourselves in the
context of $D^{-}$-stacks as we think this is the 
right context in which moduli problems should be considered, 
but also as this will simplify the proof of our
main result which would probably be more complicated
in the underived context as certain 
computations of (co)tangent complexes are involved.

\bigskip

\begin{center} \textit{The $D^{-}$-stack $\mathcal{M}_{T}$} \end{center}

We are now ready to present the main object of this work. 
For a  dg-category $T$ over some base ring $k$, 
we are looking for a geometric $D^{-}$-stack  which classifies
objects in $T$. Our first observation is that it is not reasonable to expect 
the existence of a geometric object classifying 
exactly objects in $T$. On one side, 
a general dg-category $T$ can have too few objects, as for instance 
a deformation of an object in $T$ might not stay in $T$ anymore 
(e.g. if $T$ has a unique object). On the other side, 
$T$ can also have too many objects, such as objects of infinite type. Both of these
facts prevent the moduli of objects in $T$ to be reasonably behaved. 

Our solution to this problem is to first enlarge $T$ to the bigger 
dg-category $\widehat{T}$, consisting of all $T^{op}$-dg-modules 
($\widehat{T}$ is the dg-analog of the presheaf category). The $D^{-}$-stack 
$\mathcal{M}_{T}$ will then be defined in a way that its points are in 
correspondence with  
objects in $\widehat{T}$ satisfying some finiteness conditions called
\emph{pseudo-perfectness} (see below). More precisely, 
we define a functor
$$\begin{array}{cccc}
\mathcal{M}_{T} : & sk-CAlg & \longrightarrow & SSet \\
 & A & \mapsto & Map(T^{op},\widehat{A}_{pe}),
\end{array}$$
where $sk-CAlg$ is the category of simplicial commutative
$k$-algebras, $\widehat{A}_{pe}$ is the dg-category 
of \emph{perfect $A$-modules} (see \S 2.4), and
$Map(T^{op},\widehat{A}_{pe})$ is the mapping space
of the model category of dg-categories. 
Using the main results of \cite{to}, the classical 
points of $\mathcal{M}_{T}$ can be described 
in the following way. For a commutative $k$-algebra $A$, 
the (derived) base change from $k$ to $A$ of $T$ is denoted
by $T\otimes^{\mathbb{L}} A$, which is another dg-category 
over $k$. It is then possible to identified
the simplicial set $\mathcal{M}_{T}(A)$ with the
classifying space of $T^{op}\otimes^{\mathbb{L}} A$-dg-modules
whose underlying complex of $A$-modules is perfect, also
called \emph{pseudo-perfect $T^{op}\otimes^{\mathbb{L}} A$-dg-modules}.
Concretely, this means that there is a bijection between
$\pi_{0}(\mathcal{M}_{T}(A))$ and the set of quasi-isomorphism
classes of pseudo-perfect $T^{op}\otimes^{\mathbb{L}} A$-dg-modules. 
Furthermore, the higher homotopy groups of $\mathcal{M}_{T}(A)$
are given by 
$$\pi_{1}(\mathcal{M}_{T}(A),x)\simeq Aut(x) \qquad
\pi_{i}(\mathcal{M}_{T}(A),x)\simeq Ext^{1-i}(x,x) \; \forall \; i>1,$$
where the automorphisms and ext groups are computed in the
homotopy category of $T^{op}\otimes^{\mathbb{L}} A$-dg-modules.
This shows that $\mathcal{M}_{T}$ classifies algebraic
families of pseudo-perfect $T^{op}$-dg-modules.
These descriptions remain also valid when $A$ is 
now a simplicial commutative $k$-algebra, but using the correct
notions for $T^{op}\otimes^{\mathbb{L}} A$ and perfect complexes
of $A$-modules. 

Under certain finiteness conditions on $T$, 
\emph{properness and smoothness} (see definition \ref{d3}), 
being pseudo-perfect is
equivalent to being perfect, and thus
$\mathcal{M}_{T}$ becomes a moduli for perfect $T^{op}$-dg-modules.
In particular, when $T$ is furthermore saturated (i.e. 
all perfect $T^{op}$-dg-modules are (quasi-)representable, see 
definition \ref{d3}), 
$\mathcal{M}_{T}$ is trully a moduli of objects in $T$. 

Finally, let us mention that the construction 
$T \mapsto \mathcal{M}_{T}$ is shown to be the adjoint 
to the functor sending a $D^{-}$-stack $F$ to the
dg-category $L_{pe}(F)$ of perfect complexes on $F$
(see proposition \ref{p3}).
This provides a universal property satisfied by 
$\mathcal{M}_{T}$, showing for instance that 
for any scheme $X$, a morphism of $D^{-}$-stacks 
$X \longrightarrow  \mathcal{M}_{T}$
is precisely the same thing as a morphism of dg-categories
$T\longrightarrow L_{pe}(X)$. 

\bigskip

\begin{center} \textit{The results} \end{center}

The main theorem of this work states that 
the $D^{-}$-stack $\mathcal{M}_{T}$ is locally 
geometric and locally of finite presentation
(i.e. union of open sub-$D^{-}$-stacks which are
$n$-geometric for some $n$ and locally of finite presentation), 
under the condition that $T$ is \emph{of finite type} (see definition
\ref{d3}). Precisely, being of finite type means that 
the triangulated category of all $T^{op}$-dg-module can be written
as the homotopy category of dg-module over a dg-algebra $B$, which 
is itself finitely presented in the sense of the homotopy theory
of dg-algebras. An important result states that 
a smooth and proper dg-category is always of finite type
(see corollary \ref{c2}), giving a systematic way of
producing finite type dg-categories of geometric origins. \\

The main theorem of this work is then the following. 

\begin{thm}{(Theorem \ref{t1}, Corollary \ref{c4})}\label{ti}
Let $T$ be a dg-category of finite type
(e.g. $T$ saturated). Then, 
the $D^{-}$-stack $\mathcal{M}_{T}$ is \emph{locally geometric and locally of
finite presentation} 
(i.e. union of open $n$-geometric sub-$D^{-}$-stacks locally of finite
presentation). Moreover, for any pseudo-perfect $T^{op}$-dg-module
$E$,  corresponding to a global point of $\mathcal{M}_{T}$, 
the tangent complex of $\mathcal{M}_{T}$ at $E$ is given by
$$\mathbb{T}_{\mathcal{M}_{T},E}\simeq \mathbb{R}\underline{Hom}(E,E)[1].$$
In particular, if $E$ is quasi-representable by an object $x$ in $T$, then
we have
$$\mathbb{T}_{\mathcal{M}_{T},E}\simeq T(x,x)[1].$$
\end{thm}

As a corollary of the above theorem it follows that 
the truncation $t_{0}\mathcal{M}_{T}$, which corresponds
to the  un-derived part of $\mathcal{M}_{T}$, is
itself a locally geometric stack in the sense of C. Simpson's
higher algebraic stacks \cite{si}. However, the object
$\mathcal{M}_{T}$ contains additional and non-trivial informations, 
and is better behaved from an infinitesimal point of view, as  
is shown by the above formula for its tangent complex (note 
for instance that the tangent complex of the truncation 
is not known). 

An interesting easy consequence of our main theorem 
is the following representability result. 

\begin{cor}(Corollary \ref{caut})
Let $T$ be a saturated dg-category over some field $k$. Then, the group
$aut(T)$ of auto-equivalences of $T$ is representable by 
a group scheme locally of finite type over $k$.
\end{cor}

We will present two applications of theorem
\ref{ti}. The first one is for a smooth and proper
scheme $X \longrightarrow Spec\, k$, and $T$ the dg-category 
of perfect complexes on $X$. In this case, 
$\mathcal{M}_{T}$ can be identified with the $D^{-}$-stack
$\mathbb{R}\underline{Perf}(X)$, of perfect complexes
on $X$. The second one is for
$T$ the dg-category of bounded complexes of finite dimensional
representations of a finite quiver $Q$. In this
second case, $\mathcal{M}_{T}$ can be identified with the
$D^{-}$-stack $\mathbb{R}\underline{Perf}(Q)$ of pseudo-perfect complexes
of representations of $Q$, which are nothing
else that complexes of representations whose underlying 
complexes at each node is a perfect complex of $k$-modules. 
When $Q$ has no loops, 
$T$ is furthermore saturated and thus $\mathbb{R}\underline{Perf}(Q)$
is truly the $D^{-}$-stack of perfect complexes
of representations of $Q$.

\begin{cor}{(Corollary \ref{c5} and \ref{c6})}\label{ci}
With the above notations, the $D^{-}$-stacks
$\mathbb{R}\underline{Perf}(X)$ and $\mathbb{R}\underline{Perf}(Q)$
are locally geometric and locally of finite presentation. 
\end{cor}

Passing to the truncations of the $D^{-}$-stacks
$\mathbb{R}\underline{Perf}(X)$ and $\mathbb{R}\underline{Perf}(Q)$
we obtain the following direct consequences, interesting
for their own sake.

\begin{cor}\label{ci2}
\begin{enumerate}
\item The stack of perfect complexes on a smooth and proper
scheme is locally geometric.
\item The stack of pseudo-perfect complexes of representations
of a finite quiver is locally geometric. 
\end{enumerate}
\end{cor}

For a smooth and proper scheme $X$, it is easy to deduce from
the above corollary that the usual 1-stack $\underline{Vect}(X)$, of
vector bundles on $X$ is representable by an Artin 1-stack 
locally of finite type. This provides a new construction of moduli
of vector bundles on smooth proper schemes, which does not
use Artin's representability criterion.

\bigskip

\begin{center} \textit{Future and related works} \end{center}

To start with, we would like to mention that the next
step in our general study of the moduli $\mathcal{M}_{T}$
will be to describe the sub-$D^{-}$-stacks
defined by the additional data of 
t-structures and stability conditions. We expect that any
reasonable t-structure on a saturated 
dg-category $T$ will define an open sub-$D^{-}$-stack
of perverse objects in $\mathcal{M}_{T}$, 
and thus will define a geometric $D^{-}$-stack whose
truncation will be an algebraic 1-stack (e.g. in the sense
of Artin). In the same way, semi-stable objects 
with respect to a reasonable stability condition
should define an open sub-$D^{-}$-stack in $\mathcal{M}_{T}$, 
whose truncation will again be an algebraic 1-stack in the
sense of Artin. Furthermore, a certain moduli of stable
\emph{point-like objects} is expected to be representable by a smooth
and proper scheme $X$ mapping to $\mathcal{M}_{T}$. 
This will provide a general approach 
to the \emph{realization problem}, asking whether
or not it is possible to find a smooth and proper variety
$X$ for which $T$ is equivalent to the dg-category 
$L_{pe}(X)$ of perfect complexes on $X$, those
$X$ being precisely realized as the moduli of 
\emph{point-like objects} in $T$. We think to have a 
partial solution to this problem and we hope to 
come back to this question in a future work. 
A far reaching application to this construction
would be to try to construct the mirror of a Calabi-Yau 
manifold $X$, as a moduli of point-like objects
in the Fukaya category $Fuk(X)$, theorem 
\ref{ti} being the key statement insuring that this moduli
is representable by another Calabi-Yau manifold
$X'$. 

In \cite{to2}, (derived) Hall algebras of dg-categories 
have been introduced from a purely categorical
approach. As for Hall algebras of abelian categories
(see e.g. \cite{jo}), 
it is expected that these derived Hall algebras
possess geometric counter-parts, defined in terms
of convolution product on l-adic sheaves 
on the $D^{-}$-stack $\mathcal{M}_{T}$. A consequence of
theorem \ref{ti} is most probably that such geometric
derived Hall algebras exist, but also that 
they map to the one of \cite{to2} by means 
of a \emph{faisceaux-functions} correspondence and trace formula
for higher stacks. These kind of constructions and technologies are being
studied and will appear in another work. 

Smooth and proper dg-categories have been also studied from the point 
of view of non-commutative geometry ``\`a la Kontsevich'' (see \cite{ko-so}).
The main result of this paper constructs from a smooth and proper 
dg-category, or in other words from a \emph{non-commutative algebraic manifold}, 
a $D^{-}$-stack in the sense of \cite{hagII}, and therefore provides a link
between non-commutative geometry and (commutative) derived algebraic geometry. 

>From the geometricity of the $D^{-}$-stack 
$\mathbb{R}\underline{Perf}(X)$, we easily
deduce a positive answer to the conjecture 5.4
\cite{hagdag}, stating that the $D^{-}$-stack $\mathbb{R}\underline{Vect}_{n}(X)$
of vector bundles on $X$ is geometric (as $\mathbb{R}\underline{Vect}_{n}(X)$ 
is obviously an open sub-$D^{-}$-stack of $\mathbb{R}\underline{Perf}(X)$).
We should acknowledge however that this conjecture
has been proved previously by J. Gorski in his thesis
\cite{jan}, by completely different methods
using a derived version of the Quot functor and valid for
non-smooth varieties.  

Our main result and its consequence on the geometricity
of the stack of perfect complexes (corollary \ref{ci2} $(2)$), recovers and provides
new proofs of  
the result of the recent paper \cite{l} (and thus of \cite{i}), 
at least in the smooth case. 
Indeed, the moduli functors considered in \cite{l} 
is clearly an open sub-stack of 
the truncation of $\mathbb{R}\underline{Perf}(X)$
(precisely the one consisiting of perfect complexes
$E$ such that $Ext^{-i}(E,E)=0$ for all $i>0$), 
which is thus geometric by our theorem \ref{ti}.
Contrary to \cite{l} our proof does not use 
Artin's representability criterion, and 
also provides descriptions of sub-stacks of
finite type (see corollary \ref{c6+}), which seem to be
hard to obtain using the methods of \cite{l}. Some further applications
of the existence of these moduli spaces are given in \cite{at,tv}. 
Moreover, our theorem also provides a non-trivial 
derived structure on the stack considered in \cite{l}, 
and thus is a stricly stronger result (again in the smooth
case). As an example of application, this derived structure can be used to 
prove very easily the existence of the virtual fundamental class
used in \cite{th}, as done for instance at the very end of  
\cite{ck} for the stack of stable maps (see \cite[\S 4.4, Ex. 4]{seat} for more details).

Finally, we think it is possible to weaken 
the finite type assumption on the dg-category $T$ and still
being able to prove theorem \ref{ti}.
A possible approach would probably be 
a direct application of Artin's representability
criterion as well as its extension to the context
of $D^{-}$-stacks by J. Lurie (see \cite{lu}). 
However, we think it is not that much more easy 
to check that the conditions of application of
this criterion are satisfied than to prove
directely by hand that $\mathcal{M}_{T}$ is locally geometric.
Moreover, our approach also provides 
a description of the sub-$D^{-}$-stacks of finite
type in $\mathcal{M}_{T}$, which does not seem easy to 
obtain using Artin's representability criterion.

Another approach to generalize our theorem \ref{ti}
would be to use the theory of $D$-stacks of \cite{hagII}, 
which are generalizations of $D^{-}$-stacks 
($D^{-}$-stacks sits inside $D$-stacks in the same way 
as negatively graded complexes sits inside
unbounded complexes\footnote{The notation
$D$ is to remind the unbounded derived category.}).
Our method of proof of theorem \ref{ti}
shows that $\mathcal{M}_{T}$
is a $1$-geometric $D$-stack (in the sense
of \cite{hagII}) as soon as $T$ has compact generator
(however $\mathcal{M}_{T}$ is not locally finitely presented anymore).
As being geometric as a $D$-stack is strictly weaker
then being geometric as a $D^{-}$-stack,
this does not imply that our theorem \ref{ti}
stays correct under this very weak assumption.
However, the difference between the two notions
is somehow mild, and most probably 
these are equivalent for 
locally finitely presented objects. 

\bigskip

\textbf{Acknowledgments:} We are grateful to 
T. Pantev and L. Katzarkov for motivating discussions
concerning possible applications of stacks
of perfect complexes on algebraic varieties
in order to produce some kind of compactifications
of stacks of vector bundles. The present work 
precisely started when trying to write down a proof
of the geometricity of the stack of perfect
complexes on a scheme. 

We also thank C. Simpson, G. Tabuada and G. Vezzosi for many useful
comments on ealier versions of this manuscript as well
as for their interests. 
\bigskip
\bigskip

\textbf{Conventions:}
All along this work we will fix two universes
$\mathbb{U}\in \mathbb{V}$.
We will always assume that they satisfy the
infinite axiom. By convention, all structures considered, such as
sets, groups, rings  \dots will be
elements of $\mathbb{U}$, except if the contrary
is explicitely stated. Thus, $SSets$ will denote
the category of simplicial sets in $\mathbb{U}$, whereas
$SSet_{\mathbb{V}}$ will denote the category of simplicial sets
in $\mathbb{V}$. The same notations will be used for other
categories as well (e.g. $Sets$, $Sets_{\mathbb{V}}$).
The unique exception is
with the word \emph{categories}, that will not have to  be
elements of $\mathbb{U}$. Categories which are elements of $\mathbb{U}$
will be called \emph{small categories}.

We use the notion of model categories
in the sense of \cite{ho}. The expression
\emph{equivalence} always refer to
weak equivalence in a model category.
For a model category $M$, we will denote by
$Map_{M}$ (or $Map$ if $M$ is clear) its mapping spaces
as defined in \cite{ho}. In the same way, the set of morphisms
in the homotopy category $Ho(M)$ will be denoted
by $[-,-]_{M}$, or by $[-,-]$ if $M$ is clear.
The natural $Ho(SSet)$-tensor structure
on $Ho(M)$ will be denoted
by $K\otimes^{\mathbb{L}}X$, for $K$ a simplicial set
and $X$ an object in $M$. In the same way, the
$Ho(SSet)$-cotensor structure will be denoted
by $X^{\mathbb{R}K}$. The homotopy fiber products will
be denoted by $x\times^{h}_{z}y$, and dually the
homotopy push-outs will be denoted by
$x\coprod^{\mathbb{L}}_{z}y$.

For a simplicial set $X$, we say that a sub-simplicial
set $Y\subset X$ is \emph{full}, if $Y$ is a union of
connected components of $X$. For two 
objects $x$ and $y$ in a model category $M$, 
we denote by $Map^{eq}_{M}(x,y)$ (or 
$Map^{eq}(x,y)$ if $M$ is clear) the full sub-simplicial set
of $Map_{M}(x,y)$ consisting of connected components
corresponding to equivalences (i.e. 
corresponding to isomorphisms in $Ho(M)$ throught the
identification between $\pi_{0}(Map_{M}(x,y))$ and
$[x,y]_{M}$). 

For all along this work, we fix an
associative, unital and commutative ring $k$. We denote
by $C(k)$ the category of
(un-bounded) complexes of $k$-modules. Our complexes
will always be cohomologically indexed 
(i.e. with increasing differentials). The shift
of a complex $E$ is denoted by $E[1]$, and defined
by $E[1]_{n}:=E_{n+1}$ (with the usual sign 
rule $d_{n}^{E[1]}:=-d_{n+1}^{E}$). 
The category $C(k)$ is a
symmetric monoidal model category, using
the projective model structures for which fibrations
are epimorphisms and equivalences are quasi-isomorphisms
(see e.g. \cite{ho}). The monoidal structure on $C(k)$ is the
tensor product of complexes over $k$,  and will be denoted by
$\otimes$. Its derived version will be denoted by $\otimes^{\mathbb{L}}$.

Finally, for $B$ a dg-algebra over $k$, we denote by 
$B-Mod$ its category of left dg-modules. In order to avoid
confusions, when $B$ is a (non-dg) $k$-algebra, the category 
of (non-dg) left $B$-modules will be denoted by $Mod(B)$, whereas
$B-Mod$ will denote the category of left $B$-dg-modules
(or in other words the category or complexes of left
$B$-modules).

\bigskip

\section{The 1-stack of objects in a linear category}

The purpose of this first section is to quickly review a construction of the moduli 
1-stack of objects in a linear category $C$. The construction of the moduli stack $\mathcal{M}_{T}$
of objects in a dg-category $T$ will be based on the same ideas, and we hope this
section will help to understand the definitions and results of the sequel. However,
this part is independant with the rest of the paper and the reader can skip it.  

All along this first section the word \emph{stacks} or \emph{1-stacks} refer
to the usual notion, e.g. as presented in \cite{lm}. \\

We fix a $k$-linear category $C$ and we define a moduli functor
$$m_{C} : k-CAlg \longrightarrow Gpd,$$
from the category of commutative $k$-algebras to the category of 
groupoids, in the following way. For $k'\in k-CAlg$ we set
$$m_{C}(k'):=Fun_{k}(C^{op},Mod(k')^{proj,tf}),$$
the groupoid of $k$-linear functors from $C$ to 
the category of $k'$-modules which are projective of finite type. For
$k' \longrightarrow k''$ a morphism in $k-CAlg$ we have a base change functor
$$k''\otimes_{k'} - : Mod(k')^{proj,tf} \longrightarrow Mod(k'')^{proj,tf},$$
and thus a functor
$$m_{C}(k') \longrightarrow m_{C}(k'')$$
obtained by composition. This defines a (lax) functor $m_{C}$ which is
easily seen to be a 1-stack (on the big \'etale site of affine k-schemes). 

We will say that $C$ is \emph{of finite type} (as a $k$-linear category) 
when the $k$-linear category $\widehat{C}=Mod(C^{op})$ of $C^{op}$-modules
is equivalent to $Mod(B)$ for $B$ a finitely presented associative
$k$-algebra. Equivalently, $C$ is of finite type if it is Morita equivalent to 
a finitely presented associative $k$-algebra $B$ (considered as
a linear category with a unique object).

\begin{thm}\label{tnondg}
With the above notation, and for $C$ a $k$-linear category of finite type,
the 1-stack $m_{C}$ is an Artin stack locally of finite presentation. 
\end{thm}

The proof of the theorem starts by considering a morphism of stacks
$$\pi : m_{C} \longrightarrow m_{k}=\underline{Vect},$$
where $m_{k}=\underline{Vect}$ is the stack of vector bundles, defined
as follows. The choice of an equivalence $\widehat{C}\simeq Mod(B)$, 
induces an equivalence between the stack $m_{C}$ and the
stack $m_{B^{op}}$ of left $B$-modules which are projective and of finite
type over the base $k$. In other words, the groupoid
$m_{C}(k')$ is naturally equivalent to the groupoid of $B\otimes_{k}k'$-modules
which are projective and of finite type as $k'$-modules. We define
$$\pi(k') : m_{C}(k') \longrightarrow \underline{Vect}(k')$$
by sending a $B\otimes_{k}k'$-module  to its underlying $k'$-module. 

It is well known that the stack $\underline{Vect}$ is an Artin stack locally of finite type, 
and more precisely we have
$$\underline{Vect}\simeq \coprod_{n} BGl_{n}.$$
Therefore, to prove theorem \ref{tnondg} it is enough to prove that the morphism
$\pi$ is representable (we use a key fact that 
the domain of a representable morphism is an Artin stack if the codomain is so). 
Let $X=Spec\, k'$ be an affine $k$-scheme, and
$X \longrightarrow \underline{Vect}$ a morphism corresponding to 
a $k'$-module $M$ projective and of finite type. We consider
the stack $m_{C}\times_{\underline{Vect}}X$, which can be identified with the
sheaf of $B$-module structures on $M$
$$\begin{array}{cccc}
\underline{Hom}(B,End(M)) : & k'-CAlg &  \longrightarrow & Sets \\
 & k'' &\mapsto & Hom_{k-Alg}(B,End(M)\otimes_{k'} k'').
\end{array}$$
The algebra $B$ being finitely presented, we can write $B$ as a push-out
$$\xymatrix{B_{1} \ar[r] \ar[d] & B_{2} \ar[d] \\
k \ar[r] & B,}$$ 
where $B_{i}$ are free associative algebras with a finite number of generators
$n_{i}$. The sheaf $\underline{Hom}(B,End(M))$ can then be written as
a pull-back
$$\xymatrix{\underline{End}(M)^{n_{1}} &\ar[l] \underline{End}(M)^{n_{2}} \\
Spec\, k' \ar[u] & \ar[u] \ar[l] \underline{Hom}(B,End(M)),}$$ 
where $\underline{End}(M)$ is the sheaf of endomorphisms of the $k'$-module $M$. 
The module $M$ being projective and of finite type the sheaf 
$\underline{End}(M)$ is representable by an affine $k'$-scheme of finite presentation, and thus
$\underline{Hom}(B,End(M))$ is also an affine $k'$-scheme of finite presentation. This finishes
the sketch of a proof of the theorem. \\

The stack $m_{C}$ is rather far from the starting category $C$ itself in general, as
its groupoid of global points $m_{C}(k)$ can be quite different from 
the groupoid of isomorphisms in $C$ (in general, there is not even a natural
functor between these two groupoids). However, 
under some conditions on $C$ it can be shown to actually classify  
objects in $C$.  For this we will say that $C$ is \emph{locally proper} if its $Hom$'s are 
projective $k$-modules of finite type. We will say that $C$ is 
\emph{smooth} if $C$ is a projective $C\otimes_{k}C^{op}$-modules of finite 
type. We say that $C$ \emph{has a compact and projectif generator} if it 
is Morita equivalent to an associative $k$-algebra $B$. Finally, we say that 
$C$ is \emph{saturated} if it is Karoubian, smooth, proper and
has a compact and projective generator. It can be shown that $C$ is saturated if
and only if  it is
equivalent to the category $Mod(B)^{proj,tf}$, of projective $B$-modules
of finite type, for $B$ an associative $k$-algebra  satisfying the following 
two conditions:

\begin{itemize}

\item $B$ is projective of finite type as a $k$-module.

\item $B$ is a projective $B\otimes_{k}B^{op}$-module of finite type.

\end{itemize}

\begin{rmk}\label{remsat}

The important two properties of saturated $k$-linear categories are the following.

\begin{itemize}

\item \emph{Any saturated $k$-linear category $C$ is of finite type, and thus
$m_{C}$ is an Artin stack.}

\item \emph{The groupoid of global points $m_{C}(k)$ is equivalent to the
groupoid of objects in $C$. More precisely, $m_{C}(k')$ is equivalent 
to the groupoid of $k$'-modules in $C$ (i.e. $k$-linear functors
$k' \longrightarrow C$). Therefore, $m_{C}$ is trully a moduli of 
objects in $C$. }

\end{itemize}
\end{rmk}

In the sequel we will formally follow the same ideas. For this we will 
start by replacing $k$-linear categories by $k$-dg-categories, 
associative $k$-algebras by associative $k$-dg-algebras, modules
by dg-modules \dots and so on. Our first task will be to introduce
analogs of the notions of being of finite type, projective of finite type, 
smooth, proper \dots suited in the dg world. With these notions, the 
definition of the moduli stack $\mathcal{M}_{T}$ for a dg-category $T$ will be very similar 
to the definition of the 1-stack $m_{C}$, and the proof of its algebraicity 
will be very close to the proof of theorem \ref{tnondg} presented above. 

It should be noted that being saturated is a very strong, and thus not very
convenient condition on a linear category. However, the analogous notion
for dg-categories will be much more  interesting as it will be satisfied in many
examples (e.g. an algebra might be saturated when considered as a dg-category 
without being saturated as a linear category).

\section{Preliminaries}

In this section we have gathered definitions and
results that will be used in the next section for the
definition of the $D^{-}$-stack $\mathcal{M}_{T}$ and for
the proof of its geometricity. The first paragraph
on model categories seem to be of a folklore 
knowledge and is essentially extracted from
\cite{hi}. In the second paragraph we review
some terminology and results concerning dg-categories.
Most of them are well known, but some of them 
are probably new. Finally, the two last 
paragraphs are extracted from \cite{hagII}, and
contain the basic language of derived algebraic
geometry we will need in the sequel. The notion
of locally geometric $D^{-}$-stack is the only new notion
which does not appear in \cite{hagII}, but is 
a straighforward generalization of the notion
of $n$-geometric $D^{-}$-stacks. 
 
\subsection{Compactly generated model categories}

We fix a model category
$M$. The purpose of this paragraph is to
consider the relations between
homotopically finitely presented objects and
finite cell objects, under certain conditions
on $M$. \\

Let us fix $I$, a set of generating cofibrations
in $M$.

\begin{df}\label{d1}
\begin{enumerate}
\item An object $X$ is a \emph{strict finite $I$-cell object},
if there exists a finite sequence
$$\xymatrix{
X_{0}=\emptyset \ar[r] & X_{1} \ar[r] & \dots \ar[r] & X_{n}=X,}$$
and for any $0\leq i<n$ a push-out square
$$\xymatrix{
X_{i} \ar[r] & X_{i+1} \\
A \ar[u] \ar[r]_-{u_{i}} & B, \ar[u]}$$
with $u_{i} \in I$.
\item An object $X$ is a \emph{finite $I$-cell object}
(or simply a \emph{finite cell object} when $I$ is clear)
if it is equivalent to a strict finite $I$-cell object.
\item An object $X\in M$ is
\emph{homotopically finitely presented} if for
any filtered system of objects $Y_{i}$ in $M$,
the natural morphism
$$Colim_{i}Map(X,Y_{i}) \longrightarrow Map(X,Hocolim_{i}Y_{i})$$
is an isomorphism in $Ho(SSet)$.
\item
The model category $M$ is \emph{compactly generated},
if it is cellular (in the sense
of \cite[\S 12]{hi}), and
there is a set of generating cofibrations
$I$ and generating trivial cofibrations $J$
whose domains and codomains are cofibrant and
$\omega$-compact, and $\omega$-small (with respect to the whole category $M$).
\end{enumerate}
\end{df}

The main result of this paragraph is the following proposition.
It seems to be a folklore result, but for which we do not know any
reference.

\begin{prop}\label{p1}
Let $M$ be a compactly generated model category, and
$I$  be a set of generating cofibrations
whose domains and codomains are cofibrant, $\omega$-compact
and $\omega$-small with respect to the whole category $M$.
\begin{enumerate}
\item A filtered colimit of fibrations (resp. trivial
fibrations) is a fibration (resp. a trivial fibration).
\item For any filtered diagram $X_{i}$ in $M$, the natural
morphism
$$Hocolim_{i}X_{i} \longrightarrow Colim_{i}X_{i}$$
is an isomorphism in $Ho(M)$.

\item Any object $X$ in $M$ is equivalent to a filtered
colimit of strict finite $I$-cell object.

\item If furthermore, filtered colimits are exact in $M$, an object $X$ in $M$ is homotopically finitely
presented if and only if it is equivalent to a
retract of a strict finite $I$-cell object.

\end{enumerate}
\end{prop}

\textit{Proof:} $(1)$ Note first that by definition $M$ is finitely generated
in the sense of \cite[\S 7]{ho}. Property $(1)$ is then proved in 
\cite[\S 7.4.1]{ho}. \\

$(2)$ As homotopy colimit is the left derived functor
of colimit, it is enough to show that filtered colimits preserve equivalences.
For a filtered category $A$, the functor $colim : M^{A} \longrightarrow M$
is left Quillen for the projective model structure (equivalences and fibrations
are levelwise). By $(1)$, colim preserves trivial fibrations. Therefore, the functor
colim is a left Quillen functor preserving trivial fibrations, and thus preserves
equivalences. \\

$(3)$ The small object argument gives that
any object $X$ is equivalent to a $I$-cell complex $Q(X)$.
By $\omega$-compactness of the domains and codomains of $I$,
$Q(X)$ is the filtered colimit of its finite sub-$I$-cell complexes.
This implies that $X$ is equivalent to a filtered colimit
of strict finite $I$-cell objects. \\

$(4)$ Let $A$ be a filtered category, and $Y \in M^{A}$ be a $A$-diagram.
Let $c(Y) \longrightarrow R_{*}(Y)$ be a Reedy fibrant replacement
of the constant simplicial object $c(Y)$ with values $Y$ 
(in the model category of simplicial objects in $M^{A}$, see \cite[\S 5.2]{ho}).
By $(2)$, the induced morphism
$$c(Colim_{a\in A} Y_{a}) \longrightarrow Colim_{a\in A} R_{*}(Y_{a})$$
is an equivalence of simplicial objects in $M^{A}$. Moreover, 
$(1)$ and the exactness of filtered colimit imply that 
$Colim_{a\in A} R_{*}(Y_{a})$ is a Reedy fibrant object in the model
category of simplicial objects in $M$ (as filtered colimits
commute with matching objects for the Reedy category $\Delta^{op}$, see \cite[\S 5.2]{ho}). 
This implies that for any cofibrant and $\omega$-small object $K$ in $M$, 
we have
$$Hocolim_{a\in A}Map(K,Y_{a})\simeq
Colim_{a\in A}Map(K,Y_{a}) \simeq 
Colim_{a\in A}Hom(K,R_{*}(Y_{a}))\simeq $$
$$Hom(K,Colim_{a\in A}R_{*}(Y_{a}))\simeq Map(K,Colim_{a\in A}Y_{a}).$$
This implies that the domains and codomains of 
$I$ are homotopically finitely presented. 

As filtered
colimits of simplicial sets preserve homotopy pull-backs, we deduce that
any finite cell objects is also homotopically finitely presented,
as they are constructed from domains and codomains of $I$
by iterated homotopy push-outs (we use
here that domains and codomains of $I$ are cofibrant).
This implies that
any retract of a finite cell object is
homotopically finitely presented. Conversely,
let $X$ be a homotopically finitely presented object in $Ho(M)$, and
by $(3)$ let us write it as $Colim_{i}X_{i}$, where
$X_{i}$ is a filtered diagram of finite cell objects. Then,
$[X,X]\simeq Colim_{i}[X,X_{i}]$, which
implies that this identity of $X$ factors through some $X_{i}$, or in other words
that $X$ is a retract in $Ho(M)$ of some $X_{i}$. \hfill $\Box$ \\

Examples of compactly generated model categories
include $C(k)$ the model category of complexes of $k$-modules,
$dg-Cat$ the model category of dg-categories,
and $T-Mod$ the model category of dg-modules over a dg-category $T$
(see next paragraph).

\subsection{dg-categories}

We will use notations and defintions
from \cite{to} and \cite{tab}. We refer the reader to the
overview \cite{ke} for an introduction to the concepts discussed in this
paragraph.\\

The category of (small) dg-categories $dg-Cat$ (over the base
ring $k$) has a model structure for which 
equivalences are the quasi-equivalences (see \cite{tab}). The
model category $dg-Cat$ is furthermore 
compactly generated in the sense of \ref{d1} $(4)$, as this
can be seen easily using the standard 
generating cofibrations and trivial cofibrations
of \cite{tab}. 
We recall from \cite{to}, that the homotopy category 
$Ho(dg-Cat)$ has a natural symmetric monoidal
structure $\otimes^{\mathbb{L}}$, induced by deriving
the tensor product of dg-categories. 
This symmetric monoidal structure is known to be closed, and
its corresponding internal Hom's will be denoted
by $\mathbb{R}\underline{Hom}$ (see \cite{to}). \\

Let $M$ be a $C(k)$-model category which 
we assume to be $\mathbb{V}$-small.
We define a $\mathbb{V}$-small dg-category
$Int(M)$ in the following way. The set of objects
of $Int(M)$ is the set of fibrant and cofibrant objects
in $M$. For two such objects $F$ and $E$
we set
$$Int(M)(E,F):=\underline{Hom}(E,F)\in C(k),$$
where $\underline{Hom}(E,F)$  are the $C(k)$-valued
Hom's of the category $M$.
The set of objects of the dg-category $Int(M)$ 
belongs to $\mathbb{V}$ but not to $\mathbb{U}$ anymore,
and thus $Int(M)$ is only a $\mathbb{V}$-small dg-catgeory.
However, for any $E$ and $F$ in $Int(M)$
the complex $Int(M)(E,F)$ is in fact
$\mathbb{U}$-small.

For $T$ a dg-category, we denote by $T^{op}-Mod$ the category of
$T^{op}$-dg-modules, i.e. the category of $C(k)$-enriched functors
$F : T^{op} \longrightarrow C(k)$.
The category $T^{op}-Mod$ can be endowed with a structure of
model category such that a morphism $f: F \longrightarrow G$ is
an equivalence (resp. a fibration) if for any $z \in T$ the
morphism $f_z : F(z) \longrightarrow G(z)$ is an equivalence
(resp. a fibration) in $C(k)$.
Moreover the natural structure of $C(k)$-module on $T^{op}-Mod$
makes it a $C(k)$-model category in the sense of \cite[4.2.18]{ho}.
The model category $T^{op}-Mod$ is compactly generated in the
sense of definition \ref{d1} $(4)$, as this is easily seen using the 
natural generating cofibrations and trivial cofibrations 
(see e.g. \cite[11.6]{hi}).

To any dg-algebra $B$ (over $k$), we associate 
a dg-category, denoted also by $B$,
with a unique object $*$ and such that 
the dg-algebra of endomorphisms of this object $B(*,*)$, is 
equal to $B$. The category
$B-Mod$ is then naturally isomorphic to the
category of left-$B$-dg-modules.
The dg-category associated to the dg-algebra $k$ itself is denoted by $\mathbf{1}$, 
and is the dg-category with a unique object $*$ and with $\mathbf{1}(*,*)=k$.
The dg-category $\mathbf{1}$ is also the unit for
the monoidal structure $\otimes$ on $dg-Cat$.
We recall our convention that for
a non-dg $k$-algebra, considered as dg-category as above, 
$B-Mod$ is then the category of complexes of left $B$-modules, 
and should not be confused with $Mod(B)$ the category 
of left $B$-modules.  

This construction provides a functor
$$G : dg-Alg \longrightarrow dg-Cat_{*}:=\mathbf{1}/dg-Cat,$$
from the model category of associative and unital dg-algebras
over $k$ (fibrations are
epimorphisms and equivalences are quasi-isomorphisms, 
see e.g. \cite{ss}), to the
model category of pointed dg-categories. This functor
is easily seen to be a left Quillen functor, whose
right adjoint sends a pointed dg-category $T$ to
the dg-algebra $T(t,t)$ of endomorphisms of the 
distinguished object $t$ of $T$. \\

The $C(k)$-enriched version of the Yoneda lemma provides a
morphism of dg-categories
$$\underline{h}_{-} : T \longrightarrow T^{op}-Mod,$$
defined by $\underline{h}_{x}(z) = T(z,x)$ which is quasi-fully faithful.
This morphism is a morphism of $\mathbb{V}$-small dg-categories.
As any object in $T^{op}-Mod$ is fibrant, $Int(T^{op}-Mod)$ is
simply the full sub-$C(k)$-enriched category of $T^{op}-Mod$ consisting
of cofibrant objects, and we will denote it by $\widehat{T}$. In particular, 
the dg-category $\widehat{\mathbf{1}}$ is the dg-category of 
cofibrant complexes of $k$-modules. 
There is a natural isomorphism in $Ho(dg-Cat)$
$$\widehat{T}\simeq \mathbb{R}\underline{Hom}(T^{op},\widehat{\mathbf{1}}).$$
It is easy to see that for any $x\in T$ the object 
$\underline{h}_{x}$ is cofibrant in $T^{op}-Mod$, and thus
the morphism $\underline{h}$ factors as
$$\underline{h}_{-} : T \longrightarrow \widehat{T}.$$
This last morphism is quasi-fully faithful, and by definition
is the Yoneda embedding of the dg-category $T$. 

\begin{df}\label{d2}
Let $T$ be a dg-category. A $T^{op}$-dg-module is 
\emph{perfect}, or \emph{compact}, if
it is homotopicaly finitely presented in 
the model category $T^{op}-Mod$. The full sub-dg-category
of perfect objects in $\widehat{T}$ will be denoted by
$\widehat{T}_{pe}$.
\end{df}

The model category $T^{op}-Mod$ being compactly generated in the sense
of definition \ref{d1} $(4)$,
a $T^{op}$-dg-module $F$ is in $\widehat{T}_{pe}$ if and only if
it is equivalent to a retract of some $F'$, such that there exists a finite
sequence
$$\xymatrix{F_{0}=0 \ar[r] & F_{1} \ar[r] & \dots \ar[r] & F_{i} \ar[r] & F_{i+1} \ar[r] & \dots \ar[r] & F_{n}=F',}$$
and for any $i$ a (homotopy) push-out diagram
$$\xymatrix{
F_{i} \ar[r] & F_{i+1} \\
A\otimes\underline{h}_{x} \ar[u] \ar[r] & B\otimes\underline{h}_{x}, \ar[u]}$$
where $A \rightarrow B$ is a cofibration between bounded complexes of
projective $k$-modules of finite type.
Because of this description, up to a quasi-equivalence,
the dg-category $\widehat{T}_{pe}$ is small. In the sequel we will 
proceed as though
$\widehat{T}_{pe}$ were small. Note that 
for any $x\in T$, the $T^{op}$-dg-module $\underline{h}_{x}$ is perfect, 
and therefore the Yoneda embedding factors as
$$\underline{h}_{-} : T \longrightarrow \widehat{T}_{pe} \subset \widehat{T}.$$
We also note that the above description 
of perfect $T^{op}$-dg-modules implies that 
for a $k$-algebra $B$, considered as a dg-category, 
the perfect $B$-dg-modules are precisely the
perfect complexes of $B$-modules (i.e. the complexes of $B$-modules
which are quasi-isomorphic to bounded complexes of projective $B$-modules
of finite presentation). This justifies 
the terminology \emph{perfect}.

For any dg-category $T$ we denote by $[T]$ the category whose set of
objects is the same as the one of $T$, 
and such that for any $x$ and $y$ in $T$ the set
of morphisms is defined by $[T](x,y):=H^0(T(x,y))$ (with the natural
induced composition maps).
Note that for any $C(k)$-enriched model category $M$, there is a natural
equivalence of categories $[Int(M)]\simeq Ho(M)$ (\cite{to} Proposition 3.5).
The construction $T \mapsto [T]$ is clearly functorial, and provides
a functor from $dg-Cat$ to categories, which furthermore sends
quasi-equivalences to equivalences of categories. 

For a dg-category $T$, the Yoneda embedding
$$\underline{h}_{-} : T \longrightarrow \widehat{T}$$
induces a fully faithful functor
$$\underline{h}_{-} : [T] \longrightarrow 
[\widehat{T}_{pe}] \subset [\widehat{T}]\simeq Ho(T^{op}-Mod).$$ 
As $T^{op}-Mod$ is a stable model category in the sense of \cite[\S 7]{ho},
there exists a natural triangulated structure on
$Ho(T^{op}-Mod)$ whose triangles are the 
image of the homotopy fibration sequences. 
The perfect $T^{op}$-dg-modules are precisely the
compact objects of this triangulated category  in the sense of 
\cite{ne}, 
and thus there exists a natural equivalence of
triangulated categories 
$$[\widehat{T}_{pe}]\simeq [\widehat{T}]_{c}.$$
Moreover, $[\widehat{T}_{pe}] \subset [\widehat{T}]$ is the
smallest thick triangulated sub-category (i.e. 
stable by shifts, cones and retracts) containing 
the image of $\underline{h}$. In other words, 
$[\widehat{T}_{pe}]$ is the thick closure of 
$[T]$ in $[\widehat{T}]$. Using the language of 
dg-categories, this can also be stated 
as the fact that $\widehat{T}_{pe}$ is the smallest
full sub-dg-category or $\widehat{T}$ containing
the essential image of $\underline{h}$, and which 
is stable by retracts, shifts and homotopy 
push-outs. 

\begin{df}\label{d3}
Let $T$ be a dg-category.
\begin{enumerate}
\item $T$ is \emph{locally perfect} (or \emph{locally proper})
if for any two objects $x$ and $y$ in $T$,
$T(x,y)$ is a perfect complex of $k$-modules.
\item $T$ has \emph{a compact generator} if the triangulated
category $[\widehat{T}]$ has a compact generator
(in the sense of \cite{ne}).
\item $T$ is \emph{proper} if $T$ is locally perfect
and has a compact generator.
\item 
$T$ is \emph{smooth} if $T$, considered as
a $T^{op}\otimes^{\mathbb{L}}T$-dg-module
$$\begin{array}{cccc}
T : & T^{op}\otimes^{\mathbb{L}}T & \longrightarrow & C(k) \\
 & (x,y) & \mapsto & T(x,y),
\end{array}$$
is perfect as an object in $\widehat{T\otimes^{\mathbb{L}}T^{op}}$.
\item A dg-category $T$ is \emph{triangulated} if the Yoneda embedding
$$\underline{h} : T \longrightarrow \widehat{T}_{pe}$$
is a quasi-equivalence.
\item A  dg-category
$T$ is \emph{saturated} if it is proper, smooth and triangulated.
\item A dg-category $T$ is \emph{of finite type}
if there exists a dg-algebra $B$, homotopically finitely 
presented in the model category $dg-Alg$, and such that 
$\widehat{T}$ is quasi-equivalent to $\widehat{B^{op}}$.
\end{enumerate}
\end{df}

The notions introduced in the previous definition are extensions to dg-categories of the 
notions we have used in our section \S 1 for linear categories. The main idea was to 
think of the property of being perfect as a dg-version of the property of
being projective and of finite type. Indeed, this analogy if based on the observation
that projective modules of finite type over a commutative ring $k$ are exactly the 
dualizable objects in $Mod(k)$. In the same way, the 
perfect complexes are exactly the dualizable objects in $D(k)$, the homotopy category 
of dg-$k$-modules. With this as a starting point, the above definitions $(1)-(4)$ are 
natural extensions of the notions presented in section \S 1.  

The notion of being triangulated is itself a dg-version of the notion of being
Karoubian. In this sense, $\widehat{T}_{pe}$ is the triangulated hull of $T$ (see lemma below), 
in the same way as the Karoubian hull of a linear category $C$ is  
$\widehat{C}_{proj,tf}$ (the category of projective $C^{op}$-modules of finite type).

The condition of being of finite type for the dg-category $T$ is the key property that will allow us
to prove that the moduli stack $\mathcal{M}_{T}$ is geometric. \\

\begin{ex}\label{ex}
\begin{enumerate}
\item \emph{Having a compact generator for a dg-category $T$ sounds like a
strong condition, as the analog for linear categories is to be Morita equivalent to 
an algebra. The abelian category of quasi-coherent sheaves on a scheme  
can not be written in general as the category of modules over 
an algebra. However, the dg-category of quasi-coherent complexes on 
any quasi-compact quasi-separated scheme is compactly generated (see \cite[Thm. 3.1.1]{bv}). This is
one of the many reasons why dg-categories of complexes of sheaves are in general better behaved than 
abelian category of sheaves. }

\emph{The fact that the dg-categories of complexes of sheaves are compactly generated
seems very specific to algebraic geometry. In fact, the dg-category of complexes
of coherent sheaves on a smooth compact analytic variety is not compactly generated
in general (see \cite[Thm. 5.6.1]{bv}).}

\item \emph{The smoothness notion of definition \ref{d3} is an analog of the
geometric notion for schemes. Indeed, if $A$ is a commutative algebra of finite 
type over a field $k$, then the morphism $Spec\, A \longrightarrow Spec\, k$ is smooth if and only 
if the dg-category $\widehat{A}_{pe}$ of perfect complexes of $A$-modules is 
smooth in the sense of definition \ref{d3}. This notion of smoothness has been used by 
several authors in order to define} non-commutative smooth varieties, \emph{(see e.g.
\cite{ko-so}).}

\item \emph{Our notion of being saturated does not seem exactly equivalent to the one defined
for triangulated categories in \cite{bv}. It is closer to the notion of} homologically 
smooth and compact $A_{\infty}$-algebra \emph{introduced in \cite{ko-so}.}

\emph{We think important to notice that being saturated as a dg-category is much weaker than 
being saturated as a linear category (as defined in \S 1). For instance, 
any finite dimensional $k$-algebra $B$ which is of finite global cohomological dimension
is a saturated dg-category.
In particular, the path algebra of a finite Quiver without oriented loop is 
a saturated dg-category but not a saturated linear category in general. We will also
see that any smooth proper scheme has a saturated dg-category of perfect complexes (see 
\S 3.5). }

\item \emph{In definition \ref{d3} $(7)$, we warn the reader that 
being of finite type must be understood in the derived and non-commutative
sense. This notion is therefore rather far from the geometric notion of
schemes of finite type (it is in fact stronger). Indeed, we will see that a dg-category of finite type is
smooth (see \ref{plisse}). For instance, the commutative $k$-algebra $k[\epsilon]:=k[X]/X^{2}$, of
dual numbers over $k$ is not of finite type as a dg-category.}

\emph{The typical example of homotopically finitely presented dg-algebras are given by 
free dg-algebras (in the associative sense) over a finite number of generators. In fact, 
by proposition \ref{p1} all homotopically finitely presented dg-algebras are obtained
from the free dg-algebras over a finite number of generators by taking a finite number of 
homtopy push-outs and retracts. This fact will be used in the proof of our main theorem
in order reduce the problem to the free case.}
 
\end{enumerate}
\end{ex}

The following lemma gathers some basic results 
concerning smooth and proper dg-categories. 
Some of them are most probably well-known folklore.

\begin{lem}\label{l1}
Let $T$ be a dg-category.
\begin{enumerate}
\item The dg-category $\widehat{T}_{pe}$ is triangulated.
\item The dg-category $T$ has a compact generator (resp. is
proper, resp. is smooth) if and only if $\widehat{T}_{pe}$
has a compact generator (resp. is proper, resp. is smooth).
\item The dg-category $T$ has a compact generator 
if and only if there exists a dg-algebra $B$ such that 
$\widehat{T}$ and $\widehat{B^{op}}=Int(B-Mod)$ are quasi-equivalent. 
\item 
Suppose that $T$ has a compact generator and let 
$B$ be as in $(3)$. Then, $T$ is proper if and only if
the underlying complex of $B$ is perfect.
\item Suppose that $T$ has a compact generator and let 
$B$ be as in $(3)$.
The dg-category $T$ is smooth  if and only if 
$B$ is  perfect as a $B^{op}\otimes^{\mathbb{L}}B$-dg-module.
\end{enumerate}
\end{lem}

\textit{Proof:} $(1)$ From \cite[Lem. 7.5]{to} we have a quasi-equivalence
$$\widehat{(\widehat{T}_{pe})}\simeq
\mathbb{R}\underline{Hom}(\widehat{T}_{pe}^{op},\widehat{\mathbf{1}}) \longrightarrow
\widehat{T}\simeq \mathbb{R}\underline{Hom}(T^{op},\widehat{\mathbf{1}}).$$
Passing to the sub-dg-categories of perfect objects implies the result. \\

$(2)$ For the property of having 
a compact generator this follows by the same argument as above that
$$\widehat{(\widehat{T}_{pe})} \simeq \widehat{T}.$$

Suppose that $T$ is proper. Then, by
for any $x$ and $y$ in $T$ the complex 
$$T(x,y)\simeq \widehat{T}_{pe}(\underline{h}_{x},\underline{h}_{y})$$ 
is perfect. As $[\widehat{T}_{pe}]$ is the thick closure  
of $[T]$, we see that this implies that for any $x$ and $y$
in $\widehat{T}_{pe}$, the complex $\widehat{T}_{pe}(x,y)$
is again perfect. This shows that $\widehat{T}_{pe}$ is 
locally perfect and thus proper. 
Conversely, assume that $\widehat{T}_{pe}$ is proper. As 
$T$ is quasi-equivalent to a full sub-dg-category of
$\widehat{T}_{pe}$ this implies that $T$ is locally perfect. 
Therefore $T$ is proper. 

For a dg-category $T$, being smooth is equivalent to say that 
the Yoneda embedding $\underline{h} : T \longrightarrow
\widehat{T}$, considered as an object in 
$\mathbb{R}\underline{Hom}(T,\widehat{T})\simeq \widehat{T^{op}\otimes^{\mathbb{L}}T}$, 
is perfect. By \cite[\S 7]{to}, we have
$$\mathbb{R}\underline{Hom}(\widehat{T}_{pe},\widehat{(\widehat{T}_{pe})})\simeq 
\mathbb{R}\underline{Hom}(T,\widehat{T}),$$
easily implying that $T$ is smooth if and only if 
$\widehat{T}_{pe}$ is so. \\

$(3)$ If $\widehat{T}$ is quasi-equivalent to $\widehat{B^{op}}$, for
some dg-algebra $B$, then $B$ as a dg-module over itself
is a compact generator of $[\widehat{B^{op}}]\simeq [\widehat{T}]$. 
Conversely, let $E$ be a compact generator of 
$[\widehat{T}]$, and we set $B=\widehat{T}(E,E)^{op}$ be
the opposite dg-algebra
of endomorphisms of $E$. The dg-algebra $B^{op}$, considered
as a dg-category with a unique object can be identified
with the full sub-dg-category of $\widehat{T}_{pe}$ consisting
of the unique object $E$. Therefore, the restricted 
Yoneda embedding provides a morphism
$$\widehat{T} \longrightarrow \widehat{(\widehat{T}_{pe})}
\longrightarrow \widehat{B^{op}}.$$
As $E$ is a compact generator, the induced functor  
$$[\widehat{T}] \longrightarrow [\widehat{B^{op}}]$$
is easily seen to be an equivalence (see for example
\cite{ss2}). This clearly implies that 
$$\widehat{T} \longrightarrow \widehat{B^{op}}$$
is a quasi-equivalence. \\

$(4)$ Suppose that $T$ is a proper dg-category, 
then by $(2)$ so is $\widehat{T}_{pe}$. Therefore, 
if $\widehat{T}\simeq \widehat{B^{op}}$, we find that 
$\widehat{T}_{pe}\simeq \widehat{B^{op}}_{pe}$ is a proper 
dg-category. Therefore, 
$\widehat{B^{op}}_{pe}(B,B)\simeq B$, is a perfect 
complex of $k$-modules. Conversely, 
if $\widehat{T}\simeq \widehat{B^{op}}$, with 
$B$ perfect as a complex of $k$-modules. 
The dg-category with a unique object $B$ is proper, and thus
by $(2)$ so is $\widehat{B^{op}}_{pe}\simeq \widehat{T}_{pe}$. Again
by $(2)$ this implies that $T$ is proper. \\

$(5)$ The proof is the same as for $(4)$, and is 
a consequence of $(2)$.  \hfill $\Box$ \\

Let $T$ and $T'$ be two dg-categories, then we can
associate to any object $E$ in $\widehat{T \otimes^{\mathbb{L}} T'}$
a morphism of dg-categories $F_{E} :T^{op}\longrightarrow \widehat{T'}$
by letting for $x\in T$
$$\begin{array}{cccc}
F_{E}(x) : & (T')^{op} & \longrightarrow & C(k) \\
 & y & \mapsto & E(x,y).
\end{array}$$
This morphism will be considered as a morphism in $Ho(dg-Cat_{\mathbb{V}})$.
By \cite[\S 7]{to}, there exists a natural 
isomorphism in $Ho(dg-Cat)$ 
$$\widehat{T \otimes^{\mathbb{L}} T'}\simeq
\mathbb{R}\underline{Hom}((T\otimes^{\mathbb{L}}T')^{op},\widehat{\mathbf{1}})
\simeq \mathbb{R}\underline{Hom}(T^{op},\widehat{T'}).$$
The morphism $F_{E}$ corresponds to $E$ throught this
identification.

\begin{df}\label{d4}
We say that an object 
$E \in \widehat{T \otimes^{\mathbb{L}} T'}$ is \emph{pseudo-perfect
relatively to $T'$} if the morphism $F_{E}$ factorizes, in $Ho(dg-Cat_{\mathbb{V}})$,
through
$\widehat{T'}_{pe}$:
$$\xymatrix{
T^{op} \ar[rd] \ar[rr]^{F_{E}} && \widehat{T'} \\
& \widehat{T'}_{pe} \ar@{^{(}->}[ru]}$$
The full sub-dg-category of pseudo-perfect $(T \otimes^{\mathbb{L}} T)^{op}$-dg-modules
relatively to $T'$ is denoted by
$$\widehat{T \otimes^{\mathbb{L}} T'}_{pspe} \subset \widehat{T \otimes^{\mathbb{L}} T'}.$$
When $T'=\mathbf{1}$, we will simply use the
terminology \emph{pseudo-perfect $T^{op}$-dg-modules}. 
\end{df}

The main property of smooth and proper dg-categories
is the following lemma, relating perfect and pseudo-perfect
dg-modules. 

\begin{lem}\label{l2}
Let $T$ and $T'$ be two dg-categories.
\begin{enumerate}
\item
If $T$ is locally perfect, then a perfect object $E \in 
\widehat{T\otimes^{\mathbb{L}}T'}$
is pseudo-perfect relative to $T'$.
\item If $T$ is smooth, then
a pseudo-perfect object 
$E \in \widehat{T\otimes^{\mathbb{L}}T'}$ relative to $T'$ is 
perfect.
\item If $T$ is smooth and locally perfect, then 
an object $E \in \widehat{T\otimes^{\mathbb{L}}T'}$
is perfect if and only if it is pseudo-perfect relative
to $T'$.
\end{enumerate}
\end{lem}

\textit{Proof:} $(1)$ As $[\widehat{T'}_{pe}]$ is the
thick closure of $[T']$ in $[\widehat{T'}]$, we see that it 
is enough to show that for any 
$(x,y)\in T\otimes^{\mathbb{L}}T'$, the 
$(T\otimes^{\mathbb{L}}T')^{op}$-dg-module
$\underline{h}_{(x,y)}$ is pseudo-perfect with respect 
to $T'$. But, for $z\in T$, we have
$$\begin{array}{cccc}
\underline{h}_{(x,y)}(z) : & (T')^{op} & \longrightarrow & C(k) \\
 & t & \mapsto & T(z,x) \otimes^{\mathbb{L}} T'(t,y).
\end{array}$$
In other words, the $(T')^{op}$-dg-module 
$\underline{h}_{(x,y)}(z)$ is of the form 
$T(z,x) \otimes^{\mathbb{L}} \underline{h}_{y}$. 
As $\underline{h}_{y}$ is perfect and 
$T(z,x)$ is a perfect complex by assumption on $T$, 
$\underline{h}_{(x,y)}(z)$ is a perfect $(T')^{op}$-dg-module. In other
words $\underline{h}_{(x,y)}$ is pseudo-perfect relative to $T'$. \\

$(2)$ Let us now assume that $T$ is 
smooth, and let $E$ be a pseudo-perfect
$(T\otimes^{\mathbb{L}}T')^{op}$-dg-module relative to $T'$. 
We consider $F_{E}$ as an object in 
$\mathbb{R}\underline{Hom}(T^{op},\widehat{T'})$, as well as
the composition morphism (well defined 
in $Ho(dg-Cat)$)
$$\mathbb{R}\underline{Hom}(T^{op},T^{op}) \otimes^{\mathbb{L}} \mathbb{R}\underline{Hom}(T^{op},\widehat{T'})
\longrightarrow \mathbb{R}\underline{Hom}(T^{op},\widehat{T'}).$$
Evaluating at $F_{E}$ provides a morphism
$$\mathbb{R}\underline{Hom}(T^{op},T^{op}) \longrightarrow 
\mathbb{R}\underline{Hom}(T^{op},\widehat{T'})$$
sending the identity to $F_{E}$. Using \cite[\S 7]{to}, this
last morphism can also be written as
a continuous morphism (i.e.  preserving arbitrary direct sums)
$$\widehat{T\otimes^{\mathbb{L}}T^{op}}\simeq \mathbb{R}\underline{Hom}(T^{op},\widehat{T^{op}})
\longrightarrow
\widehat{T\otimes^{\mathbb{L}}T'},$$
sending $T(-,-)$
to the object $E$. As $T$ is smooth, 
the object $T(-,-) \in [\widehat{T\otimes^{\mathbb{L}}T^{op}}]$
belongs to the smallest thick triangulated category 
containing the objects of the form 
$\underline{h}_{(x,y)}$ for some $(x,y)\in T\otimes^{\mathbb{L}}T^{op}$. 
Therefore, it only remains to prove that the image
of $\underline{h}_{(x,y)}$ is a perfect object in 
$\widehat{T\otimes^{\mathbb{L}}T'}$.  But, by construction this
image is the $(T\otimes^{\mathbb{L}}T')^{op}$-dg-module
sending the object $(a,b)\in T\otimes^{\mathbb{L}}T'$
to $T(a,x)\otimes^{\mathbb{L}}E(y,b)$. 
By assumption on $E$, the $(T')^{op}$-dg-module
$E(y,-)$ is perfect. Therefore, the 
$(T\otimes^{\mathbb{L}}T')^{op}$-dg-module
$T(-,x)\otimes^{\mathbb{L}}E(y,-)$ being the external 
product of two perfect dg-modules is itself perfect.  \\

$(3)$ Follows from $(1)$ and $(2)$. \hfill $\Box$ \\

\begin{cor}\label{c1}
\begin{enumerate}
\item 
Let $T$ be a smooth and proper dg-category with a compact generator, and let
$B$ be a dg-algebra with $\widehat{T}\simeq \widehat{B^{op}}$. Then, an object
$x\in \widehat{T}$ is perfect if and only if
the corresponding $B$-dg-module is perfect as a complex
of $k$-modules.
\item Let $T$ be a saturated dg-category. 
Then an object in $[\widehat{T}]$ is 
pseudo-perfect if and only if it quasi-representable. 
\end{enumerate}
\end{cor}

\textit{Proof:} This is the case $T'=\mathbf{1}$ in lemma \ref{l2} 
$(3)$ and the definition of being saturated. \hfill $\Box$ \\

We will see in the next section that the moduli stack $\mathcal{M}_{T}$
does not classify exactly objects in $T$. However, by point $(2)$ of the above
corollary, this will be the case when $T$ is saturated. This is a dg-analog 
of our remark \ref{remsat}. \\

\begin{lem}\label{l3}
Let $\{T_{\alpha}\}_{\alpha\in A}$ a 
filtered diagram of objects in $dg-Cat$, with
colimit $T$, then the natural morphism
$$Colim _{\alpha\in A} (\widehat{T_{\alpha}})_{pe} \longrightarrow 
\widehat{T}_{pe}$$
is an isomorphism in $Ho(dg-Cat)$.
\end{lem}

\textit{Proof:}
Let $u_{\alpha}$ be the map from $T_{\alpha}$ in $T$,
$$\xymatrix{
T_{\alpha} \ar[r] \ar@/^1pc/[rrr]^{u_{\alpha}} &
T_{\alpha '} \ar[r] & \ldots \ar[r] & T}$$
and let
$$\xymatrix{
T_{\alpha}^{op} -Mod \ar@<2pt>[r]^{{u_{\alpha}}_!} &
T^{op}-Mod \ar@<2pt>[l]^{u_{\alpha}^*}}$$
be the Quillen adjunction defined by $u_{\alpha}$.
This Quillen adjunction induces a morphism
of dg-categories 
$${u_{\alpha}}_{!} : \widehat{T_{\alpha}} \longrightarrow
\widehat{T},$$
and moreover ${u_{\alpha}}_{!}$ preserves the
perfect objects and induces a morphism
$${u_{\alpha}}_{!} : (\widehat{T_{\alpha}})_{pe} \longrightarrow 
\widehat{T}_{pe}.$$
The morphisms ${u_{\alpha}}_{!}$ being functorial 
in $\alpha$ (up to the usual strictification
procedure) 
give a well defined morphism of dg-categories
$$\phi : Colim_{\alpha \in A} (\widehat{T}_{\alpha})_{pe} \longrightarrow 
\widehat{T}_{pe}.$$

Let us start to prove that the morphism 
$\phi$ is quasi-fully faithfull. For this, 
we denote by $T_{0}:=Colim_{\alpha} \widehat{T_{\alpha}}$. 
Let $E$ and $F$ be two objects
in $T_{0}$, and 
we choose $\alpha_{0} \in A$, and two objects
$E_{\alpha_{0}}$, $F_{\alpha_{o}} \in (\widehat{T_{\alpha _{o}}})_{pe}$
representing $E$ and $F$. By definition we have
$$T_{0}(E,F) = Colim_{\alpha \in \alpha_{0}/A} (\widehat{T_{\alpha}})
\bigl ( {i_{\alpha}}_!(E_{\alpha _o}),
{i_{\alpha}}_!(F_{\alpha _o})\bigr ),$$
where $({i_{\alpha}}_! , i_{\alpha}^*)$
is the Quillen adjunction defined by the morphism
$i_{\alpha} : T_{\alpha _o} \longrightarrow T_{\alpha}$, 
for $\alpha \in \alpha _o/A$.
On the other hand
$${\widehat T}(\phi (E), \phi (F)) = {\widehat T} \bigl (
{u_{\alpha _o}}_!(E_{\alpha _o}) ,
{u_{\alpha _o}}_!(F_{\alpha _o}) \bigr ).$$
By adjunction, as $E_{\alpha _o}$ is a perfect object, we have
$$T_{0}(E,F) \simeq \underline{Hom}_{T^{op}_{\alpha_{0}}} \bigl ( E_{\alpha _o} ,
Colim _{\alpha\in \alpha_{0}/A} 
i_{\alpha}^* {i_{\alpha}}_! (F_{\alpha _o}) \bigr )$$
and 
$${\widehat T}(\phi (E), \phi (F)) \simeq \underline{Hom}_{T^{op}_{\alpha_{0}}}
\bigl (E_{\alpha _o} , u_{\alpha _o}^* {u_{\alpha _o}}_! (F_{\alpha _o})
\bigr ),$$ 
where $\underline{Hom}_{T^{op}_{\alpha_{0}}}$ denotes the $C(k)$-enriched $Hom$'s of the
model category $T^{op}_{\alpha_{0}}-Mod$. 
Therefore,  it 
is enough to show that the natural morphism
$$\beta : Colim_{\alpha\in \alpha_{0}/A} i_{\alpha}^*{i_{\alpha}}_! (F_{\alpha _o})
\longrightarrow u_{\alpha _o}^* {u_{\alpha _o}}_! (F_{\alpha _o})$$
is an equivalence in $T_{\alpha _o}^{op}-Mod$.

The two morphisms $Colim_{\alpha} i_{\alpha}^*{i_{\alpha}}_!$ and
$u_{\alpha _o}^* {u_{\alpha _o}}_!$, from 
$[\widehat{T_{\alpha_{0}}}]$
to itself, are triangulated and preserve direct sums.
Furthermore, there exists a natural transformation 
$\beta$ between them, which is easily seen to induce 
isomorphisms when 
evaluated at any representable objects $\underline{h}_{x}
\in [\widehat{T_{\alpha_{0}}}]$. As the representable
objects are compact generators, this implies that 
$\beta$ is an isomorphism. In particular, evaluated at 
the object $F_{\alpha_{0}}$, we get that the natural morphism
$$Colim_{\alpha} i_{\alpha}^*{i_{\alpha}}_! (F_{\alpha _o})
\longrightarrow u_{\alpha _o}^* {u_{\alpha _o}}_! (F_{\alpha _o})$$
is an equivalence. 

We have proved that $\phi$ is fully faithful, and 
it remains to show that $\phi$ is also quasi-essentially surjective.
Let $E$ be a perfect object in $\widehat T$, then $E$ is equivalent
to a retract of an object $F$ which is obtained as the finite colimit
of push-out square
$$\xymatrix{
F_{i} \ar[r] & F_{i+1} \\
A\otimes \underline{h}_{x} \ar[u] \ar[r] &
\ar[u] B\otimes \underline{h}_{x}}$$
for some $x\in T$, and some cofibration $A\rightarrow B$ in $C(k)$ with
$A$ and $B$ bounded complexes of projective modules of finite type
(because of proposition \ref{p1} $(4)$).

As $\phi$ is defined by $\phi ((E_{\alpha})) = {u_{\alpha}}_!(E_{\alpha})$
for some $\alpha$ and is quasi-fully faithful, and as ${u_{\alpha}}_!$
is a left Quillen functor which preserves homotopy push-out and
tensor product by $A\in C(k)$,
it is enough to show that any representable module $\underline{h}^{x}$ in
$[\widehat T]$ is in the image of $[\phi ]$. But, 
this is a consequence of the fact that
${u_{\alpha}}_! (\underline{h}_{x _{\alpha}}) \simeq
\underline{h}_{{u_{\alpha}}(x _{\alpha})}$.
\hfill $\Box$ \\

An important consequence of lemma \ref{l3} is the following finiteness
statement. 

\begin{lem}\label{lll}  
Let $B$ be a dg-algebra, then $B$ is homotopically finitely presented    
in $dg-Alg$ if and only if for any filtered diagram   
$\{T_{\alpha}\}_{\alpha\in A}$ in $dg-Cat$, with colimit $T$, 
the natural morphism 
$$Colim _{\alpha\in A} Map_{dg-Cat}(B,(\widehat{T_{\alpha}})_{pe}) 
 \longrightarrow Map_{dg-Cat}(B,(\widehat{T})_{pe})$$
is an equivalence.  
\end{lem}  

\textit{Proof:}
The functor which associates to a dg-algebra $B$ the
dg-category with one object defines a Quillen adjunction
$$\xymatrix{dg-Alg \ar@<2pt>[r] & dg-Cat_* \ar@<2pt>[l]^E}$$
between the model category of dg-algebras and the model category of
pointed dg-categories $dg-Cat_{*}:=\mathbf{1}/dg-Cat$.
The right adjoint $E$ is defined by sending
a pointed dg-category $T$ to the dg-algebra $T(t,t)$ of
endomorphisms of the distinguished object $t$ of $T$.
Then, for any dg-algebra $B$ and any pointed dg-category $T$ there is
an equivalence $Map_{dg-Alg}(B,E(T)) \simeq Map_{dg-Cat_*}(B,T)$.

Let $B$ be a dg-algebra, then the unit of $B$ defines a morphism
$\mathbf{1} \longrightarrow B$ in $dg-Cat$, and for any dg-category $D$
we obtain a homotopy fibration
$$Map_{dg-Cat _*}(B,D) \longrightarrow Map_{dg-Cat}(B,D)
\longrightarrow Map_{dg-Cat}(\mathbf{1},D)$$
where $D$ is  pointed via the base point
$u : \mathbf{1} \rightarrow D$.
By adjunction, we find a homotopy fibration
$$Map_{dg-Alg}(B,E(D)) \longrightarrow Map_{dg-Cat}(B,D)
\longrightarrow Map_{dg-Cat}(\mathbf{1},D) .$$

Let $D_{\alpha}$ be the dg-category $(\widehat{T_{\alpha}})_{pe}$,
then by lemma \ref{l3} we have $Colim_{\alpha}D_{\alpha}$ equivalent
to $D={\widehat T}_{pe}$ and $Colim_{\alpha}E({D_{\alpha}})$ 
equivalent to $E(D)$, 
and we get the commutative diagram 
 
$$\xymatrix{
Colim_{\alpha}Map_{dg-Alg}(B,E({D_{\alpha}})) \ar[d]^-{a} \ar[r] &
Colim_{\alpha}Map_{dg-Cat}(B,D_{\alpha}) \ar[d]^-{b}\ar[r] &
Colim_{\alpha}Map_{dg-Cat}(\mathbf{1},D_{\alpha}) \ar[d]^-{c} \\
Map_{dg-Alg}(B,E(D)) \ar [r] &
Map_{dg-Cat}(B,D)\ar[r] &
Map_{dg-Cat}(\mathbf{1},D), }$$
in which each horizontal row is a fibration sequence, 
and as $\mathbf{1}$ is homopically finitely presented
in $dg-Cat$, the vertical arrow $c$ is an equivalence.   

If $B$ is homotopically finitely presented in $dg-Alg$, the vertical
arrow $a$ is also an equivalence, then $b$ is an equivalence too. 
Conversely we consider the same diagram with $T_{\alpha}$ the
dg-category associated to a dg-algebra $C_{\alpha}$ then 
$E({T_{\alpha}})=C_{\alpha}$.  
Then if the vertical arrow $b$ is an equivalence we deduce that
$B$ is homotopically finitely presented in $dg-Alg$.  
\hfill $\Box$ \\

\begin{cor}\label{ccc} 
Let $T$ be a dg-category with a compact generator,   
then $T$ is of finite type if and only if any dg-algebra $B$ 
such that $\widehat T$ is quasi-equivalent to $\widehat{B ^{op}}$,  
is homotopically finitely presented in $dg-Alg$.
\end{cor}  

\textit{Proof:} 
We have to show that if $B$ and $C$ are two dg-algebras such 
that the dg-categories $\widehat{B ^{op}}$ and $\widehat{C ^{op}}$  
are equivalent, and if $C$ is homotopically finitely 
presented in $dg-Alg$ then $B$ is also homotopically 
finitely presented.  
This a consequence of the lemma \ref{lll} and of the main result  
of \cite[\S 7]{to}.    
\hfill $\Box$ \\

\begin{cor}\label{c2}
A smooth and proper dg-category $T$ is of finite type.
\end{cor}

\textit{Proof:}  
By the lemma \ref{lll} and the main result of \cite[\S 7]{to} 
it is enough to show that for any filtered diagram 
$\{D_{\alpha}\} _{\alpha\in A}$ in $dg-Cat$, the natural map  
$$Colim _{\alpha\in A} Map_{dg-Cat}(T^{op},(\widehat{D_{\alpha}})_{pe})
 \longrightarrow Map_{dg-Cat}(T^{op},(\widehat{D})_{pe})$$
is an equivalence.

By definition and by results of \cite[\S 7]{to}, 
for any dg-categories $T$ and $T'$ there exists a natural 
equivalence of simplicial sets
$$Map(T^{op}, {\widehat D}_{pe}) \simeq
Map(\mathbf{1}, (\widehat{T\otimes^{\mathbb{L}} D})_{pspe}).$$
Moreover, if $T$ is smooth and proper we deduce from 
lemma \ref{l2} the equivalence  
$$Map(T^{op}, {\widehat D}_{pe}) \simeq
Map(\mathbf{1}, (\widehat{T\otimes^{\mathbb{L}} D})_{pe}).$$

Let $D_{\alpha}$ be a filtered diagram in $dg-Cat$ with colimit $D$,
then we have
$$T\otimes^{\mathbb{L}} D \simeq
T\otimes^{\mathbb{L}} (Colim_{\alpha}D_{\alpha}) \simeq
Colim_{\alpha} (T\otimes^{\mathbb{L}} D_{\alpha}) ,$$
and by lemma \ref{l3}
$$(\widehat{T\otimes^{\mathbb{L}} D })_{pe} \simeq
Colim_{\alpha} ((\widehat{T\otimes^{\mathbb{L}} D_{\alpha}})_{pe}).$$

Next, as $\mathbf{1}$ is homopically finitely presented
in $dg-Cat$, we have
$$\begin{array}{rcl}
Map(T^{op}, {\widehat D }_{pe})  & \simeq &
Map(\mathbf{1}, Colim _{\alpha}
((\widehat{T\otimes^{\mathbb{L}} D_{\alpha}})_{pe})) \\  \\
& \simeq & Colim_{\alpha} Map(\mathbf{1},
(\widehat{T\otimes^{\mathbb{L}} D_{\alpha}})_{pe})
\ \simeq \ Colim_{\alpha} Map(T^{op}, ({\widehat D_{\alpha}})_{pe}).
\end{array}$$
\hfill $\Box$ \\

Corollary \ref{c2} provides a very useful tool to construct 
dg-categories of finite type, as in practice it is often more easy to 
check that a dg-category is saturated than to check directly that 
it is of finite type. \\

Corollary \ref{c2} has the following converse. 

\begin{prop}\label{plisse}
Any dg-category of finite type is smooth.
\end{prop}

\textit{Proof:} It is enough to show that for any
homotopically finitely presented dg-algebra $B$, 
$B$ is homotopically finitely presented as a
$B\otimes^{\mathbb{L}}B^{op}$-dg-module. For this, 
we consider the exact triangle of $B\otimes^{\mathbb{L}}B^{op}$-dg-modules
$$\xymatrix{I_{B} \ar[r] & B\otimes^{\mathbb{L}}B^{op} \ar[r] & B,}$$
where the morphism of the right hand side is given by multiplication in $B$. 
This triangle shows that it is enough to show that 
$I_{B}$ is homotopically finitely presented as a $B\otimes^{\mathbb{L}}B^{op}$-dg-module. 
But, for any $B\otimes^{\mathbb{L}}B^{op}$-dg-module $M$ we have
$$Map_{B\otimes^{\mathbb{L}}B^{op}-Mod}(I_{B},M)\simeq Map_{dg-Alg/B}(B,B\oplus M),$$
where $B\oplus M$ is the trivial square zero extension of $B$ by $M$
(see \cite{la}). 
Therefore, the fact that $I_{B}$ is homotopically finitely presented follows
from the fact that $B$ is homotopically finitely presented as a dg-algebra. \hfill $\Box$ \\

\subsection{Geometric $D^{-}$-stacks}

We will use the theory of $D^{-}$-stacks, as presented
in \cite[\S 2.2]{hagII}, from which we recall the basic
notations.

We let $sk-CAlg$ be the category of simplicial commutative
(associative and unital) $k$-algebras. By definition
$k-D^{-}Aff$ is the opposite category of $sk-CAlg$.
The category  $sk-CAlg$ is endowed with
a model structure for which fibrations and equivalences
are defined on the underlying simplicial sets. The category
$k-D^{-}Aff$ is endowed with the opposite model structure.
For any $A\in sk-CAlg$, we denote by $A-Mod_{s}$
the category of simplicial modules over the
simplicial ring $A$. It is endowed with its natural
model structure for which fibrations and equivalences
are defined on the underlying simplicial sets.

We define natural extensions of the notions
of \'etale, smooth and flat morphisms to the
case of simplicial algebras as follows. A
morphism $A \longrightarrow B$ in $sk-CAlg$
is called \emph{\'etale} (resp. \emph{smooth}, resp.
\emph{flat}) if

\begin{itemize}

\item The induced morphism $\pi_{0}(A) \longrightarrow
\pi_{0}(B)$  is an \'etale (resp. smooth, resp. flat)
morphism of commutative rings.

\item For any $i$, the induced morphism
$$\pi_{i}(A)\otimes_{\pi_{0}(A)}\pi_{0}(B) \longrightarrow
\pi_{i}(B)$$
is an isomorphism.

\end{itemize}

A morphism $A \longrightarrow B$ will be called
an \'etale covering if it is \'etale and if furthermore
the morphism of affine schemes $Spec\, \pi_{0}(B) \longrightarrow
Spec\, \pi_{0}(A)$ is surjective. This notion
endows the model category $k-D^{-}Aff$ with a structure
of a ($\mathbb{V}$-small) model site, and from
the general theory of \cite{hagI} we can construct
a corresponding model category of stacks $k-D^{-}Aff^{\sim,et}$.
The underlying category of $k-D^{-}Aff^{\sim,et}$ is the category
of functors $sk-CAlg \longrightarrow SSet_{\mathbb{V}}$, from
commutative simplicial $k$-algebras to $\mathbb{V}$-small
simplicial sets. The model structure on $k-D^{-}Aff^{\sim,et}$
is a certain left Bousfield localization of the projective
levelwise model category structure on simplicial presheaves,
with respect to equivalences of simplicial rings and
nerves of \'etale hyper coverings.

By definition, the homotopy category of $D^{-}$-stacks
is $Ho(k-D^{-}Aff^{\sim,et})$, and is denoted 
as in \cite{hagII} by $D^{-}St(k)$. In the same way, for
any object $F \in D^{-}St(k)$, we will set
$$D^{-}St(F):=Ho(k-D^{-}Aff^{\sim,et}/RF),$$
were $RF$ is a fibrant replacement of $F$ (note that 
$D^{-}St(F)$ is not equivalent to the $D^{-}St(k)/F$).

For any $A\in sk-CAlg$, we define its spectrum
$\underline{Spec}\, A \in D^{-}St(k)$, by the formula
$$\begin{array}{cccc}
\underline{Spec}\, A : & sk-CAlg & \longrightarrow & SSet \\
 & B & \mapsto & \underline{Hom}(A,B),
\end{array}$$
where $\underline{Hom}(A,B)$ are the natural
simplicial $Hom$'s of the category $sk-CAlg$. The functor
$A \mapsto \underline{Hom}(A,B)$ is almost a right Quillen
functor\footnote{This functor is not right Quillen because it 
does not possess a left adjoint (because of universes issues). However, 
it does preserve fibrations and trivial fibrations and thus we can define 
a total right derived functor.}
$$\underline{Spec} : sk-CAlg^{op} \longrightarrow k-D^{-}Aff^{\sim,et},$$
and its right derived functor
$$\mathbb{R}\underline{Spec} : Ho(sk-CAlg)^{op}=Ho(k-D^{-}Aff) \longrightarrow
D^{-}St(k)$$
is fully faithful. The $D^{-}$-stacks in the essential image
of the functor $\mathbb{R}\underline{Spec}$ are
by definition the \emph{representable $D^{-}$-stacks}.
The notion of \'etale, smooth and flat morphisms
between simplicial $k$-algebras extend in a unique
way to morphisms between representable $D^{-}$-stacks.
A $D^{-}$-stack is then called \emph{$(-1)$-geometric} if
it is isomorphic in $D^{-}St(k)$
to a representable $D^{-}$-stack.

For an integer $n\geq 0$, an
object $F\in D^{-}St(k)$ is an \emph{$n$-geometric
$D^{-}$-stack} if it satisfies the following
two conditions.

\begin{itemize}

\item The morphism $F \longrightarrow F\times^{h} F$
is $(n-1)$-representable.

\item There exists a family of
representable $D^{-}$-stacks $\{X_{i}\}$
and a covering
$$\coprod_{i} X_{i} \longrightarrow F$$
such that each morphism $X_{i} \longrightarrow F$
is smooth. Such a family
$\{X_{i}\}$ together with the morphism
$\coprod_{i}X_{i} \longrightarrow F$
is called an \emph{$n$-atlas for $F$}.
\end{itemize}

In order for the previous notion to make sense
we need to finish the induction on $n$
by setting the following.

\begin{itemize}

\item A morphism $f : F \longrightarrow G$
between $D^{-}$-stacks is $n$-representable
if for any representable $X$ and any morphism
$X \longrightarrow G$, the $D^{-}$-stack
$F\times^{h}_{G}X$ is $n$-geometric.

\item An $n$-representable morphism $f : F \longrightarrow
G$ is smooth, if for any representable $X$ and any morphism
$X \longrightarrow G$, there exists an $n$-atlas
$\{Y_{i}\}$ of the $D^{-}$-stack
$F\times^{h}_{G}X$, such that each morphism
$Y_{i} \longrightarrow X$ is a smooth morphism between
representable $D^{-}$-stacks.

\end{itemize}

We will also use the following terminology.

\begin{itemize}

\item A $D^{-}$-stack $F$ is \emph{quasi-compact} if
there exists a representable $D^{-}$-stack $X$ and
a covering $X \longrightarrow F$.

\item A morphism of $D^{-}$-stacks $f : F \longrightarrow G$
is \emph{quasi-compact} if for any representable
$X$ and any morphism $X \longrightarrow G$, the
$D^{-}$-stack $F\times^{h}_{G}X$ is quasi-compact.

\item By induction on $n$, an $n$-geometric $D^{-}$-stack $F$ is
\emph{strongly quasi-compact}
if it is quasi-compact, and if for any
two representable $X$ and $Y$ and any morphisms
$X \rightarrow F$, $Y \rightarrow F$, the
$(n-1)$-geometric $D^{-}$-stack $X\times^{h}_{F}Y$ is strongly quasi-compact.

\item An $n$-representable morphism $F \longrightarrow G$ is
\emph{strongly quasi-compact} if
for any representable $D^{-}$-stack $X$, and any morphism $X\rightarrow F$,
$F\times^{h}_{G}X$ is a strongly quasi-compact $n$-geometric $D^{-}$-stack.

\item A representable $D^{-}$-stack $F\simeq \mathbb{R}\underline{Spec}\, A$
is
\emph{finitely presented} if for any
filtered system of objects $B_{i}$ in $sk-CAlg$,
the morphism
$$Colim_{i}Map(A,B_{i}) \longrightarrow
Map(A,Colim_{i}B_{i})$$
is an equivalence.

\item An $n$-geometric $D^{-}$-stack $F$ is
\emph{locally of finite presentation} if
it has an $n$-atlas $\{X_{i}\}$ such that
each $X_{i}$ is finitely presented.

\item An $n$-geometric $D^{-}$-stack $F$
is \emph{strongly of finite presentation} if
it is locally of finite presentation and
strongly quasi-compact. 

\end{itemize}

\begin{lem}\label{l4}
The $n$-geometric $D^{-}$-stacks (resp.
$n$-geometric $D^{-}$-stacks locally of finite
presentation, resp. $n$-geometric $D^{-}$-stacks strongly of finite
presentation) are stable by homotopy pull-backs and
retracts.
\end{lem}

\textit{Proof:} The stability by homotopy pull-backs is easily reduced by induction on $n$ to the
case of representable $D^{-}$-stacks, for which the proof is
straighforward (see e.g. \cite[Prop. 1.3.3.3]{hagII}). The stability by retracts 
requires a separate treatment. 

Let $F$ by an $n$-geometric $D^{-}$-stack, and let 
$F_{0}$ be a retract in $D^{-}St(k)$ of $F$. Up to equivalences
we can arrange things so that there exists 
a diagram of fibrant objects in $k-D^{-}Aff^{\sim,et}$
$$\xymatrix{F_{0} \ar[r]^-{i} & F \ar[r]^-{r} & F_{0}}$$
such that $r\circ i=id$ (i.e. that $F_{0}$ is a retract of $F$ 
as a simpicial presheaf). Let $p=i\circ r$ be the corresponding
projector on $F$. As $p^{2}=p$, this defines an action 
of the monoid $M$ freely generated by a projector on $F$. 
We then have an isomorphism in $D^{-}St(k)$
$$F_{0}\simeq lim_{BM}F \simeq Holim_{BM}F$$
where $F$ is considered as a diagram over $BM$, the category
with one object and $M$ as its endomorphisms. The homotopy 
limit on the right hand side can also be described 
as a homotopy limit of a cosimplicial diagram $G_{*}$ whose
object in degree $k$ is $G_{k}:=Hom(M^{k},F)$ and with transitions maps
defined using the action of $M$ on $F$ and the multiplication in $M$. 
The fact that $F_{0}$ is $n$-geometric now follows from 
the following general result.

\begin{sublem}\label{slret}
Let $G_{*} : \Delta \longrightarrow k-D^{-}Aff^{\sim,et}$ be a cosimplicial
diagram such that $G_{k}$ is $n$-geometric for any $k$. Then, 
$Holim_{\Delta}G_{*}$ is also $n$-geometric.
\end{sublem}

\textit{Proof of the sub-lemma \ref{slret}:} Let $\Delta^{*}$ be the 
cosimplicial object in $k-D^{-}Aff^{\sim,et}$ sending $k$ to the 
constant simplicial presheaf $\Delta^{k}$. We have 
$$Holim_{\Delta}G_{*}\simeq Map(\Delta^{*},G_{*})$$
where the mapping spaces on the right are considered in the model category of
cosimplicial simplicial sets (e.g. with the Reedy model
structure of \cite[\S 5.2]{ho}). In other words, for two 
cosimplicial $D^{-}$-stacks $X_{*}$ and $Y_{*}$ we denote
by $Map(X_{*},Y_{*})$ the $D^{-}$-stack sending $A\in sk-CAlg$ to 
$Map_{SSet^{\Delta}}(X_{*}(A),Y_{*}(A))$.

We consider the 
subsimplicial presheaf $\Delta^{*}_{m}$ of $\Delta^{*}$ sending 
$k$ to the $m$-th skeleton $Sq_{m}\Delta^{k}$ of $\Delta^{k}$. 
We clearly have 
$$\Delta^{*}\simeq Colim_{m} \Delta^{*}_{m}$$
and thus
$$Holim_{\Delta}G_{*}\simeq Map(\Delta^{*},G_{*})\simeq Holim_{m}Map(\Delta^{*}_{m},G_{*}).$$
Moreover, for any $m$ we have a homotopy cocartesian square of cosimplicial objects
$$\xymatrix{
\Delta^{*}_{m} \ar [r] & \Delta^{*}_{m+1} \\
h_{m+1}\times \partial \Delta^{m+1} \ar[u] \ar[r] & h_{m+1}\times  \Delta^{m+1} \ar[u],}$$
where $h_{m+1}$ is the cosimplicial object sending $k$ to $\Delta^{k}(m+1)=Hom([m+1],[k])$.
Therefore, there exists a homotopy cartesian square of $D^{-}$-stacks
$$\xymatrix{
Map(\Delta^{*}_{m+1},G_{*}) \ar [r] \ar[d] & Map(\Delta^{*}_{m},G_{*}) \ar[d] \\
G_{m+1} \ar[r] & Map(\partial \Delta^{m+1},G_{m+1}).}$$
>From this description we deduce by induction on $m$ that each $Map(\Delta^{*}_{m},G_{*})$
is an $n$-geometric $D^{-}$-stack. Moreover, we also see by induction on 
$m$ that the morphism
$$Map(\Delta^{*}_{m+1},G_{*}) \longrightarrow Map(\Delta^{*}_{m},G_{*})$$
is $(n-m-1)$-representable, and thus is $(-1)$-representable as soon as 
$m\geq n$ (the case $m=0$ is true because the diagonal
$G_{1} \longrightarrow G_{1}\times^{h} G_{1}$ is $(n-1)$-geometric by definition). 
Using the fact that the $(-1)$-representable $D^{-}$-stacks are stable
by homotopy limits we deduce from this that the morphism
$$Holim_{\Delta}G_{*} \longrightarrow Map(\Delta^{*}_{m},G_{*})$$
is $(-1)$-representable as soon as $m\geq n$. This implies that 
$Holim_{\Delta}G_{*}$ is $n$-geometric. 
\hfill $\Box$ \\

By the previous sub-lemma we know that $F_{0}$ is $n$-geometric. If moreover
$F$ is quasi-compact, the fact that the morphism $r : F \longrightarrow F_{0}$ is
a covering implies that $F_{0}$ is also quasi-compact. By induction on $n$ we also 
see that 
$F_{0}$ is strongly quasi-compact if $F$ is so. To finish the proof
of the lemma it only remains to show that $F_{0}$ is locally of finite presentation
if $F$ is so. We keep the same notations as in sub-lemma \ref{slret}. 
We first consider the truncated stacks $t_{0}(F_{0})$ and $t_{0}(F)$ (see below
or \cite[\S 2.2.4]{hagII}). 
The stack 
$t_{0}(F_{0})$ is now a retract of $t_{0}(F)$, and we can use the same description 
of $t_{0}(F_{0})$ as the homotopy limit of $t_{0}(G_{*})$ as before. But now, 
the stacks  $t_{0}(F_{0})$ and $t_{0}(F)$ being $n$-geometric are also $n$-truncated, 
and thus we see that for $m>n$ the natural map
$$Holim_{\Delta}t_{0}(G_{*}) \longrightarrow Map(\Delta^{*}_{m},t_{0}(G_{*}))$$
is an isomorphism in the homotopy category of stacks. By induction on $m$ we see
that $t_{0}(F_{0})$ is thus locally of finite presentation as a stack. Finally, 
for any $A\in sk-CAlg$, with $X:=\mathbb{R}\underline{Spec}\, A$, and any
morphism $u : X\longrightarrow F_{0}$, the cotangent complex 
$\mathbb{L}_{u}$ is a retract of the cotangent complex $\mathbb{L}_{i\circ u}$
(see \cite{hagII} for the definition of cotangent complexes).
As $F$ is locally of finite presentation the cotangent complex
$\mathbb{L}_{i\circ u}$ is perfect, and thus so is $\mathbb{L}_{u}$. The $D^{-}$-stack
$F_{0}$ is therefore $n$-geometric and such that the truncation $t_{0}(F_{0})$ 
is locally of finite presentation, and its cotangent complexes are
perfect. By \cite[Prop. 2.2.2.4]{hagII} this implies that $F_{0}$ is locally of finite presentation.
\hfill $\Box$ \\

We finish by the following extended definition of geometric
stacks.

\begin{df}\label{d6}
A $D^{-}$-stack $F$ is called \emph{locally geometric}
if it can be written as a filtered colimit
$$F \simeq Hocolim_{i}F_{i}$$
such that the following conditions are satisfied.
\begin{enumerate}
\item Each $D^{-}$-stack $F_{i}$ is $n$-geometric for some
$n$ (depending on $i$).

\item Each morphism $F_{i} \longrightarrow F$ is a monomorphism
(recall this means that the natural morphism
$F_{i} \rightarrow F_{i}\times^{h}_{F}F_{i}$ is an 
isomorphism in $D^{-}St(k)$, or equivalentely that
the morphism induces an isomorphism between $F_{i}$ and
a full sub-$D^{-}$-stack of $F$).

\end{enumerate}
A locally geometric $D^{-}$-stack $F$ is
\emph{locally of finite presentation} if
each of $F_{i}$ as above can be chosen to be
locally of finite presentation.
\end{df}

Note that in the definition above, each morphism
$F_{i} \longrightarrow F$ is a monomorphism and thus 
so are the transition morphisms $F_{i} \longrightarrow F_{j}$.
In particular, the morphisms $F_{i} \longrightarrow F_{j}$
are always formally \'etale. If furthermore $F$
is locally of finite presentation, so that each
$F_{i}$ can itself be chosen
locally of finite presentation, then the morphisms
$F_{i} \longrightarrow F_{j}$ will be
a formally \'etale monomorphism locally of finite presentation,
and thus will be a Zariski open immersion (see 
\cite[Def. 2.2.3.5]{hagII}). In other words, a locally geometric
$D^{-}$-stack which is locally of finite presentation is always
the colimit of a filtered system of Zariski open immersions between
$n$-geometric $D^{-}$-stacks locally of finite presentation. \\

\begin{lem}\label{l5}
Let $F$ be a locally geometric $D^{-}$-stack.
Then, $F$ is $n$-geometric if and only if its diagonal is
$(n-1)$-representable.
\end{lem}

\textit{Proof:} Indeed, by definition it is enough to show that
$F$ an a smooth $n$-atlas. For this, we write
$F$ as a union of $F_{i}$ as in definition \ref{d6}. As each $F_{i}$ is
a full sub-$D^{-}$-stack of $F$ its diagonal is again
$(n-1)$-representable. Therefore, each $F_{i}$ is 
in fact $n$-geometric. We let 
$\{U_{i,j}\}_{j}$ be a smooth $n$-atlas of $F_{i}$, then 
clearly the $\{U_{i,j}\}_{i,j}$ is a smooth 
$n$-atlas for $F$.  \hfill $\Box$ \\

To finish this part, we recall from \cite[\S 2.1]{hagII} that the existence of the homotopy category $St(k)$
of stacks over $k$, being defined as the homotopy category of the model category 
of simplicial presheaves on the site of affine k-schemes endowed with the
\'etale topology. Its objects can be described as functors
$$F : k-CAlg \longrightarrow SSet,$$
from the category $k-CAlg$ of commutative $k$-algebras to the category of simplicial
sets, and satisfying the descent condition for \'etale hypercoverings of affine $k$-schemes. 
There exists a natural embedding $k-CAlg \longrightarrow sk-CAlg$ which consists
of considering a commutative k-algebra as a constant 
simplicial commutative k-algebra. This embedding induces a restriction functor, 
also called the truncation functor
$$t_{0} : D^{-}St(k) \longrightarrow St(k).$$
The functor $t_{0}$ has a left adjoint 
$$i : St(k) \longrightarrow D^{-}St(k)$$
which is fully faithful. Therefore, the 
category $St(k)$ will be considered as a full sub-category 
of $D^{-}St(k)$. However, we warn the reader that the functor
$i$ does not preserve homotopy fiber products except along
flat morphisms. 

As for $D^{-}$-stacks, there exists a notion of $n$-geometric stack 
obtained by using similar definitions (see \cite[\S 2.1]{hagII} for
details). The inclusion functor is moreover compatible with these definitions 
in the sense that a stack $F$ is n-geometric
if and only if the $D^{-}$-stack $i(F)$ is so. By definition, a stack $F$ is
locally geometric if the $D^{-}$-stack $i(F)$ is so. Be careful that 
$i$ does not preserve monomorphisms, so the definition \ref{d6}
must be translated carefully in the language of stacks. In any case, 
a stack is locally geometric and locally finitely presented if and only if
it  can be written as a union of open substacks which are
n-geometric and locally of finite presentation.

By definition, an Artin n-stack is a stack which is n-truncated (i.e. 
for any $A\in k-CAlg$ the simplicial set $F(A)$ is n-truncated) and m-geometric for
some integer $m\geq 0$ (see \cite[2.1]{hagII}). Finally, we mention
the following two results.

\begin{lem}\label{l5'}
Let $F$ be a locally geometric stack which is n-truncated for some
$n\geq 0$. Then $F$ is (n+1)-geometric, and in particular is
an Artin n-stack.
\end{lem}

\textit{Proof:}   This follows easily by induction on $n$. For $n=0$ 
this is \cite[Rem. 2.1.1.5]{hagII}. \hfill $\Box$ \\

\begin{lem}\label{lrep}
Let $F \longrightarrow G$ be an $n$-representable morphism
of $D^{-}$-stacks. If $G$ is $n$-geometric then so is
$F$. 
\end{lem}

\textit{Proof:} This follows from stability of $n$-representable 
morphisms by composition (see \cite[Prop. 1.3.3.3]{hagII}) and from 
the local character of being $n$-geometric (see \cite[Prop. 1.3.3.4]{hagII}). 
\hfill $\Box$ \\

\subsection{Tor amplitude of perfect dg-modules}

Let us fix $A \in sk-CAlg$. We will denote by
$N(A)$ the normalized cochain complex of $k$-modules associated to
$A$. Recall that the shuflle products produce
for two simplicial $k$-modules $X$ and $Y$, a natural morphism
of complexes
$$N(X)\otimes N(Y) \longrightarrow N(X\otimes Y)$$
which is furthermore associative, unital and commutative.
This can be used to 
define a natural structure of a commutative differential
graded $k$-algebra on $N(A)$. This construction provides a functor
$$N : sk-CAlg \longrightarrow k-cdga,$$
which sends equivalences of simplicial algebras to
quasi-isomorphisms of dg-algebras. Furthermore,
for $A\in sk-CAlg$, the normalization provides a functor
$$A-Mod_{s} \longrightarrow N(A)-Mod,$$
from the category of simplicial $A$-modules to the category of
$N(A)$-dg-modules. This last functor also sends equivalences
of simplicial $A$-modules to quasi-isomorphisms. It even induces a
fully faithful functor
$$Ho(A-Mod_{s}) \longrightarrow Ho(N(A)-Mod),$$
whose essential image consists of all $N(A)$-dg-modules
$E$ such that $H^{i}(E)=0$ for $i>0$ (see \cite{ss3}).

Any $A\in sk-CAlg$ possesses a natural augmentation
$A \longrightarrow \pi_{0}(A)$, and thus we can define a base
change morphism
$$Ho(N(A)-Mod) \longrightarrow Ho(\pi_{0}(A)-Mod)\simeq D(\pi_{0}(A)).$$
Note that, by our conventions,
$\pi_{0}(A)$ is considered here as a dg-algebra, and thus
that $\pi_{0}(A)-Mod$ denotes the category of $\pi_{0}(A)$-dg-modules,
which is nothing else than the category of unbounded
complexes of $\pi_{0}(A)$-modules. The category
$D(\pi_{0}(A))$ is thus the unbounded derived category
of $\pi_{0}(A)$.

\begin{df}\label{d7}
An  $N(A)$-dg-module $P$ is
\emph{of Tor amplitude contained in $[a,b]$},
if for any $\pi_{0}(A)$-module (non-dg) $M \in Mod(\pi_{0}(A))$,
we have
$$H^{i}(P\otimes^{\mathbb{L}}_{N(A)}M)=0 \qquad \forall \; i\notin [a,b].$$
\end{df}

The following proposition characterizes the perfect dg-modules
of Tor amplitude contained in $[a,b]$.

\begin{prop}\label{p2}
Let $A\in sk-CAlg$, and $P$ and $Q$ be two perfect $N(A)$-dg-modules.
\begin{enumerate}
\item If $P$ (resp. $Q$) has Tor amplitude contained in
$[a,b]$ (resp. $[a',b']$), then
$P\otimes^{\mathbb{L}}_{N(A)}Q$ is a perfect $N(A)$-dg-module
of Tor amplitude contained in $[a+a',b+b']$.
\item  If $P$ and $Q$ have Tor amplitude contained in
$[a,b]$, then for any morphism $f : P \longrightarrow Q$,
the homotopy fiber of $f$ has Tor aplitude contained
in $[a,b+1]$.
\item The $N(A)$-dg-module $P$ has Tor amplitude contained
in $[a,b]$ if and only if $P\otimes^{\mathbb{L}}_{N(A)}\pi_{0}(A)$
is a perfect complex of $\pi_{0}(A)$-modules whose Tor amplitude
is contained in $[a,b]$.
\item If $A \longrightarrow A'$ is a morphism in $sk-CAlg$, and
if $P$ has Tor amplitude contained in $[a,b]$, then
$P\otimes^{\mathbb{L}}_{N(A)}N(A')$ is a perfect $N(A')$-dg-module
of Tor amplitude contained in $[a,b]$.
\item There exists $a\leq b$ such that
$P$ has Tor amplitude contained in $[a,b]$.
\item If $P$ has Tor amplitude contained in $[a,a]$, then
it is isomorphic in $Ho(N(A)-Mod)$ to
$E[-a]$, for some projective $N(A)$-module of finite type (see \cite[\S 1.2.4]{hagII}).
\item If $P$ has Tor amplitude contained in $[a,b]$, then
there exists a morphism
$$E[-b] \longrightarrow P,$$
with $E$ a projective $N(A)$-dg-module of finite type 
whose homotopy cofiber as Tor amplitude contained in $[a,b-1]$.

\end{enumerate}
\end{prop}

\textit{Proof:} $(1)$, $(2)$, $(3)$, $(4)$ and $(5)$ are all easily deduced
from the definition and the corresponding statement for perfect complexes
of $\pi_{0}(A)$-modules, for which this is well known.

For $(6)$, we can shift and thus reduce to the case
where $a=0$. Let $P$ be a perfect $N(A)$-dg-module
such that the complex of $\pi_{0}(A)$-modules $P\otimes^{\mathbb{L}}_{N(A)}\pi_{0}(A)$ is quasi-isomorphic to
a $\pi_{0}(A)$-module projective and of finite type $E$ (considered as a complex concentrated in degree $0$).
As $E$ is projective, the quasi-isomorphism
$$E[0] \longrightarrow P\otimes^{\mathbb{L}}_{N(A)}\pi_{0}(A)\simeq \pi_{0}(P)$$
can be lifted to a morphism of $N(A)$-dg-modules
$$u : E' \longrightarrow P,$$
with $E'$ a projective $N(A)$-dg-module of finite type. Let 
$C$ be the homotopy cofiber of the morpism $u$, which is a perfect
$N(A)$-dg-module such that 
$$H^{i}(C\otimes^{\mathbb{L}}_{N(A)}\pi_{0}(A))\simeq 0 \; \forall \; i.$$
This easily implies by an induction on $i$ that $H^{i}(C)\simeq 0$
for all $i$. Therefore, $C\simeq 0$ and thus 
$u$ is an equivalence, showing that $P$ is projective and of finite type. 

Finally,
for $(7)$, let $E_{0}$ be a projective $\pi_{0}(A)$-module of finite type and
$$u : E_{0}[-b] \longrightarrow P\otimes^{\mathbb{L}}_{N(A)}\pi_{0}(A)$$
be a morphism in $Ho(\pi_{0}(A)-Mod)$, whose homotopy cofiber
is of amplitude contained in $[a,b-1]$. As $E_0$ is a projective
$\pi_{0}(A)$-module, it is possible to lift $u$ to a morphism
$$v : E[-b] \longrightarrow P,$$
with $E$ a projective $N(A)$-dg-module of finite type.
The homotopy cofiber $C$ of $v$ is a perfect $N(A)$-dg-module, whose base
change $C\otimes^{\mathbb{L}}_{N(A)}\pi_{0}(A)$ has Tor amplitude
contained in $[a,b-1]$. By $(3)$ this implies that $C$ has Tor amplitude
contained in $[a,b-1]$.
\hfill $\Box$ \\

\subsection{Zariski opens and perfect modules}

We fix a simplicial commutative $k$-algebra $A\in sk-CAlg$, and
a perfect $N(A)$-dg-module $P$. We define a sub-$D^{-}$-stack
$V_{P}$ of $\mathbb{R}\underline{Spec}\, A$ as follows. For 
$A \longrightarrow B$ a morphism in $sk-CAlg$, $V_{P}(B)$ is the
full sub-simplicial set of $Map(A,B)$, consisting of morphisms
$u : A\longrightarrow B$ such that $P\otimes_{N(A)}^{\mathbb{L}}N(B)$
is quasi-isomorphic to $0$. The $D^{-}$-stack $V_{P}$ is 
called the \emph{complementary of the support of $P$}. The purpose of this
paragraph is to prove the following proposition, that will be needed
in the sequel. 

\begin{prop}\label{psupp}
The natural inclusion 
$$V_{P} \longrightarrow \mathbb{R}\underline{Spec}\, A$$
is a quasi-compact Zariski open immersion. 
\end{prop}

\textit{Proof:} We consider $K=P\otimes_{N(A)}^{\mathbb{L}}\pi_{0}(A)$, which 
is a perfect complex of $\pi_{0}(A)$-modules. By the semi-continuity of the dimension of 
the cohomology groups of $K$, we know
that there exists a quasi-compact Zariski open sub-scheme $U_{K} \subset Spec\, \pi_{0}(A)$
with the property that for any commutatif $k$-algebra $k'$ and any 
morphism $u : Spec\, k' \longrightarrow Spec\, \pi_{0}(A)$, the complex
$K\otimes_{\pi_{0}(A)}^{\mathbb{L}}k'$ is quasi-isomorphic to zero if and only if
$u$ factors throught $U_{K}$.

Let $(f_{1},\dots,f_{n})$ be elements in $\pi_{0}(A)$ such that 
$U_{K}$ is the union of the standard Zariski open subschemes 
$$U_{i}:=Spec\, \pi_{0}(A)[f_{i}^{-1}] \subset Spec\, \pi_{0}(A).$$
We lift the elements $f_{i}$ to elements $g_{i}\in A_{0}$, and we consider
the localized simplicial rings $A[g_{i}^{-1}]$, e.g. as defined in \cite[\S 1.2.9]{hagII}.
We let $V_{i}:=\mathbb{R}\underline{Spec}\, A[g_{i}^{-1}]$, which is a Zariski open 
sub-$D^{-}$-stack of $\mathbb{R}\underline{Spec}\, A$. We define $V$ to be the quasi-compact
Zariski open sub-$D^{-}$-stack of $\mathbb{R}\underline{Spec}\, A$ which is the
union of the $V_{i}$. 

We claim that the two sub-objects $V_{P}$ and $V$ of $\mathbb{R}\underline{Spec}\, A$ coincide, which 
will imply our proposition. Indeed, a morphism $u : A \longrightarrow B$ in $sk-CAlg$
belongs to $V(B)$ if and only if there exists an index $i$ such that 
the induced morphism $\pi_{0}(A) \longrightarrow \pi_{0}(B)$ sends $f_{i}$ to
an invertible element (see \cite[\S 1.2.9]{hagII}), at least locally 
for the etale topology on $B$. By the choice of the $f_{i}$, this
is also equivalent to the fact that the complex
$$K\otimes_{\pi_{0}(A)}^{\mathbb{L}}\pi_{0}(B)\simeq 
P\otimes^{\mathbb{L}}_{N(A)}\pi_{0}(B)\simeq (P\otimes^{\mathbb{L}}_{N(A)}N(B))
\otimes^{\mathbb{L}}_{N(B)}\pi_{0}(B)$$
is quasi-isomorphic to zero. But, as $P$ is a perfect $N(A)$-dg-module, 
it is easy to check that  $(P\otimes^{\mathbb{L}}_{N(A)}N(B))
\otimes^{\mathbb{L}}_{N(B)}\pi_{0}(B) \simeq 0$ if and only if 
$P\otimes^{\mathbb{L}}_{N(A)}N(B)\simeq 0$. Indeed, $P\otimes^{\mathbb{L}}_{N(A)}N(B)$ 
being perfect is bounded
as a complex of $k$-modules. Therefore, if  $P\otimes^{\mathbb{L}}_{N(A)}N(B)$ were not 
quasi-isomorphic to zero, we chose $j$ to be the maximal index for which 
$H^{j}(P\otimes^{\mathbb{L}}_{N(A)}N(B))\neq 0$, and we would have
$$0\neq H^{j}((P\otimes^{\mathbb{L}}_{N(A)}N(B))
\otimes^{\mathbb{L}}_{N(B)}\pi_{0}(B))\simeq H^{j}(P\otimes^{\mathbb{L}}_{N(A)}N(B))=0.$$
\hfill $\Box$ \\

\section{Moduli $D^{-}$-stacks associated to dg-categories}

This third section is the main body of the paper. 
We will start by the construction of a $D^{-}$-stack
$\mathcal{M}_{T}$ associated to any dg-category $T$, 
and classifying pseudo-perfect $T^{op}$-dg-modules. 
We will also provide a universal property of 
the $D^{-}$-stack $\mathcal{M}_{T}$ by showing that 
the construction $T \mapsto \mathcal{M}_{T}$
has an adjoint sending a $D^{-}$-stack $F$ to 
its dg-category of perfect complexes $L_{pe}(F)$. 
The second paragraph is devoted to the proof of the
fact that $\mathcal{M}_{T}$ is locally geometric
and locally of finite type when $T$ is a smooth and
proper dg-category. In a third paragraph we will
study the sub-$D^{-}$-stack strongly of finite type
of $\mathcal{M}_{T}$, and provide for each choice
of a compact generator of $T$ an exhaustive
family of such sub-stacks. Finally, in the last 
paragraph we present two 
examples of applications, perfect complexes
on a smooth and proper scheme, and 
complexes of representations of a Quiver. 

\subsection{Construction}

For any $A\in sk-CAlg$, we denote by $N(A)$
the $k$-dg-algebra obtained by normalization.
We let
$N(A)-Mod$ be the $C(k)$-model category of (un-bounded)
dg-modules over $N(A)$.
For a morphism $A \longrightarrow B$ in $sk-CAlg$, there is
an induced morphism $N(A) \longrightarrow N(B)$ of dg-algebras,
and an induced base change functor
$$N(B)\otimes_{N(A)} - : N(A)-Mod \longrightarrow N(B)-Mod$$
which is a $C(k)$-enriched left Quillen functor.
For simplicity we will also denote this functor by
$$B\otimes_{A} - : N(A)-Mod \longrightarrow N(B)-Mod.$$
The construction $A \mapsto N(A)-Mod$ is not
functorial in $A$, and is only
a lax functor. However, applying the standard strictification
procedure, we can suppose that, up to an equivalence of categories,
$A \mapsto N(A)-Mod$ is
a genuine functor from $sk-CAlg$ to $C(k)$-model categories
and $C(k)$-enriched left Quillen functor. We will still
denote by $A \mapsto N(A)-Mod$ this strictified object.
We apply the construction $M \mapsto Int(M)$ levelwise,
and thus obtain a presheaf of ($\mathbb{V}$-small) dg-categories
(note that all objects in $N(A)-Mod$ are fibrant, so the base change
functors perserve cofibrant and fibrant objects)
$$\begin{array}{ccc}
sk-CAlg & \longrightarrow & dg-Cat_{\mathbb{V}} \\
A & \mapsto & Int(N(A)-Mod) \\
(A\rightarrow B) & \mapsto & B\otimes_{A}-
\end{array}$$
For any $A$, we let $\widehat{A}_{pe}$ be the full sub-dg-category
of $Int(N(A)-Mod)$ consisting of all perfect objects in
the sense of \S 2.1 (also called homotopically finitely presented).
We first note that perfect objects are stable by the base change functors $B\otimes_{A}-$.
Indeed, this follows formally from the definition of being 
homotopically finitely presented (see \ref{d1}), from the existence of the Quillen adjunction
$$B\otimes_{A}- : N(A)-Mod \longrightarrow N(B)-Mod \\
N(A)-Mod \longleftarrow N(B)-Mod : F,$$
where $F$ is the forgetful functor which commutes with all homotopy colimits. 
We thus obtain a new presheaf of dg-categories
$$\begin{array}{ccc}
sk-CAlg & \longrightarrow & dg-Cat_{\mathbb{V}} \\
A & \mapsto & \widehat{A}_{pe} \\
(A\rightarrow B) & \mapsto & B\otimes_{A}-.
\end{array}$$

We now let $T$ be any dg-category and
$T^{op}$ be its opposite dg-category.
We define a
simplicial presheaf
$$\mathcal{M}_{T} : sk-CAlg \longrightarrow SSet$$
by the formula
$$\mathcal{M}_{T}(A):=Map_{dg-Cat}(T^{op},\widehat{A}_{pe}).$$
Here $Map_{dg-Cat}$ denotes the mapping spaces of
the model categories $dg-Cat$, which as any dg-category is
fibrant will be taken to be
$$Map_{dg-Cat}(T,T'):=Hom(\Gamma^{*}(T),T')$$
where $\Gamma^{*}$ is a co-simplicial resolution functor
in the sense of \cite[\S 5]{ho}. For $A \longrightarrow B$
a morphism in $sk-CAlg$, the transition morphism
$$\mathcal{M}_{T}(A) \longrightarrow \mathcal{M}_{T}(B)$$
is of course defined as the composition
with $B\otimes_{A}- : \widehat{A} _{pe} \longrightarrow \widehat{B} _{pe}$.
This defines the simplicial presheaf $\mathcal{M}_{T}$.

\begin{lem}\label{l6}
The simplicial presheaf $\mathcal{M}_{T}$ is a $D^{-}$-stack
in the sense of \cite[Def. 1.3.2.1]{hagII}.
\end{lem}

\textit{Proof:} Using the properties of the mapping spaces
in model categories we see that it is enough to check the
following three conditions (see
\cite{hagII} for more details).

\begin{itemize}

\item For any equivalence $A \longrightarrow B$ in $sk-CAlg$,
the induced morphism $\widehat{A}_{pe} \longrightarrow \widehat{B}_{pe}$
is a  quasi-equivalence of dg-categories.

\item For any two objects $A$ and $B$ in $sk-CAlg$, the
natural morphism
$$\widehat{A\times B}_{pe} \longrightarrow \widehat{A}_{pe}\times
\widehat{B}_{pe}$$
is a quasi-equivalence of dg-categories.

\item For any \'etale hyper covering $X_{*} \longrightarrow Y$
in $D^{-}Aff$, correspoding to a co-augmented
co-simplicial object
$A \longrightarrow B_{*}$
in $sk-CAlg$, the induced morphism
$$\widehat{A}_{pe} \longrightarrow Holim_{[n]\in \Delta}\left(
\widehat{B_{n}}_{pe}\right)$$
is a quasi-equivalence of dg-categories.

\end{itemize}

The first of this property is satisfied because for an
equivalence $A \longrightarrow B$, the base change
functor $B\otimes_{A}-$ is a Quillen equivalence, so induced
a quasi-equivalence $Int(N(A)-Mod) \longrightarrow Int(N(B)-Mod)$,
and thus a quasi-equivalence on the full sub-dg-categories
of perfect objects. For the second an third properties, 
we first notice that
it is enough to show that the induced morphisms
$$Int(N(A)\times N(B)-Mod) \longrightarrow Int(N(A)-Mod)\times
Int(N(B)-Mod)$$
$$Int(N(A)-Mod) \longrightarrow Holim_{[n]\in \Delta}\left(
Int(N(B_{n})-Mod)\right)$$
are quasi-equivalences, as being
perfect is a local condition for the \'etale
topology (see \cite[Cor. 1.3.7.4]{hagII}). The first of these two quasi-equivalences
is clear.

The fact that
the second one is a quasi-equivalence is a consequence of the main result
of \cite{to} and of the strictification theorem of \cite[Appendix B]{hagII}.
Indeed, we need to prove that for any dg-category $C$,
the morphism of simplicial sets
$$Map(C,Int(N(A)-Mod)) \longrightarrow
Holim_{[n]\in \Delta}\left( Map(C,Int(N(B_{n})-Mod)) \right),$$
is an equivalence.
Applying the main theorem of \cite{to}, the simplicial set
$Map(C,Int(N(A)-Mod))$ is equivalent to
$N((C\otimes^{\mathbb{L}}N(A))-Mod^{cof})$, the nerve of
the category of equivalences between cofibrant $C\otimes^{\mathbb{L}}N(A)$-dg-modules.
In the same way, $Map(C,Int(N(B_{n})-Mod))$ is equivalent to
$N((C\otimes^{\mathbb{L}}N(B_{n}))-Mod^{cof})$,
the nerve of
the category of equivalences between cofibrant $C\otimes^{\mathbb{L}}N(B_{n})$-dg-modules. Finally, the fact that
morphism
$$N((C\otimes^{\mathbb{L}}N(A))-Mod^{cof}) \longrightarrow
Holim_{[n]\in \Delta}\left( N((C\otimes^{\mathbb{L}}N(B_{n}))-Mod^{cof})
\right)$$
is an equivalence follows easily from the
strictification theorem
of \cite[Appendix B]{hagII} and faithfully flat cohomological
descent (see \cite[\S 1.3.7]{hagII} for the similar example
of the stack of quasi-coherent modules). \hfill $\Box$ \\

\begin{df}\label{d8}
The $D^{-}$-stack $\mathcal{M}_{T}$ is called
\emph{the moduli stack of pseudo-perfect $T^{op}$-modules}.
\end{df}

By definition, the simplicial set
$\mathcal{M}_{T}(k)$ is $Map(T^{op},\widehat{k}_{pe})$, which
by \cite[\S 7]{to} is equivalent to the nerve of the
category of equivalences between $T^{op}$-modules $E$, such
that for any $x\in T^{op}$, $E(x)$ is a perfect complex
of $k$-modules. In other words, $\mathcal{M}_{T}(k)$
is a classifying space for pseudo-perfect $T^{op}$-modules
in the sense of definition \ref{d4}. In particular, the set
$\pi_{0}(\mathcal{M}_{T}(k))$ is in natural bijection with
the isomorphism classes of pseudo-perfect $T^{op}$-modules
in $Ho(T^{op}-Mod)$.  Moreover,
for $x \in Ho(T^{op}-Mod)$ pseudo-perfect,
one has natural group isomorphisms
$$\pi_{1}(\mathcal{M}_{T},x)\simeq Aut(x,x) \qquad
\pi_{i}(\mathcal{M}_{T},x)\simeq Ext^{1-i}(x,x),$$
where $Aut$ and $Ext$'s are computed in the triangulated
category $Ho(T^{op}-Mod)$.
This explains our choice of terminology in definition \ref{d8}.

Assume furthermore that $T$ is saturated in the sense of definition
\ref{d3}, then the Yoneda embedding $\underline{h} : T \longrightarrow
\widehat{T}$ is quasi-fully faithful and its image consists
on all pseudo-perfect (or equivalentely perfect) $T^{op}$-modules.
In this case, corollary \ref{c1} implies that the simplicial set
$\mathcal{M}_{T}(k)$ is equivalent to $Map(\mathbf{1},T)$, which
is a model for the classifying space of
objects in $T$. More generally, for any
$A\in sk-CAlg$, the simplicial set
$\mathcal{M}_{T}(A)$ is equivalent to $Map(\mathbf{1},\widehat{T\otimes_{k}^{\mathbb{L}}N(A)}_{pe})$, and
is a model for the classifying space of
perfect $T^{op}\otimes_{k}^{\mathbb{L}}N(A)$-dg-modules. In particular,
the set $\pi_{0}(\mathcal{M}_{T}(A))$ is in natural bijection
with the set of isomorphism classes of perfect objects in
$Ho((T\otimes_{k}^{\mathbb{L}}N(A))^{op}-Mod)$. Furthermore, for
$E$ a perfect $T^{op}\otimes_{k}^{\mathbb{L}}N(A)$-module, 
there is a natural group isomorphisms
$$\pi_{1}(\mathcal{M}_{T}(A),E)\simeq Aut(E,E) \qquad
\pi_{i}(\mathcal{M}_{T}(A),E)\simeq Ext^{1-i}(E,E),$$
where $Aut$ and $Ext$'s are computed in the triangulated
category $Ho((T^{op}\otimes_{k}^{\mathbb{L}}N(A))-Mod)$. Therefore, when
$T$ is saturated, the $D^{-}$-stack $\mathcal{M}_{T}$
classifies objects in $T$.  \\

Clearly, the construction $T \mapsto \mathcal{M}_{T}$ is
contravariantly functorial in $T$, and gives rise
to a functor
$$\begin{array}{ccc}
dg-Cat^{op} & \longrightarrow & k-D^{-}Aff^{\sim,et} \\
T & \mapsto & \mathcal{M}_{T} \\
(u : T \rightarrow T') & \mapsto & \left( u^{*} : \mathcal{M}_{T'} \rightarrow \mathcal{M}_{T} \right),
\end{array}$$
from the model category of small dg-categories to the
model category of $D^{-}$-stacks. For
$u : T \longrightarrow T'$, and $A\in sk-CAlg$, the
morphism $\mathcal{M}_{T'} \longrightarrow \mathcal{M}_{T}$
evaluated at $A$ is simply given by composition with $u$
$$u^{*} : Map((T')^{op},\widehat{A}_{pe}) \longrightarrow
Map(T^{op},\widehat{A}_{pe}).$$
By the properties of the mapping spaces (see \cite[\S 5]{ho}) 
the functor $\mathcal{M}_{-}$ sends quasi-equivalences
between dg-categories to equivalences of $D^{-}$-stacks, and thus
passes throught the homotopy categories
$$\begin{array}{ccc}
Ho(dg-Cat)^{op} & \longrightarrow & D^{-}St(k) \\
T & \mapsto & \mathcal{M}_{T} \\
(u : T \rightarrow T') & \mapsto & \left( u^{*} : \mathcal{M}_{T'} \rightarrow \mathcal{M}_{T} \right).
\end{array}$$
As this last functor is obtained from an equivalence preserving functor,
it is naturally compatible with the $Ho(SSet)$-enrichement
of $Ho(dg-Cat)^{op}$ and $D^{-}St(k)$.
Explicitely, for $T$ and $T'$ two dg-categories, the
morphism in $Ho(SSet)$
$$Map(T,T') \longrightarrow Map(\mathcal{M}_{T'},\mathcal{M}_{T})$$
is adjoint to
$$Map(T,T') \times \mathcal{M}_{T'} \longrightarrow \mathcal{M}_{T},$$
which evaluated at $A\in sk-CAlg$ is given by the composition
of mapping spaces for the model category $dg-Cat$
$$Map(T,T')\times Map((T')^{op},\widehat{A}_{pe}) \longrightarrow
Map(T^{op},\widehat{A}_{pe}).$$

The following lemma shows Morita invariance
of the construction $T \mapsto \mathcal{M}_{T}$.

\begin{lem}\label{l7}
Let $T$ be any dg-category, and let us
consider $\underline{h} : T \longrightarrow \widehat{T}_{pe}$
its Yoneda embedding. Then, the induced morphism
in $D^{-}St(k)$
$$\underline{h}^{*} : \mathcal{M}_{\widehat{T}_{pe}} \longrightarrow
\mathcal{M}_{T}$$
is an isomorphism.
\end{lem}

\textit{Proof:} This is a direct consequence of
the definition of $\mathcal{M}_{T}$ and of
the main results of \cite[\S 7]{to}. \hfill $\Box$ \\

By definition, if $T$ and $T'$ are two $\mathbb{V}$-small 
dg-categories, the simplicial set 
$Map(T,T')$ is itself $\mathbb{V}$-small. In particular, 
the object $\mathcal{M}_{T}$ exists in $D^{-}St(k)$ even 
when $T$ is only a $\mathbb{V}$-small dg-categories. In other words, 
the functor $T \mapsto \mathcal{M}_{T}$ extends to a functor
defined on all $\mathbb{V}$-small dg-categories
$$\mathcal{M}_{-} : Ho(dg-Cat_{\mathbb{V}})^{op} \longrightarrow D^{-}St(k).$$

\begin{prop}\label{p3}
The $Ho(SSet)$-enriched functor
$$\mathcal{M}_{-} : Ho(dg-Cat_{\mathbb{V}})^{op} \longrightarrow D^{-}St(k)$$
has a $Ho(SSet)$-enriched left adjoint
$$L_{pe} : D^{-}St(k) \longrightarrow Ho(dg-Cat_{\mathbb{V}})^{op}.$$
Furthermore, for any $A\in sk-CAlg$, there exists
a natural isomorphism in $Ho(dg-Cat_{\mathbb{V}})$
$$L_{pe}(\mathbb{R}\underline{Spec}\, A)\simeq \widehat{A}_{pe}.$$
\end{prop}

\textit{Proof:} Let $F \in k-D^{-}Aff^{\sim,et}$,
and let us write $F$ as a homotopy colimit
of representable objects
$$F \simeq Hocolim_{i}h_{A_{i}},$$
where $A_{i}\in sk-CAlg$ and $h_{A_{i}} \in k-D^{-}Aff^{\sim,et}$
is defined by $h_{A_{i}}(A):=Hom(A_{i},A)$.
Note that this can be done functorially in $F \in k-D^{-}Aff^{\sim,et}$
by using for instance the standard free resolution of $F$.
The diagram $i\mapsto h_{A_{i}}$ in $k-D^{-}Aff^{\sim,et}$ gives
a diagram $i \mapsto A_{i}$ in $sk-CAlg$ by the Yoneda lemma.
Applying the construction $A \mapsto \widehat{A}_{pe}$ (suitably
strictified) we obtain a diagram $i \mapsto \widehat{A_{i}}_{pe}$
in $dg-Cat$. We set
$$L_{pe}(F):=\left( Holim_{i}\widehat{A_{i}}_{pe} \right)^{op} \in
Ho(dg-Cat_{\mathbb{V}}).$$
By definition, we have for any dg-category $T$
$$Map_{dg-Cat^{op}}(L_{pe}(F),T) \simeq
Holim_{i}Map_{dg-Cat^{op}}(\widehat{A_{i}}_{pe},T^{op})\simeq
Holim_{i}\mathcal{M}_{T}(A_{i})\simeq $$
$$Holim_{i}Map(h_{A_{i}},\mathcal{M}_{T})
\simeq Map(Hocolim_{i}h_{A_{i}},\mathcal{M}_{T})\simeq
Map(F,\mathcal{M}_{T}).$$
These are  isomorphisms in $Ho(SSet)$, and are functorial in $T$.
This implies that $F \mapsto L_{pe}(F)$ as defined above is the
left adjoint to $T \mapsto \mathcal{M}_{T}$.  \hfill $\Box$ \\

Note that for a general object $F\in D^{-}St(k)$,
the dg-category $L_{pe}(F)$ is not small in general, and
only belongs to $\mathbb{V}$. However, when $F$ is equivalent to
some small homotopy colimit of representable
$D^{-}$-stacks, then
$L_{pe}(F)$ becomes equivalent to a small homotopy limit
of small dg-categories, and thus is itself equivalent to a small
dg-category. This is in particular the case when
$F$ is $n$-geometric for some $n$.

\begin{df}\label{d9}
For a $D^{-}$-stack $F \in D^{-}St(k)$,
the dg-category $L_{pe}(F) \in Ho(dg-Cat_{\mathbb{V}})$ is called
\emph{the dg-category of perfect complexes on $F$}.
\end{df}

\subsection{Geometricity}

The purpose of this paragraph is to provide a proof
of the following theorem. Its proof will follow the main
two steps we have used in the proof the algebraicity of the 
1-stack $m_{C}$ in $\S 1$. We will first show that $\mathcal{M}_{\mathbf{1}}$
is locally geometric and of finite presentation. Then, we will 
construct a morphism $\mathcal{M}_{T}  \longrightarrow \mathcal{M}_{\mathbf{1}}$
and show that it is representable by identifying its fiber with 
dg-modules structures on a given perfect complex of $k$-modules.

\begin{thm}\label{t1}
Let $T$ be dg-category of finite type (see definition \ref{d3}).
Then, the $D^{-}$-stack $\mathcal{M}_{T}$ is
locally geometric and locally of finite presentation.
\end{thm}

\textit{Proof:}
The proof of the theorem will take us some time,
and will be devided in several propositions and lemmas.  \\

The following proposition is a particular case
where $T$ is the trivial dg-category $\mathbf{1}$, which 
is obviously of finite type. In this case, we note that 
$\mathcal{M}_{\mathbf{1}}$ is simply the \emph{$D^{-}$-stack of perfect
modules}. Indeed, for any $A\in sk-CAlg$, the simplicial set
$\mathcal{M}_{\mathbf{1}}(A)$ is the nerve of the category of
quasi-isomorphisms between cofibrant and perfect $N(A)$-dg-modules. 
In other words, $\mathcal{M}_{\mathbf{1}}(A)$ is a classifying space
of perfect $N(A)$-dg-modules up to quasi-isomorphisms. 
The underived version of $\mathcal{M}_{\mathbf{1}}$
has been considered in \cite[\S 21]{hs}, and the
following proposition is a derived analog of 
theorem \cite[Thm. 21.5]{hs}.

\begin{prop}\label{p4}
The $D^{-}$-stack $\mathcal{M}_{\mathbf{1}}$ is
locally geometric and locally of finite
presentation.
\end{prop}

\textit{Proof of proposition \ref{p4}:} The result of this proposition is somehow contained
in \cite{hagII}, but we will reproduce the proof here for the reader's
convenience.

Let $a,b\in \mathbb{Z}$ be two integers with
$a\leq b$.
We define a full sub-$D^{-}$-stack $\mathcal{M}_{\mathbf{1}}^{[a,b]}$
of $\mathcal{M}_{\mathbf{1}}$ in the following way.
For any $A\in sk-CAlg$, $\pi_{0}(\mathcal{M}_{\mathbf{1}}(A))$
can be identified with the quasi-isomorphism classes
of perfect $N(A)$-dg-modules.
We define $\mathcal{M}_{\mathbf{1}}^{[a,b]}(A)$ to be the
full sub-simplicial set of
$\mathcal{M}_{\mathbf{1}}(A)$ consisting of connected components  
corresponding to perfect $N(A)$-dg-modules
of Tor amplitude contained in $[a,b]$ (see definition \ref{d7}).
As the Tor amplitude is stable by base change (see 
proposition \ref{p2}),
this defines a full sub-$D^{-}$-stack
$\mathcal{M}_{\mathbf{1}}^{[a,b]} \subset \mathcal{M}_{\mathbf{1}}$.
Furthermore, as any perfect $N(A)$-dg-module has
a finite Tor amplitude (see proposition \ref{p2}),
we find
$$\mathcal{M}_{\mathbf{1}}=\cup_{a\leq b}\mathcal{M}_{\mathbf{1}}^{[a,b]}.$$
In order to prove Proposition \ref{p4} it is then enough to show that
each $\mathcal{M}_{\mathbf{1}}^{[a,b]}$ is
$n$-geometric and locally of finite presentation, for $n=b-a+1$. \\

We start by studying the diagonal.

\begin{lem}\label{l8}
Let $X$ be a representable $D^{-}$-stack, and
$$x : X \longrightarrow \mathcal{M}_{\mathbf{1}}^{[a,b]} \qquad
y : X \longrightarrow \mathcal{M}_{\mathbf{1}}^{[a,b]}$$
be two morphisms. Then, the $D^{-}$-stack
$X\times^{h}_{\mathcal{M}_{\mathbf{1}}^{[a,b]}}X$
is $(b-a)$-geometric and strongly of finite presentation.
\end{lem}

\textit{Proof of lemma \ref{l8}:} We start by
dealing directly with the case where $a=b$. In this case,
$\mathcal{M}_{\mathbf{1}}^{[a,a]}$ is equivalent to the $D^{-}$-stack
of vector bundles (denoted
by $\mathbf{Vect}=\coprod_{n}\mathbf{Vect}_{n}$ in \cite{hagII}). Therefore,
we know by \cite{hagII} that $\mathcal{M}_{\mathbf{1}}^{[a,a]}$
has a $(-1)$-representable diagonal of finite presentation.
Therefore, $X\times^{h}_{\mathcal{M}_{\mathbf{1}}^{[a,b]}}X$
is a $(-1)$-representable $D^{-}$-stack of finite presentation
and thus is $0$-geometric and strongly of finite presentation.

Let us now assume $a<b$.
Let $A$  be an object in $sk-CAlg$ representing $X$.
The morphisms $x$ and $y$ correspond to
$P$ and $Q$, two  perfect $N(A)$-dg-modules
of Tor amplitude contained in $[a,b]$.
By \cite[Appendix B]{hagII}, the
$D^{-}$-stack $X\times^{h}_{\mathcal{M}_{\mathbf{1}}^{[a,b]}}X \rightarrow X$
can be described by the functor
$$\begin{array}{ccc}
A/sk-CAlg & \longrightarrow & SSet \\
B & \mapsto & Map^{eq}_{N(A)-Mod}(P,Q\otimes^{\mathbb{L}}_{N(A)}N(B)).
\end{array}$$
Let us consider the $D^{-}$-stack over $X$
$$\begin{array}{cccc}
F : & A/sk-CAlg & \longrightarrow & SSet \\
& B & \mapsto & Map_{N(A)-Mod}(P,Q\otimes^{\mathbb{L}}_{N(A)}N(B)),
\end{array}$$
together with the natural monomorphism
$$j : X\times^{h}_{\mathcal{M}_{\mathbf{1}}^{[a,b]}}X \hookrightarrow
F.$$
We first claim that $j$ is a $0$-representable morphism. Indeed,
for $B\in A/sk-CAlg$, and $u \in F(B)$, corresponding to
a morphism $u : P\otimes^{\mathbb{L}}_{N(A)}N(B) \longrightarrow
Q\otimes^{\mathbb{L}}_{N(A)}N(B)$ of perfect $N(B)$-dg-modules,
the $D^{-}$-stack
$$j^{-1}(u):=(X\times^{h}_{\mathcal{M}_{\mathbf{1}}^{[a,b]}}X)\times^{h}_{F}
\mathbb{R}\underline{Spec}\, B$$
is the full sub-$D^{-}$-stack of $\mathbb{R}\underline{Spec}\, B$
where $u$ becomes a quasi-isomorphism. 
This is the sub-$D^{-}$-stack of $\mathbb{R}\underline{Spec}\, B$
where the cone of $u$ is quasi-isomorphic to zero, which by 
Prop. \ref{psupp} is 
a quasi-compact
Zariski open sub-$D^{-}$-scheme of $\mathbb{R}\underline{Spec}\, B$,
This implies that $j^{-1}(u)$ is a  $0$-geometric $D^{-}$-stack.

In order to prove the lemma, it only remains to show that
$F$ is an $(n-1)$-geometric $D^{-}$-stack. But $F$ can be written as
$$\begin{array}{cccc}
F : & A/sk-CAlg & \longrightarrow & SSet \\
& B & \mapsto & Map_{N(A)-Mod}(R,N(B)),
\end{array}$$
where $R$ is the $N(A)$-dg-module
$P\otimes_{N(A)}^{\mathbb{L}}Q^{\vee}$, of derived morphisms from $Q$ to $P$.
Note that $R$ is a perfect $N(A)$-dg-module, and
has Tor amplitude contained in $[a-b,b-a]$,
as by assumptions $P$ and $Q$ are both of Tor amplitude contained in $[a,b]$
(see proposition \ref{p2}).
Therefore the lemma \ref{l8} will be proved if we prove the
following general sub-lemma.

\begin{sublem}\label{sl1}
Let $A\in sk-CAlg$, and
$R$ be a $N(A)$-dg-module which is perfect and of Tor amplitude
contained in $[a,b]$, with $a,b \in \mathbb{Z}$. Then,
the $D^{-}$-stack over $\mathbb{R}\underline{Spec}\, A$
$$\begin{array}{cccc}
F : & A/sk-CAlg & \longrightarrow & SSet \\
& B & \mapsto & Map_{N(A)-Mod}(R,N(B)),
\end{array}$$
is $b$-geometric and strongly of finite presentation.
\end{sublem}

\textit{Proof of sub-lemma \ref{sl1}:}
The proof is by induction on $b$. Let us assume that
$b\leq 0$, then we have
$$F\simeq \mathbb{R}\underline{Spec}\, B$$
where $B$ is the object of $A/sk-CAlg$ defined
as the derived free commutative $A$-algebra
$$B:=\mathbb{R}Symm_{A}(D(R)),$$
where $D(R)$ is the simplicial $A$-module obtained by
denormalizing $R$. As $D(R)$ is
homotopy finitely presented as a simplicial $A$-module,
$B$ is homotopically finitely presented as
a commutative $A$-algbera. This implies that
$F$ is representable and strongly of finite presentation.

Let us now assume that $b>0$. By proposition \ref{p2}, 
one can find a homotopy cofibration
sequence of $N(A)$-dg-modules
$$R \longrightarrow R_{1} \longrightarrow P[-b+1],$$
where $P=N(A)^{r}$ is a free $N(A)$-dg-module,
and $R_{1}$ is a perfect $N(A)$-dg-module of
Tor amplitude contained in $[a,b-1]$. This gives rise to
a fibration sequence of $D^{-}$-stacks
$$K(\mathbb{G}_{a},b-1)^{r} \longrightarrow F_{1} \longrightarrow F,$$
where $F_{1}$ is defined as
$$\begin{array}{cccc}
F_{1} : & A/sk-CAlg & \longrightarrow & SSet \\
& B & \mapsto & Map_{N(A)-Mod}(R_{1},N(B)).
\end{array}$$
As the $D^{-}$-stack $K(\mathbb{G}_{a},b-1)^{r}$ is $r$-geometric, strongly
of finite presentation and smooth, 
the projection $F_{1} \longrightarrow F$ is a covering,
$(b-1)$-representable, strongly of finite presentation and smooth.
By induction $F_{1}$ is known to be $(b-1)$-geometric and
strongly of finite presentation. We deduce from this that
$F$ can be obtained as the classifying $D^{-}$-stack of
a $(b-1)$-smooth Segal groupoid $X_{*}$ such that each
$X_{i}$ is strongly of finite presentation. This implies
that $F$ is $b$-geometric and strongly of finite presentation
(see \cite{hagII}). \hfill $\Box$ \\

This finishes the proof of lemma \ref{l8}
\hfill $\Box$ \\

We come back to the proof of the proposition \ref{p4}. It
remains to show that the $D^{-}$-stack
$\mathcal{M}_{\mathbf{1}}^{[a,b]}$ has an $n$-atlas
(for $n=b-a+1$)
$$U \longrightarrow \mathcal{M}_{\mathbf{1}}^{[a,b]}$$
with $U$ locally of finite presentation.

For this, we proceed by induction on the amplitude $n$.
First of all, let us assume that $a=b$. Then,
$\mathcal{M}_{\mathbf{1}}^{[a,a]}$ is equivalent to
the $D^{-}$-stack of vector bundles. The fact that this defines
a $1$-geometric $D^{-}$-stack locally of finite presentation,
is well known (see \cite{hagII}).

Let us now assume that $n>1$. By induction
on $n$, we see
that $\mathcal{M}_{\mathbf{1}}^{[a,b-1]}$
is $(n-1)$-geometric locally of finite presentation.
We define
a $D^{-}$-stack $U$ in the following way. For
$A\in sk-CAlg$, we consider $N(A)-Mod^{(1)}$, the
model category of morphisms in $N(A)-Mod$, endowed with
its projective model structure (so that fibrations and
equivalences are defined on the underlying objects in $N(A)-Mod$, see
\cite[Thm. 11.6.1]{hi}).
We restrict to cofibrant objects $u : Q \rightarrow R$ in $N(A)-Mod^{(1)}$
such that $Q$ belongs to $\mathcal{M}_{\mathbf{1}}^{[a,b-1]}(A)$, and
$R$ belongs to $\mathcal{M}_{\mathbf{1}}^{[b-1,b-1]}(A)$. This defines
a full sub-category $\mathcal{C}(A)$ in
$(N(A)-Mod^{(1)})^{cof}$. Passing to the nerves of the
sub-categories of equivalences gives
a full sub-simplicial set $U(A):=N(w\mathcal{C}(A))$ of $N(w(N(A)-Mod^{(1)})^{cof})$.
When $A$ varies in $sk-CAlg$, this defines (up to the standard
strictification procedure) a simplicial presheaf
$$\begin{array}{ccc}
sk-CAlg & \longrightarrow & SSet \\
A & \mapsto & U(A).
\end{array}$$

\begin{lem}\label{l9}
The $D^{-}$-stack $U$ defined above is $(n-1)$-geometric and
locally of finite presentation.
\end{lem}

\textit{Proof of lemma \ref{l9}:} We consider the morphism
$$p : U \longrightarrow
\mathcal{M}_{\mathbf{1}}^{[a,b-1]}\times^{h} \mathbf{Vect}$$
sending a morphism $u : Q \rightarrow R$ to
the pair $(Q,R[b-1])$. We know by induction that
$\mathcal{M}_{\mathbf{1}}^{[a,b-1]}\times^{h} \mathbf{Vect}$
is $(n-1)$-geometric and locally of finite presentation.

\begin{sublem}\label{sl1'}
Let $M$  be the $C(k)$-enriched model category of $N(A)$-dg-modules, 
$M^{(1)}$ be the model category of morphisms in $M$, and $x$ and $y$ two objects in $M$. 
Then, there is a natural homotopy fiber square of simplicial sets
$$\xymatrix{
N(wM^{(1)}) \ar[r]^-{a} & N(wM)\times N(wM) \\
Map_{M}(x,y) \ar[r] \ar[u] & \bullet \ar[u]_-{x\times y}},$$
where $a$ is induced by the functor sending a morphism in $M$ to 
its domain and codomain. 
\end{sublem}

\textit{Proof:} This follows from \cite[Thm. 8.3]{re}.
\hfill $\Box$ \\

The above sub-lemma implies that 
for $A\in sk-CAlg$, and
$(Q,E)\in \mathcal{M}_{\mathbf{1}}^{[a,b-1]}(A)\times \mathbf{Vect}(A)$,
the $D^{-}$-stack
$$p^{-1}(Q,E):=F\times^{h}_{\mathcal{M}_{\mathbf{1}}^{[a,b-1]}\times^{h} 
\mathbf{Vect}}\mathbb{R}\underline{Spec}\, A$$
can be described as
$$\begin{array}{cccc}
p^{-1}(Q,E) : & A/sk-CAlg & \longrightarrow & SSet \\
& B & \mapsto & Map_{N(A)-Mod}(Q,E\otimes^{\mathbb{L}}_{N(A)}N(B)[-b+1]).
\end{array}$$
The $N(A)$-dg-module $R$ being of Tor amplitude contained
in $[a,b-1]$, 
$R\otimes^{\mathbb{L}}_{N(A)}E^{\vee}[b-1]$
is perfect and of Tor amplitude contained
in $[a,0]$. We deduce from
Sub-lemma \ref{sl1} that $p^{-1}(Q,R)$ is a representable $D^{-}$-stack
locally of finite presentation. In other words, the morphism
$p$ is representable and locally of finite presentation.
This implies that the
$D^{-}$-stack $U$ is itself $(n-1)$-representable and locally
of finite presentation.
\hfill $\Box$ \\

We now consider for any $A\in sk-CAlg$, the morphism
$$U(A) \longrightarrow \mathcal{M}_{\mathbf{1}}(A)$$
sending a morphism of $N(A)$-dg-module
$u : Q \longrightarrow R$ to its homotopy fiber
$Hofib(u)$. By the definition of $U(A)$, it is
clear that $Hofib(u)$ belongs to
$\mathcal{M}_{\mathbf{1}}^{[a,b]}(A) \subset
\mathcal{M}_{\mathbf{1}}(A)$. When $A$ varies
in $sk-CAlg$, and using a model for the homotopy fiber
which is compatible with base change,
one gets a morphism of $D^{-}$-stacks
$$\pi : U \longrightarrow \mathcal{M}_{\mathbf{1}}^{[a,b]}.$$

\begin{lem}\label{l10}
The morphism $\pi$ is an $(n-1)$-representable smooth covering.
\end{lem}

\textit{Proof of lemma \ref{l10}:} We already know that $U$ is
$(n-1)$-geometric locally of finite presentation, and that the diagonal of
$\mathcal{M}_{\mathbf{1}}^{[a,b]}$ is $(n-1)$-representable
and locally of finite presentation.
Therefore, the morphism $\pi$ is itself $(n-1)$-representable
and locally of finite presentation.

Let us first prove that $\pi$ is a covering. Let $A\in sk-CAlg$, and
$P\in \mathcal{M}_{\mathbf{1}}^{[a,b]}(A)$. By definition,
one can find a vector bundle $E$ on
$\mathbb{R}\underline{Spec}\, A$, and
a morphism $E[-b] \longrightarrow P$, such that its
homotopy cofiber $Q$ is of Tor amplitude contained in
$[a,b-1]$. Therefore, there exists a homotopy fibration sequence
$$\xymatrix{P \ar[r] & Q \ar[r] & E[-b+1].}$$
This implies by definition of $U$ and $\pi$ that $P$ is
in the image of $\pi$. This shows that $\pi$ is a covering.

It remains to show that $\pi$ is smooth. For this, let
$A\in sk-CAlg$, and $P$ be a point in
$\mathcal{M}_{\mathbf{1}}^{[a,b]}$, corresponding to a morphism
$$P : X:=\mathbb{R}\underline{Spec}\, A \longrightarrow \mathcal{M}_{\mathbf{1}}^{[a,b]}.$$
Let us denote by $\pi^{-1}(P)$ the $D^{-}$-stack over $X$ defined
by
$$\pi^{-1}(P):=U\times^{h}_{\mathcal{M}_{\mathbf{1}}^{[a,b]}}X
\longrightarrow X.$$
The diagonal of $\mathcal{M}_{\mathbf{1}}^{[a,b]}$ being
$(n-1)$-representable and locally of finite presentation,
the $D^{-}$-stack $\pi^{-1}(P)$ is itself $(n-1)$-geometric
and locally of finite presentation.
We need to show that $\pi^{-1}(P)$ is smooth over $X$.
We consider the morphism
$$p : U \longrightarrow
\mathcal{M}_{\mathbf{1}}^{[a,b-1]}\times^{h} \mathbf{Vect}$$
as well as the induced morphism
$$q : \pi^{-1}(P) \longrightarrow U\times^{h} X \longrightarrow
\mathcal{M}_{\mathbf{1}}^{[a,b-1]}\times^{h} \mathbf{Vect} \times^{h}X
\longrightarrow \mathbf{Vect}\times^{h} X.$$
As we know that the family of morphisms
$$u_{r} : X \longrightarrow \mathbf{Vect}\times^{h} X,$$
corresponding to the trivial rank $r$ vector bundle on $X$
(for various $r$), form
a smooth atlas,
we see that it is enough to check that
the morphism
$$V:=\pi^{-1}(P)\times_{\mathbf{Vect}\times^{h} X}^{h}X \longrightarrow
X$$
is smooth, for each of the morphism $u_{r}$.
Let the rank $r$ be fixed, then
by construction, the $D^{-}$-stack
$V$ can be described in the following way
$$\begin{array}{cccc}
V : & A/sk-CAlg & \longrightarrow & SSet \\
 & (A\rightarrow B) & \mapsto & Map_{N(B)-Mod}(B^{r}[-b+1],P\otimes^{\mathbb{L}}_{N(A)}N(B)[1]).
\end{array}$$
In other words, $V$ can be written as
$$\begin{array}{cccc}
V : & A/sk-CAlg & \longrightarrow & SSet \\
 & (A\rightarrow B) & \mapsto & Map_{N(A)-Mod}(Q,N(B)).
\end{array}$$
where $Q:=A^{r}\otimes^{\mathbb{L}}_{N(A)}P^{\vee}[-b]$.
Now, let $u : Q \longrightarrow N(B)$ be a morphism
of $N(A)$-dg-modules, corresponding to
commutative diagram of of $D^{-}$-stacks
$$\xymatrix{
& V \ar[d] \\
Y \ar[r] \ar[ru]^-{u} & X}$$
where $Y:=\mathbb{R}\underline{Spec}\, B$. For any
simplicial $B$-module $M$, we have
$$\mathbb{D}er_{u}(V/X,M)\simeq Map_{N(A)-Mod}(Q,N(M)),$$
where $N(M)$ is the normalized $N(A)$-dg-module associated
to $M$. This shows that
the relative cotangent complex of $V$ over $X$ at the point $u$, given by
$$\mathbb{L}_{V/X,u}\simeq Q\otimes^{\mathbb{L}}_{N(A)}N(B).$$
As $Q=A^{r}\otimes^{\mathbb{L}}_{N(A)}P^{\vee}[-b]$
is of Tor amplitude contained in $[0,b-a]$ (see proposition \ref{p2}),
the infinitesimal criterion for smoothness (see \cite[\S 2.2.5]{hagII}) tells us that
$V \rightarrow X$ is smooth.
\hfill $\Box$ \\

We have finished the proof of lemma \ref{l10} and thus of proposition \ref{p4}. \hfill $\Box$ \\

We come back to the proof of theorem \ref{t1}.
We start by recalling
from lemma \ref{l7} that the Yoneda embedding
$\underline{h} : T \longrightarrow \widehat{T}_{pe}$
induces an isomorphism
in $D^{-}St(k)$
$$\underline{h}^{*} : \mathcal{M}_{\widehat{T}_{pe}} \simeq
\mathcal{M}_{T}.$$
In particular, we can assume that
$T$ is triangulated (see \ref{l1} $(1)$).
We then chose a perfect generator
$E$ of $\widehat{T}$, which as
$T$ is triangulated is considered as
an object in $T$. We consider the
morphism of dg-categories
$$i_{E} : \mathbf{1} \longrightarrow T$$
which classifies the object $E$.
This morphism
induces a morphism on the correspoding $D^{-}$-stacks
$$\pi:=i^{*}_{E} : \mathcal{M}_{T} \longrightarrow \mathcal{M}_{\mathbf{1}}.$$
By our proposition \ref{p4} $\mathcal{M}_{\mathbf{1}}$ is a locally
geometric $D^{-}$-stack, locally of finite presentation.
It can be written a union
$$\mathcal{M}_{\mathbf{1}}=\cup_{a\leq b}\mathcal{M}_{\mathbf{1}}^{[a,b]}$$
where $\mathcal{M}_{\mathbf{1}}^{[a,b]}$ denotes the full sub-$D^{-}$-stack
of objects whose Tor amplitude is contained in $[a,b]$.
We denote by $\mathcal{M}_{T}^{[a,b]}$ the full sub-$D^{-}$-stack
of $\mathcal{M}_{T}$ defined by the following homotopy cartesian square
$$\xymatrix{
\mathcal{M}_{T} \ar[r]^-{\pi} & \mathcal{M}_{\mathbf{1}} \\
\mathcal{M}_{T}^{[a,b]} \ar[u] \ar[r] & \mathcal{M}_{\mathbf{1}}^{[a,b]} \ar[u].}$$
Clearly, 
$$\mathcal{M}_{T}=\cup_{a\leq b}\mathcal{M}_{T}^{[a,b]}.$$
Therefore, in order to prove that $\mathcal{M}_{T}$ is locally
geometric and locally of finite presentation, it
is enough to prove that each $D^{-}$-stack $\mathcal{M}_{T}^{[a,b]}$
is $n$-geometric and locally of finite presentation, for
some integer $n$ (depending on $a$ and $b$).

\begin{prop}\label{p5}
Let $a\leq b$. There exists an integer $n$ such that
the morphism
$$\pi : \mathcal{M}_{T}^{[a,b]} \longrightarrow \mathcal{M}_{\mathbf{1}}^{[a,b]}$$
is $n$-representable and strongly of finite presentation.
\end{prop}

\textit{Proof of proposition \ref{p5}:} As $T$ is of finite
type, we write $T$ as
$\widehat{B^{op}}_{pe}$, for the dg-algebra $B=T(E,E)$ of endormophisms of the
generator $E$. From our hypothesis on $T$ and corollary \ref{ccc}, $B$
is homotopically finitely presented in the model
category $dg-Alg$. 

We let $A\in sk-CAlg$, $X:=\mathbb{R}\underline{Spec}\, A$, and
$X \longrightarrow \mathcal{M}_{\mathbf{1}}^{[a,b]}$ be the morphism
corresponding to a perfect $N(A)$-dg-module $P$ of Tor amplitude contained
in $[a,b]$. We denote by $\mathcal{E}(P)$ the $N(A)$-dg-algebra
of derived endomorphisms of $P$

$$\mathcal{E}(P):=\mathbb{R}\underline{Hom}_{N(A)}(P,P).$$
We let $\pi^{-1}(P)$ be the $D^{-}$-stack over $X$ defined by
$$\pi^{-1}(P):=\mathcal{M}_{T}^{[a,b]} \times^{h}_{\mathcal{M}_{\mathbf{1}}^{[a,b]}}X.$$
We define a $D^{-}$-stack over $X$ in the following way
$$\begin{array}{cccc}
\underline{Map}(B,\mathcal{E}(P)) : & A/sk-CAlg & \longrightarrow & SSet \\
 & (A\rightarrow A') & \mapsto & Map_{dg-Alg}(B,\mathcal{E}(P)\otimes_{N(A)}^{\mathbb{L}}N(A')). 
\end{array}$$
Here $\mathcal{E}(P)\otimes_{N(A)}^{\mathbb{L}}N(A')$
is a $N(A')$-dg-algebra obtained by base change the $N(A)$-dg-algebra
$\mathcal{E}(P)$, and considering it as a dg-algebra over $k$ by the forgetful functor.

\begin{lem}\label{l11}
There exists an isomorphism in $D^{-}St(X)$
$$\pi^{-1}(P)\simeq \underline{Map}(B,\mathcal{E}(P)).$$
\end{lem}

\textit{Proof of lemma \ref{l11}:} As explained
in \S 2.2, we will consider dg-algebras as
dg-categories with a unique object. Then, because of lemma \ref{l7} 
$\mathcal{M}_{T}\simeq \mathcal{M}_{B}$
and the projection $\pi : \mathcal{M}_{B} \longrightarrow \mathcal{M}_{\mathbf{1}}$
becomes induced by the unit morphism $k \longrightarrow B$. The dg-algebra
$\mathcal{E}(P)\otimes_{N(A)}^{\mathbb{L}}N(A')$, considered as a dg-category, can be identified with the full
sub-dg-category of $\widehat{A'}_{pe}$ with $P\otimes_{N(A)}^{\mathbb{L}}N(A')$
a unique object. Therefore, for
any $A' \in A/sk-CAlg$, we find a diagram with homotopy cartesian squares
$$\xymatrix{
Map_{dg-Cat}(B,\widehat{A'}_{pe}) \ar[r] & Map_{dg-Cat}(k,\widehat{A'}_{pe}) \\
Map_{dg-Cat}(B,\mathcal{E}(P)\otimes_{N(A)}^{\mathbb{L}}N(A')) \ar[r] \ar[u] &
Map_{dg-Cat}(k,\mathcal{E}(P)\otimes_{N(A)}^{\mathbb{L}}N(A')) \ar[u] \\
 \pi^{-1}(P)(A') \ar[u]  \ar[r] & \bullet. \ar[u] }$$   

Let 
$$\xymatrix{dg-Alg \ar@<2pt>[r]^G & dg-Cat_* \ar@<2pt>[l]^E}$$
be the Quillen adjunction between the model category of dg-algebras
and the model category of pointed dg-categories. 
It is easy to see that the left derived functor
$$\mathbb{L}G : Ho(dg-Alg) \longrightarrow Ho(dg-Cat)$$
is fully faithful. Therefore, the functor $G$ also induces natural
equivalences
$$Map_{dg-Alg}(B_{1},B_{2}) \simeq Map_{dg-Cat_{*}}(G(B_{1}),G(B_{2})).$$
The diagram with homotopy cartesian squares above tells us that
$\pi^{-1}(P)(A')$ can be naturally identified with
$Map_{dg-Cat_{*}}(B,\mathcal{E}(P)\otimes_{N(A)}^{\mathbb{L}}N(A'))$, and therefore
with $Map_{dg-Alg}(B,\mathcal{E}(P)\otimes_{N(A)}^{\mathbb{L}}N(A'))$. When $A'$ varies
in $A/sk-CAlg$, this provides an equivalence between
$\pi^{-1}(P)$ and $\underline{Map}(B,\mathcal{E}(P))$. \hfill $\Box$ \\

We come back to the proof of proposition \ref{p5}. Lemma \ref{l11} tells us that
it is enough to show that the $D^{-}$-stack
$\underline{Map}(B,\mathcal{E}(P))$ is $n$-geometric for some $n$, and
that the projection $\underline{Map}(B,\mathcal{E}(P)) \longrightarrow X$
is strongly of finite presentation.

Recall that the dg-algebra $B$ is homotopically finitely presented    
in the model category $dg-Alg$ of dg-algebras. 
Therefore, using proposition \ref{p1} $B$ is (equivalent to) a retract
of some dg-algebra $B'$, for which there exists a finite diagram
$$\xymatrix{
B_{0}=k \ar[r] & B_{1} \ar[r] & \dots & B_{i} 
\ar[r] & B_{i+1} \ar[r] & \dots \ar[r] & B_{m}=B',}$$
and such that for any $i$ there exists a homotopy push-out diagram
$$\xymatrix{
B_{i} \ar[r] & B_{i+1} \\
C_{i} \ar[r] \ar[u] & k,\ar[u]}$$
where $C_{i}$ is the free dg-algebra over the complex
$k[p_{i}]$ for some integer $p_{i}$. Clearly, the morphism of $D^{-}$-stacks
$$\underline{Map}(B,\mathcal{E}(P)) \longrightarrow X$$
is retract of
$$\underline{Map}(B',\mathcal{E}(P)) \longrightarrow X.$$
As $n$-geometric $D^{-}$-stacks strongly of finite presentation are
stable by retracts (see lemma \ref{l4}), it is enough to show that
$\underline{Map}(B',\mathcal{E}(P)) \longrightarrow X$ is $n$-representable
for some $n$ and strongly of finite presentation. But, there exists a diagram
of $D^{-}$-stacks over $X$
$$\xymatrix{\underline{Map}(B',\mathcal{E}(P)) \ar[r] & \dots \ar[r] &
\underline{Map}(B_{i+1},\mathcal{E}(P)) \ar[r] & \underline{Map}(B_{i},\mathcal{E}(P)) \ar[r] &
\dots  & \underline{Map}(k,\mathcal{E}(P))=X,}$$
as well as homotopy pull-back squares of $D^{-}$-stacks over $X$
$$\xymatrix{
\underline{Map}(B_{i+1},\mathcal{E}(P)) \ar[r] \ar[d]& \underline{Map}(B_{i},\mathcal{E}(P)) \ar[d] \\
\bullet \ar[r] & \underline{Map}(C_{i},\mathcal{E}(P)).}$$
As $n$-geometric $D^{-}$-stacks strongly of finite presentations
are stable by homotopy pull-backs (see lemma \ref{l4}), we see that it is enough to
show that $\underline{Map}(C_{i},\mathcal{E}(P)) \longrightarrow X$ is
$n$-geometric and strongly of finite presentation for some $n$.
But, as $C_{i}$ is free over $k[p_{i}]$, $\underline{Map}(C_{i},\mathcal{E}(P))$
can be described, up to an equivalence, in the following terms
$$\begin{array}{cccc}
\underline{Map}(C_{i},\mathcal{E}(P)) : & A/sk-CAlg & \longrightarrow & SSet \\
 & (A\rightarrow A') & \mapsto & Map_{C(k)}(k[p_{i}],\mathcal{E}(P)\otimes^{\mathbb{L}}_{N(A)}N(A')).
\end{array}$$
In other words, if we let $K$ be the $N(A)$-dg-module
$\mathcal{E}(P)^{\vee}[p_{i}]$, one finds
$$\begin{array}{cccc}
\underline{Map}(C_{i},\mathcal{E}(P)) : & A/sk-CAlg & \longrightarrow & SSet \\
 & (A\rightarrow A') & \mapsto & Map_{C(k)}(K,N(A')).
\end{array}$$
Now, $K$ is of Tor amplitude contained in $[a-b-p_{i},b-a-p_{i}]$, and thus
sub-lemma \ref{sl1} tells us that
$$\underline{Map}(C_{i},\mathcal{E}(P)) \longrightarrow X$$
is $(b-a-p_{i})$-representable and strongly of finite presentation.
This finishes the proof of the proposition \ref{p5}, 
with the additional information
that $n$ can be taken to be $(b-a-p_{i_{0}})$, where $p_{i_{0}}$ is the lowest of the
integers $p_{i}$. \\

Propositions \ref{p4} and \ref{p5}, together with lemma \ref{lrep}
achieve the proof of theorem \ref{t1}. \hfill $\Box$ \\

\begin{cor}\label{c3-}
Let $T$ be a smooth and proper dg-category 
then $\mathcal{M}_{T}$ is a locally geometric
$D^{-}$-stack locally of finite presentation.
\end{cor}

\textit{Proof:} This follows from
corollary \ref{c2} and theorem \ref{t1}. \hfill $\Box$ \\

The following result is a corollary of the proof 
of theorem \ref{t1} (more precisely
of proposition \ref{p5}).

\begin{cor}\label{c3}
Let $T$ be a triangulated dg-category of finite type, and $E\in T$ a
perfect generator classified by
a morphism $i_{E} : \mathbf{1} \longrightarrow T$, and let
$$\pi:=i^{*}_{E} : \mathcal{M}_{T} \longrightarrow
\mathcal{M}_{\mathbf{1}}$$
be the induced morphism. Then, for any $A\in sk-CAlg$,
any morphism
$$x : X:=\mathbb{R}\underline{Spec}\, A \longrightarrow \mathcal{M}_{\mathbf{1}},$$
the induced morphism
$$\mathcal{M}_{T}\times^{h}_{\mathcal{M}_{\mathbf{1}}}X \longrightarrow X$$
is $n$-representable and strongly of finite presentation, for
some $n$ depending on the choice of $X$ and the morphism
$x$.
\end{cor}

\textit{Proof:} This follows from proposition \ref{p5}, as
any morphism $X\longrightarrow \mathcal{M}_{\mathbf{1}}$ factors through
one of the $\mathcal{M}_{\mathbf{1}}^{[a,b]}$ for some $a\leq b$. \hfill $\Box$ \\

The above corollary almost says that
$\pi$ is $n$-representable and strongly of finite presentation,
except that the integer $n$ depends on the morphism
$X\longrightarrow \mathcal{M}_{\mathbf{1}}$, and can not
be chosen uniformly. \\

Another important corollary of theorem \ref{t1} is the following
description of the tangent complexes of the $D^{-}$-stack
$\mathcal{M}_{T}$. Recall from \cite[2.2]{hagII} that $n$-geometric
$D^{-}$-stacks have an obtruction theory, and in particular
a (co)tangent complex. It is immediate to deduce that
any locally geometric $D^{-}$-stack also has  
an obstruction theory and a global (co)tangent complex.

\begin{cor}\label{c4}
Let $T$ be a triangulated dg-category of finite type, and
$$E : Spec\, k \longrightarrow \mathcal{M}_{T}$$
be a morphism corresponding to an
object in $[T]$. Then, the tangent complex of
$\mathcal{M}_{T}$ at the point $E$ is given by
$$\mathbb{T}_{\mathcal{M}_{T},E}\simeq T(E,E)[1].$$
\end{cor}

\textit{Proof:} Let
$$E : *=Spec\, k \longrightarrow \mathcal{M}_{T}$$
as in the statement. We consider
$$\Omega_{E}\mathcal{M}_{T}:=*\times^{h}_{\mathcal{M}_{T}}*,$$
the loop $D^{-}$-stack taken at $E$. Using \cite[Appendix B]{hagII},
and \cite{to},
one has the following description
$$\begin{array}{cccc}
\Omega_{E}\mathcal{M}_{T} : & sk-CAlg & \longrightarrow & SSet \\
 & A & \mapsto & Map'(k,T(E,E)\otimes_{k}^{\mathbb{L}}N(A)),
\end{array}$$
where the $Map'$ is the full sub-simplicial set
of $Map_{C(k)}(k,T(E,E)\otimes_{k}^{\mathbb{L}}N(A))$
corresponding to isomorphisms in 
$$\pi_{0}(Map_{C(k)}(k,T(E,E)\otimes_{k}^{\mathbb{L}}N(A)))\simeq
[\widehat{T\otimes^{\mathbb{L}}N(A)}](E\otimes^{\mathbb{L}}N(A),
E\otimes^{\mathbb{L}}N(A))).$$
Therefore, $\Omega_{E}\mathcal{M}_{T}$ is a Zariski open sub-$D^{-}$-stack
of
$$\begin{array}{cccc}
\underline{End}(E) : & sk-CAlg & \longrightarrow & SSet \\
 & A & \mapsto & Map(k,T(E,E)\otimes_{k}^{\mathbb{L}}N(A)).
\end{array}$$
Using the definition of derivations as in \cite[\S 1.2]{hagII},
it is easy to check that
$$\mathbb{T}_{\underline{End}(E),id}
\simeq T(E,E),$$
and so
$$\mathbb{T}_{\Omega_{E}\mathcal{M}_{T},*}\simeq
T(E,E).$$
But, we have (see \cite{hagII})
$$\mathbb{T}_{\Omega_{E}\mathcal{M}_{T},*}\simeq
\mathbb{T}_{\mathcal{M}_{T},E}[-1],$$
and the required formula follows.  \hfill $\Box$ \\

In order to finish this paragraph we mention 
the following easy, but interesting, consequence
of theorem \ref{t1}. We let $T$ be a dg-category, and
$I(1)$ be the free
dg-category with two objects $0$ and $1$ and
$I(1)(0,1)=k$ (i.e. $I(1)$ is the free dg-category 
of the category $\Delta(1)$, having two objects
and a unique morphisms between them). We set
$T_{(1)}:=T\otimes I(1)$, and we consider 
the two morphisms
$$s : T \longrightarrow T_{(1)}  \qquad 
t : T \longrightarrow T_{(1)}$$
induced by the two objects $0 : \mathbf{1} \rightarrow I(1)$ and
$1 : \mathbf{1} \rightarrow I(1)$. Passing to the associated $D^{-}$-stacks, 
one gets two morphisms
$$s : \mathcal{M}_{T_{(1)}} \longrightarrow \mathcal{M}_{T}
\qquad t : \mathcal{M}_{T_{(1)}} \longrightarrow \mathcal{M}_{T}.$$

\begin{df}\label{d10-}
The $D^{-}$-stack $\mathcal{M}_{T_{(1)}}$ 
is called the $D^{-}$-stack of \emph{morphisms
between pseudo-perfect dg-modules}. It is denoted by 
$\mathcal{M}_{T}(1)$. The two morphisms
$$s : \mathcal{M}_{T}(1) \longrightarrow \mathcal{M}_{T}
\qquad t : \mathcal{M}_{T}(1) \longrightarrow \mathcal{M}_{T}$$
are called \emph{source and target}.
\end{df}

It is easy to check that if $T$ is of finite type then 
so is $T_{(1)}$. Therefore, 
from theorem \ref{t1} we see that 
if $T$ is a dg-category of finite type then 
$\mathcal{M}_{T}(1)$ is a locally geometric $D^{-}$-stack.
Furthermore, sub-lemma \ref{sl1} easily implies the
following corollary.

\begin{cor}\label{cd10-}
Let $T$ be a  dg-category of finite type. 
Let us consider the morphism
$$s\times t : \mathcal{M}_{T}(1) \longrightarrow 
\mathcal{M}_{T}\times^{h} \mathcal{M}_{T}.$$
Then, for any $A\in sk-CAlg$,
any morphism
$$x : X:=\mathbb{R}\underline{Spec}\, A \longrightarrow 
\mathcal{M}_{T}\times^{h} \mathcal{M}_{T},$$
the induced morphism
$$\mathcal{M}_{T}(1)\times^{h}_{\mathcal{M}_{T}\times^{h} \mathcal{M}_{T}}X \longrightarrow X$$
is $n$-representable and strongly of finite presentation, for
some $n$ depending on the choice of $X$ and the morphism
$x$.
\end{cor}

\textit{Proof:} Indeed, if the morphism $x$ corresponds to 
two $T^{op}\otimes^{\mathbb{L}}N(A)$-dg-modules
$P$ and $Q$, then the $D^{-}$-stack
$\mathcal{M}_{T}(1)\times^{h}_{\mathcal{M}_{T}\times^{h} \mathcal{M}_{T}}X$
is given by 
$$\begin{array}{ccc}
A/sk-CAlg & \longrightarrow & SSet \\
(A\rightarrow A') & \mapsto & Map_{N(A)-Mod}(K,N(A')),
\end{array}$$
where $K:=P\otimes^{\mathbb{L}}_{N(A)}Q^{\vee}$ is the
$N(A)$-dg-module of derived morphisms from $P$ to $Q$ (this follows
from sub-lemma \ref{sl1'}). 
Sub-lemma \ref{sl1} implies the result. \hfill $\Box$ \\

\subsection{Sub-stacks strongly of finite presentation}

This section describes an exhaustive family of open sub-$D^{-}$-stacks
$\mathcal{M}_{T}^{\nu} \subset \mathcal{M}_{T}$ which are
strongly of finite presentation. This description is
not canonical, as it depends on the choice of
a perfect generator in $\widehat{T}$, but it will be useful
in order to prove that certain $D^{-}$-stacks are
strongly of finite presentation.  \\

We let $\nu : \mathbb{Z} \longrightarrow \mathbb{N}$
be a function with finite support (i.e. there is only
a finite number of $i\in \mathbb{Z}$ with $\nu(i)\neq 0$).
For $A\in sk-CAlg$, we define
a full sub-simplicial set $\mathcal{M}_{\mathbf{1}}^{\nu}(A) \subset
\mathcal{M}_{\mathbf{1}}(A)$ in the following way. A
point $P\in \mathcal{M}_{\mathbf{1}}(A)$ belongs to
$\mathcal{M}_{\mathbf{1}}^{\nu}(A)$ if for any
field $K$ and any morphism $A \rightarrow \pi_{0}(A) \rightarrow K$
in $sk-CAlg$, we have
$$Dim_{K}H^{i}(P\otimes^{\mathbb{L}}_{N(A)}K)\leq \nu(i) \; \forall \; i.$$
This condition is stable by base change, and thus defines
a full sub-$D^{-}$-stack
$$\mathcal{M}_{\mathbf{1}}^{\nu} \subset \mathcal{M}_{\mathbf{1}}.$$
We start by noticing that
if $a\leq b$ are such that $\nu(i)=0$ for all
$i\notin [a,b]$, then
$$\mathcal{M}_{\mathbf{1}}^{\nu} \subset \mathcal{M}_{\mathbf{1}}^{[a,b]}.$$
Indeed, if $P\in \mathcal{M}_{\mathbf{1}}^{\nu}(A)$, then
its base change $P':=P\otimes^{\mathbb{L}}_{N(A)}\pi_{0}(A)$
is a perfect complex of $\pi_{0}(A)$-modules whose
fibers at all points in $Spec\, \pi_{0}(A)$ are of Tor amplitude
contained in $[a,b]$. This implies that $P'$ has
Tor amplitude contained in $[a,b]$, and  thus that
$P$ itself has Tor amplitude contained in $[a,b]$ (see proposition \ref{p2}).

Let $T$ be a triangulated dg-category of finite type, and let
us chose a perfect generator $E\in T$. It is
classified by a morphism $\mathbf{1} \longrightarrow T$, which
provides a morphism of $D^{-}$-stacks
$$\pi : \mathcal{M}_{T} \longrightarrow \mathcal{M}_{\mathbf{1}}.$$
We will denote by $\mathcal{M}_{T}^{\nu}$ the full sub-$D^{-}$-stack
of $\mathcal{M}_{T}$ defined by the homotopy pull-back square
$$\xymatrix{
\mathcal{M}_{T} \ar[r] & \mathcal{M}_{\mathbf{1}} \\
\mathcal{M}_{T}^{\nu} \ar[u] \ar[r] & \mathcal{M}_{\mathbf{1}}^{\nu} \ar[u].}$$
Note that $\mathcal{M}_{T}^{\nu}$ depends on
the parameter $\nu$, but also of the choice of the generator
$E$.

\begin{prop}\label{p6}
Let $T$ be a triangulated dg-category of finite type, 
and $E$ and $\nu$ be as above.
The $D^{-}$-stack $\mathcal{M}_{T}^{\nu}$ is
$n$-geometric for some $n$, and strongly
of finite presentation.
\end{prop}

\textit{Proof:} There exists a diagram with homotopy cartesian squares
$$\xymatrix{
\mathcal{M}_{T} \ar[r]^-{\pi} & \mathcal{M}_{\mathbf{1}} \\
\mathcal{M}_{T}^{[a,b]} \ar[r]^-{\pi^{[a,b]}} \ar[u] & \mathcal{M}_{\mathbf{1}}^{[a,b]} \ar[u]\\
\mathcal{M}_{T}^{\nu} \ar[r] \ar[u] & \mathcal{M}_{\mathbf{1}}^{\nu}, \ar[u]}$$
for some $a\leq b$ with $\nu(i)=0$ for all $i\notin [a,b]$. As we already
know by \ref{p5} that $\pi^{[a,b]}$ is $n$-representable for some
$n$ and strongly of finite presentation, then
it is enough to show that
the $D^{-}$-stack $\mathcal{M}_{\mathbf{1}}^{\nu}$, which
is known to be $(b-a+1)$-geometric, is strongly of finite
presentation. By lemma \ref{l8} we already know that the diagonal of
$\mathcal{M}_{\mathbf{1}}^{\nu}$ is strongly of finite presentation, and thus
it only remains to show that $\mathcal{M}_{\mathbf{1}}^{\nu}$
is quasi-compact. Equivalentely, it is enough to find
a quasi-compact and $n$-geometric $D^{-}$-stack
$X$ and a covering $X \longrightarrow
\mathcal{M}_{\mathbf{1}}^{\nu}$.

We proceed by induction on the lenght $l$ of $\nu$, where
$l$ is the number of non-zero terms in the sequence $\nu$.
For $l=1$, and say $\nu(i)\neq 0$,
$\mathcal{M}_{\mathbf{1}}^{\nu}$ is equivalent to the
$D^{-}$-stack of vector bundles of rank less than $\nu(i)$.
This is quasi-compact as it is covered by $\nu(i)$ copies of the
final object $*$. Let us assume that all the $D^{-}$-stacks
$\mathcal{M}_{\mathbf{1}}^{\nu'}$ are quasi-compact for all
$\nu'$ of lenght less than $l-1$, and let $\nu$ be of lenght
$l$.

We define $\nu'$ by
$$\nu'(i):=\nu(i) \; \forall \; i<b-1 \qquad
\nu'(b-1)=\nu(b-1)+\nu(b) \qquad \nu'(i)=0 \; \forall i>b-1,$$
where $b$ is the largest integer such that $\nu(b)\neq 0$.
We define a $D^{-}$-stack $X$ in the following way.
We first consider the $D^{-}$-stack of morphisms
$\mathcal{M}_{\mathbf{1}}(1)$ of definition \ref{d10-}. We consider the full
sub-$D^{-}$-stack $X$ of $\mathcal{M}_{\mathbf{1}}(1)$, whose
objects over $A\in sk-CAlg$ consists of all
objects $Q \rightarrow R$ in
$\mathcal{M}_{\mathbf{1}}(1)(A)$, with
$Q\in \mathcal{M}_{\mathbf{1}}^{\nu'}$ and $R$ of the form
$E[-b+1]$, for $E$ a vector bundle
of rank $\nu(b)$. We construct a morphism of $D^{-}$-stacks
$$X \longrightarrow \mathcal{M}_{\mathbf{1}}^{[a,b]}$$
by sending a morphism $Q\rightarrow R$ as above to its
homotopy fiber. By \ref{sl1}, it is easy to see that
the natural morphism sending $Q\rightarrow R$ to the pair
$(Q,R[b-1])$
$$X \longrightarrow \mathcal{M}_{\mathbf{1}}^{\nu'}\times^{h}
\mathbf{Vect}_{\nu(b)}$$
is $n$-representable, for some $n$, and strongly of finite presentation.
By the induction assumption, this implies that $X$ is itself
$n$-geometric for some $n$, and strongly of finite presentation.

Let us now consider the homotopy pull-back square
$$\xymatrix{
\mathcal{M}^{\nu}_{\mathbf{1}} \ar[r]^-{j} & \mathcal{M}^{[a,b]}_{\mathbf{1}} \\
U \ar[r] \ar[u] & X. \ar[u]}$$
We first note that the morphism $j$ is a quasi-compact Zariski open
immersion. Indeed, for $A\in sk-CAlg$, and
$P$ a perfect $N(A)$-dg-module of Tor amplitude
contained in $[a,b]$, the full sub-$D^{-}$-stack
of $\mathbb{R}\underline{Spec}\, A$ where $P$ lies in
$\mathcal{M}_{\mathbf{1}}^{\nu}(A)$ is a Zariski open
sub-$D^{-}$-scheme, corresponding to open sub-scheme
of $Spec\, \pi_{0}(A)$ where $P\otimes^{\mathbb{L}}_{N(A)}\pi_{0}(A)$
lies in $\mathcal{M}_{\mathbf{1}}^{\nu}(\pi_{0}(A))$. But, for
a given perfect complex $Q$ of $k'$-modules, for some
commutative ring $k'$, the sub-set of
points $x\in Spec\, k'$ such that $Dim_{k'(x)}H^{i}(Q\otimes_{k}^{\mathbb{L}}k'(x))\leq \nu(i)$
for all $i$, is a quasi-compact open sub-set by semi-continuity. 
Therefore, the $D^{-}$-stack $U$ is also $n$-geometric and
strongly quasi-compact. It remains to show that
the morphism
$$U \longrightarrow \mathcal{M}^{\nu}_{\mathbf{1}}$$
is a covering, or equivalently that
any $P\in \mathcal{M}_{\mathbf{1}}^{\nu}(A)$ is, locally on $A$, in the image
of $X \longrightarrow \mathcal{M}_{\mathbf{1}}^{[a,b]}$.

Let $P$ be a point in $\mathcal{M}_{\mathbf{1}}^{\nu}(A)$, and
$A \longrightarrow \pi_{0}(A) \longrightarrow K$ a morphism
with $K$ a field, corresponding to a point
$x\in Spec\, \pi_{0}(A)$. We can chose a morphism
$$K^{\nu(b)}[-b] \longrightarrow P\otimes^{\mathbb{L}}_{N(A)}K$$
inducing a surjective morphism
$$g : K^{\nu(b)} \longrightarrow H^{b}(P\otimes^{\mathbb{L}}_{N(A)}K).$$
As $H^{i}(P)=0$ for all $i>b$, we have
$$H^{b}(P)\otimes_{\pi_{0}(A)}K\simeq H^{b}(P\otimes^{\mathbb{L}}_{N(A)}K),$$
and thus the morphism $g$ extends to a morphism of $N(A)$-dg-modules
$$u : A^{\nu(b)}[-b] \longrightarrow P.$$
We denote by $Q$ the homotopy cofiber of $u$, and
we consider the induced morphism
$Q \rightarrow A^{\nu(b)}[-b+1]$. The morphism $H^{b}(u)$ being surjective
at the point $x$, the object $Q$ lies in
$\mathcal{M}_{\mathbf{1}}^{\nu'}(A')$, for $A \rightarrow A'$ a Zariski
open neighborhood around $x$. This shows that
$P\otimes^{\mathbb{L}}_{N(A)}N(A')$ lies in the image
of $X(A')$. As the point $x$ was arbitrary, this finishes the proof that
$U \longrightarrow \mathcal{M}_{\mathbf{1}}^{\nu}$
is a covering, and thus that $\mathcal{M}_{\mathbf{1}}^{\nu}$
is quasi-compact.  \hfill $\Box$ \\

\subsection{Some consequences}

Let $T$ be a dg-category of finite type and 
$\mathcal{M}_{T}$ be the $D^{-}$-stack of pseudo-perfect $T^{op}$-dg-modules. We consider
the truncation $t_{0}\mathcal{M}_{T}$, and we define for any integer $n> 0$ a substack
$t_{ 0}\mathcal{M}_{T}^{n-rid} \subset t_{0}\mathcal{M}_{T}$ of \emph{$n$-rigid objects} 
in the following way:
for any $k$-algebra $k' \in k-CAlg$, the simplicial set 
$t_{0}\mathcal{M}_{T}^{n-rig}(k')$ is defined to be the full sub-simplicial set 
of $t_{0}\mathcal{M}_{T}(k')=\mathcal{M}_{T}(k')$ consisting of 
$T^{op}\otimes^{\mathbb{L}}k'$-dg-modules $E$ such that for any morphism $k' \longrightarrow k''$ in 
$k-CAlg$ we have
$$Ext^{-i}_{T^{op}\otimes^{\mathbb{L}}k''-Mod}(E\otimes_{k'}^{\mathbb{L}}k'',E\otimes_{k'}^{\mathbb{L}}k'')
=\pi_{i}(Map_{T^{op}\otimes^{\mathbb{L}}k''-Mod}(E\otimes_{k'}^{\mathbb{L}}k'',E\otimes_{k'}^{\mathbb{L}}k''))=0 
\qquad \forall \; i \geq  n.$$

\begin{cor}\label{cn-rig}
The stack $t_{0}\mathcal{M}_{T}^{n-rid}$ is an Artin n-stack.  
In particular $t_{0}\mathcal{M}_{T}^{1-rid}$ is
an Artin 1-stack. 
\end{cor}

\textit{Proof:} By definition and by the description of homotopy groups of
nerves of categories of equivalences in model categories (see \cite{hagII}), the stack 
$t_{0}\mathcal{M}_{T}^{n-rid}$ is an n-stack. Moreover, by semi-continuity of the dimension of the cohomology
groups of perfect complexes, it is easy to see that 
the inclusion morphism $t_{0}\mathcal{M}_{T}^{n-rid} \longrightarrow t_{0}\mathcal{M}_{T}$
is a Zariski open immersion. This implies that $t_{0}\mathcal{M}_{T}^{n-rid}$ is 
a locally geometric stack which is n-truncated, and thus that it is an Artin n-stack
by lemma \ref{l5'}. \hfill $\Box$ \\

We define a full substack $t_{0}\mathcal{M}_{T}^{simp}\subset t_{0}\mathcal{M}_{T}^{1-rig}$, 
of \emph{simple objects}, 
as follows: for any $k$-algebra $k' \in k-CAlg$, the simplicial set 
$t_{0}\mathcal{M}_{T}^{simp}(k')$ is defined to be the full sub-simplicial set 
of $t_{0}\mathcal{M}_{T}(k')^{1-rig}$ consisting of 
$T^{op}\otimes^{\mathbb{L}}k'$-dg-modules $E$ such that for any morphism $k' \longrightarrow k''$ in 
$k-CAlg$ the natural morphism
$$k'' \longrightarrow 
End_{Ho(T^{op}\otimes^{\mathbb{L}}k''-Mod)}(E\otimes_{k'}^{\mathbb{L}}k'')$$
is an isomorphism.  

\begin{cor}\label{csimp}
The stack $t_{0}\mathcal{M}_{T}^{simp}$ is an Artin 1-stack. Moreover, the
sheaf $\pi_{0}(t_{0}\mathcal{M}_{T}^{simp})$ is  
an algebraic space locally of finite presentation, and the natural morphism
$$t_{0}\mathcal{M}_{T}^{simp} \longrightarrow 
\pi_{0}(t_{0}\mathcal{M}_{T}^{simp})$$
is a locally trivial fibration (for the etale topology) 
with fibers $K(\mathbb{G}_{m},1)$.
\end{cor}

\textit{Proof:} It is easy to see that $t_{0}\mathcal{M}_{T}^{simp}$ is an open substack 
in $t_{0}\mathcal{M}_{T}^{1-rig}$ and so is an Artin 1-stack by corollary
\ref{cn-rig}. Moreover, let $M:=\pi_{0}(t_{0}\mathcal{M}_{T}^{simp})$. 
It is easy to see that the natural morphism
$$p : t_{0}\mathcal{M}_{T}^{simp} \longrightarrow M$$
is a torsor over the group stack $K(\mathbb{G}_{m},1)$. This implies that
for any affine k-scheme $Y$ and any morphism 
$Y \longrightarrow M$, the induced morphism 
$$t_{0}\mathcal{M}_{T}^{simp}\times_{M}^{h}Y \longrightarrow Y$$
is locally (for the etale topology) equivalent on $Y$ to the 
projection $Y\times K(\mathbb{G}_{m},1)\longrightarrow Y$. This implies that 
the morphism $p$ is 1-representable, smooth and surjective. 
By \cite{hagII} this implies that $M$ is n-geometric for some n, and thus
is an Artin 0-stack, or in other words an algebraic space. Moreover, 
as $M$ is locally a retract of $t_{0}\mathcal{M}_{T}^{simp}$, $M$ is also 
locally of finite presentation. \hfill $\Box$  \\

\begin{rmk}
\emph{For any n-geometric $D^{-}$-stack $F$, the natural closed
embedding $t_{0}(F) \longrightarrow F$ is a formal thickening and
induces an equivalence between the small Zariski sites
of $F$ and of $t_{0}(F)$. Therefore, the open substacks
$$t_{0}\mathcal{M}_{T}^{simp} \subset t_{0}\mathcal{M}_{T}^{n-rig}\subset 
t_{0}\mathcal{M}_{T}$$
correspond to open sub-$D^{-}$-stacks
$$\mathcal{M}_{T}^{simp} \subset \mathcal{M}_{T}^{n-rig}\subset 
\mathcal{M}_{T}.$$
This provides natural derived version of 
the stacks $t_{0}\mathcal{M}_{T}^{simp}$ and $t_{0}\mathcal{M}_{T}^{n-rig}$.}
\end{rmk}

Let $T$ be a saturated dg-category. We define a presheaf of groups on the site of
affine $k$-schemes with the \'etale topology as follows:
$$\begin{array}{cccc}
k-CAlg & \longrightarrow & Gp \\
 k' & \mapsto & Aut_{Ho(dg-Cat_{k'})}(\widehat{T^{op}\otimes^{\mathbb{L}}k'})\simeq Aut_{Ho(dg-Cat_{k'})}(\widehat{T^{op}\otimes^{\mathbb{L}}k'}_{pe}),
 \end{array}$$
where $Ho(dg-Cat_{k'})$ the homotopy category of $k'$-dg-categories, and
$k-CAlg$ is the category of commutative k-algebras. We will denote by 
$aut(T)$ the associated sheaf of groups. 

In order to state the next corollary remind from \cite{to} that for any dg-category
$T$ it is possible to define a Hochschild complex $HH(T)$ whose
cohomology groups are given by
$$HH^{i}(T)=Ext^{i}_{T\otimes^{\mathbb{L}}T^{op}}(T,T)\simeq 
H^{i}(\mathbb{R}\underline{Hom}(T,T)(id,id)).$$

\begin{cor}\label{caut}
Assume that one of the two conditions below are satisfied
\begin{enumerate}
\item We have $HH^{i}(T)=0$ for all $i<0$ and the natural morphism $k \longrightarrow HH^{0}(T)$
is an isomorphism. 
\item The ring $k$ is a field.
\end{enumerate} 
Then the sheaf of groups
$aut(T)$ is representable by an algebraic space locally of finite presentation over $k$
(and thus by a scheme if $k$ is a field). 
\end{cor}

\textit{Proof:}  We start to define a stack
$$\begin{array}{cccc}
\underline{End}(T) : & k-CAlg & \longrightarrow & SSet \\
 & k' & \mapsto & Map_{dg-Cat_{k'}}(\widehat{T\otimes^{\mathbb{L}}k'}_{pe},
\widehat{T\otimes^{\mathbb{L}}k'}_{pe}).
\end{array}$$
Using that $T$ is saturated and \cite{to}, it is easy to see that $\underline{End(T)}$
is equivalent to $t_{0}\mathcal{M}_{T\otimes^{\mathbb{L}} T^{op}}$. Therefore, it 
is a locally geometric stack by theorem \ref{t1}. We can also define 
a substack $\underline{Aut}(T)$ of $\underline{End}(T)$ as follows:
$$\begin{array}{cccc}
\underline{Aut}(T) : & k-CAlg & \longrightarrow & SSet \\
 & k' & \mapsto & Map_{dg-Cat_{k'}}^{eq}(\widehat{T\otimes^{\mathbb{L}}k'}_{pe},
\widehat{T\otimes^{\mathbb{L}}k'}_{pe}).
\end{array}$$

\begin{lem}\label{lcsimp}
Let $B$ be a smooth and proper dg-algebra over $k$.
Let $E$ be a perfect $B\otimes^{\mathbb{L}}B^{op}$-dg-module. 
Then, the induced morphism of dg-categories
$$F_{E} : \widehat{B^{op}}_{pe} \longrightarrow \widehat{B^{op}}_{pe}$$
is a quasi-equivalence if and only if the two natural morphisms of complexes of
$k$-modules
$$B \longrightarrow \mathbb{R}\underline{End}_{B-Mod}(E,E)$$
$$E\otimes_{B}^{\mathbb{L}}\mathbb{R}\underline{End}_{B-Mod}(E,B) \longrightarrow B$$
are quasi-isomorphisms. 
\end{lem}

\textit{Proof of the lemma:} The functor $F_{E}$ sends a $B$-dg-module $X$ to 
$E\otimes_{B}^{\mathbb{L}}X$. This implies that
$F_{E}$ is fully faithful if and only if the morphism
$$B \longrightarrow \mathbb{R}\underline{End}_{B-Mod}(E,E)$$
is a quasi-isomorphism (because $B$ is a generator of 
the category of $B$-dg-modules). The functor $F_{E}$ has a 
right adjoint $G_{E}$ sending a $B$-dg-module $X$ to 
$\mathbb{R}\underline{End}_{B-Mod}(E,X)$. Therefore, 
$G_{E}$ is fully faithful if and only if the morphism
$$E\otimes_{B}^{\mathbb{L}}\mathbb{R}\underline{End}_{B-Mod}(E,B) \longrightarrow B$$
is a quasi-isomorphism. 
\hfill $\Box$ \\

The lemma \ref{lcsimp}, and the semi-continuity of the dimension of
the cohomology groups of perfect complexes, imply that 
the natural embedding 
$$\underline{Aut}(T) \longrightarrow \underline{End}(T)$$
is a Zariski open immersion (this follows proposition \ref{psupp}
applied to cones of the two morphisms of lemma \ref{lcsimp}). In particular, $\underline{Aut}(T)$
is a locally geometric stack. 

By definition $aut(T)$ is the sheaf $\pi_{0}(\underline{Aut}(T))$. Let 
$p : \underline{Aut}(T) \longrightarrow aut(T)$ be the natural morphism, and let us
prove that $p$ is n-representable for some n, and smooth. By 
\cite[Cor. 1.3.4.5]{hagII} this will imply that $aut(T)$ is an algebraic space locally of finite presentation.
For this, let $X=Spec\, k'$ be an affine $k$-scheme and $X\longrightarrow aut(T)$
be a morphism. We want to prove that the stack 
$\underline{Aut}(T)\times_{aut(T)}^{h}X$ is n-geometric for some $n$ and that 
the morphism $\underline{Aut}(T)\times_{aut(T)}^{h}X \longrightarrow X$
is smooth. As this is local on $X$, we can assume that the morphism 
$X \longrightarrow aut(T)$ lifts to a morphism $X\longrightarrow  \underline{Aut}(T)$.

The morphism $X=Spec\, k' \longrightarrow \underline{Aut}(T)$ correspond to 
a certain perfect $T\otimes^{\mathbb{L}}T^{op}$-dg-module $E$. By definition, 
the stack $\underline{Aut}(T)\times_{aut(T)}^{h}X$ can be 
written as $K(\underline{Aut}(E),1)$, where $\underline{Aut}(E)$ is the 
group stack defined by 
$$\begin{array}{cccc}
\underline{Aut}(E) : & k'-CAlg  & \mapsto & SSet \\
 & k'' & \mapsto & Map^{eq}_{T\otimes^{\mathbb{L}}T^{op}\otimes^{\mathbb{L}}k''}
(E\otimes^{\mathbb{L}}k'',E\otimes^{\mathbb{L}}k'').
\end{array}$$
Therefore, it is enough to show that $\underline{Aut}(E)$ is an
n-geometric stack for some n and is smooth over $X$. The fact that 
$\underline{Aut}(E)$ is n-geometric for some n follows from 
the local geometricity of $t_{0}\mathcal{M}_{T\otimes^{\mathbb{L}}T^{op}}$, as
we have
$$\underline{Aut}(E) \simeq X\times_{t_{0}\mathcal{M}_{T\otimes^{\mathbb{L}}T^{op}}}^{h}X.$$
It remains to show that it is smooth over $X$.

As $\underline{Aut}(T)$ is a group stack, the stack 
$\underline{Aut}(E)$ is equivalent to $\underline{Aut}(I)$, where
$I$ is $T\otimes^{\mathbb{L}}k'$ considered as a 
$T\otimes^{\mathbb{L}}T^{op}\otimes^{\mathbb{L}}k'$-dg-module. 
The stack $\underline{Aut}(I)$ is a Zariski open substack 
of $\underline{End}(I)$, defined naturally by
$$\begin{array}{cccc}
\underline{End}(I) : & k'-CAlg  & \mapsto & SSet \\
 & k'' & \mapsto & Map_{T\otimes^{\mathbb{L}}T^{op}\otimes^{\mathbb{L}}k''}
(T\otimes^{\mathbb{L}}k'',T\otimes^{\mathbb{L}}k'').
\end{array}$$
It only remains to show that $\underline{End}(I)$ is smooth over $X$. 
But, under condition $(1)$, we have 
$\underline{End}(I)\simeq \mathbb{G}_{a}$. An under condition $(2)$
it is easy to see that $\underline{End}(I)$ is equivalent, as a stack over $X$, to 
$$\underline{End}(I)\simeq \prod_{i\leq 0} \mathbb{A}^{n_{i}}_{X},$$
where $n_{i}:=dim_{k}HH^{i}(T)$ (this uses that the complex of $k$-vector spaces
$HH(T)$ is quasi-isomorphic to a direct sum of its cohomology groups). 
\hfill $\Box$ \\

It follows from proposition \ref{p6} and corollary \ref{caut}, that for any saturated
dg-category $T$ over a field 
$k$, the group scheme $aut(T)$ has only countably many connected components. 
We therefore obtain the following more precise statement.

\begin{cor}\label{caut2}
Let $T$ be a dg-category saturated over a field $k$. Then, the group 
scheme $aut(T)$ is an extension 
$$\xymatrix{
1 \ar[r] & aut(T)_{e} \ar[r] & aut(T) \ar[r] & \Gamma \ar[r] & 1,}$$
where $\Gamma$ is a countable discrete group, and 
$aut(T)_{e}$ is an algebraic group of finite type over $k$. 
\end{cor}

\textit{Proof:} This follows from the fact that $aut(T)$ has only 
countably many connected components, and by the fact
a connected group scheme locally of finite type over a field $k$ is
quasi-compact (see \cite[Exp. $VI_{B}$, Cor. 3.6]{sga3}). \hfill $\Box$ \\

\subsection{Two examples}

To finish this section, we will present two
fundamental examples of dg-categories of finite type
coming from algebraic geometry and representation theory.
They will give by theorem \ref{t1} two
locally geometric $D^{-}$-stacks, one classifying perfect
complexes on a smooth and proper schemes, and a second one
classifying complexes of representation of a finite Quiver. They both seem
to us important examples to study in the future. \\

Let $p : X\longrightarrow Spec\, k$ be a smooth and proper morphism of schemes.
We consider the category $QCoh(X)$ of quasi-coherent sheaves on $X$, as
well as $CQCoh(X)$ the category of unbounded complexes in
$QCoh(X)$. As the scheme $X$ is quasi-compact and separated, it is known by
\cite{ho2} that $CQCoh(X)$ is a model category for which cofibrations are monomorphisms
and equivalences are quasi-isomorphisms. This model category has a natural
enrichment in $C(k)$, by tensoring a complex of quasi-coherent sheaves
by a complex of $k$-modules. This makes $CQCoh(X)$ into a $C(k)$-model category, 
and we set
$$L_{qcoh}(X):=Int(CQCoh(X)).$$
Note that $L_{qcoh}(X)$ is only a $\mathbb{V}$-small dg-category, but with
$\mathbb{U}$-small complexes of morphisms between two fixed objects.
We denote by $L_{pe}(X)$ the full sub-dg-category of
$L_{qcoh}(X)$ consisting of perfect complexes. The dg-category
$L_{pe}(X)$ is still only $\mathbb{V}$-small, but is quasi-equivalent to
a $\mathbb{U}$-small dg-category as the set of isomorphisms classes
of objects in $[L_{pe}(X)]$ is $\mathbb{U}$-small. We will therefore
do as if $L_{pe}(X)$ were a $\mathbb{U}$-small dg-category.

The triangulated category $[L_{qcoh}(X)]$ is naturally equivalent to the
unbounded derived category $D_{qcoh}(X)$ of quasi-coherent sheaves on $X$.
In the same way, $[L_{pe}(X)]$ is naturally equivalent
to $D_{pe}(X)$, the full sub-category of $D_{qcoh}(X)$ consisting of
perfect complexes. By \cite{bv}, we know that
$D_{pe}(X)$ is precisely the sub-category of compact objects in
$D_{qcoh}(X)$, and furthermore that $D_{qcoh}(X)$ is compactly generated.
This easily implies that the restricted Yoneda embedding
$$L_{qcoh}(X) \longrightarrow \widehat{L_{pe}(X)}$$
is a quasi-equivalence. Furthermore, $L_{pe}(X)$ is a triangulated
dg-category, and therefore 
$$L_{pe}(X)\simeq \widehat{L_{pe}(X)}_{pe}\simeq (L_{qcoh}(X))_{pe}.$$

\begin{lem}\label{l12}
The dg-category $L_{pe}(X)$ is saturated.
\end{lem}

\textit{Proof:} We already know that
$L_{pe}(X)$ is triangulated. As the scheme
$X$ is proper and flat over $Spec\, k$, for two
perfect complexes $E$ and $F$ on $X$,
the complex of $k$-modules $\mathbb{R}\underline{Hom}(E,F)$
is perfect. This implies that
$L_{pe}(X)$ is locally perfect. By \cite{bv}, it is also
known that $L_{qcoh}(X)=\widehat{L_{pe}(X)}$ has a compact generator,
and thus that $L_{pe}(X)$ is a proper dg-category. Finally,
by \cite[\S 8.3]{to}, we have $\widehat{L_{pe}(X)\otimes^{\mathbb{L}} L_{pe}(X)^{op}}\simeq
L_{qcoh}(X\times_{k}X)$, and the diagonal $(L_{pe}(X)^{op}\otimes^{\mathbb{L}} L_{pe}(X))$-dg-module
$$\begin{array}{cccc}
\delta : & L_{pe}(X)^{op}\otimes^{\mathbb{L}} L_{pe}(X) & \longrightarrow & C(k) \\
 & (x,y) & \mapsto & L_{pe}(X)(x,y)
\end{array}$$
corresponds by this quasi-equivalence to $\Delta_{X}$, the class of the diagonal in
$L_{qcoh}(X\times_{k}X)$. As $X$ is a smooth and separated over $Spec\, k$, we know
that $\Delta_{X}$ is a perfect complex on $X\times_{k}X$, and thus
lies in $L_{qcoh}(X\times_{k}X)_{pe}$. By definition, this means that
$L_{pe}(X)$ is a smooth dg-category. \hfill $\Box$ \\

\begin{df}\label{d10}
The $D^{-}$-stack $\mathcal{M}_{L_{pe}(X)}$ is called the \emph{$D^{-}$-stack
of perfect complexes on $X$}. It is denoted by
$\mathbb{R}\underline{Perf}(X)$, and its truncation by
$\underline{Perf}(X):=t_{0}(\underline{Perf}(X))$.
\end{df}

For any
commutative $k$-algebra $A$, the simplicial set
$\mathbb{R}\underline{Perf}(X)(A)$ is by definition
$Map_{dg-Cat}(L_{pe}(X)^{op},\widehat{A}_{pe})$, which in turns is equivalent
to the nerve of the category of quasi-isomorphisms between
perfect complexes on $X\times^{h} Spec\, A$ (see \cite[\S 8.3]{to}). 
In other words, if 
$\mathbb{R}\underline{Perf}$ denotes 
the $D^{-}$-stack $\mathcal{M}_{\mathbf{1}}$, which is
the $D^{-}$-stack of perfect complexes, we have
$$\mathbb{R}\underline{Perf}(X)\simeq \mathbf{Map}(X,\mathbb{R}\underline{Perf}),$$
where $\mathbf{Map}$ denotes the internal $Hom$ of the category
$St(k)$ (see \cite{hagII}). This justifies the
terminology of definition \ref{d10}. 
In particular,
the truncation $\underline{Perf}(X):=t_{0}\mathbb{R}\underline{Perf}(X)$ is the
(un-derived) stack of perfect complexes on $X$, whose absolute
version has been considered
previously in \cite[\S 21]{hs}. \\

Lemma \ref{l12} and theorem \ref{t1} imply  the following corollary.

\begin{cor}\label{c5}
The $D^{-}$-stack $\mathbb{R}\underline{Perf}(X)$ is
locally geometric and locally of finite type. For any
global point $E\in \mathbb{R}\underline{Perf}(X)(k)$, we have
$$\mathbb{T}_{\mathbb{R}\underline{Perf}(X),E}\simeq \mathbb{R}\underline{Hom}(E,E)[1].$$
\end{cor}

\begin{rmk}
\begin{enumerate}
\item \emph{From corollaries \ref{cn-rig}, \ref{csimp} and \ref{c5}
we find a new proof of the existence of the Artin 1-stack $\underline{Coh}(X)$ of
coherent sheaves on a smooth and proper scheme (e.g. as 
presented in \cite{lm}). More generally, we also find a new
proof of the existence of the Artin 1-stack of 1-rigid perfect
complexes on $X$, previously constructed in \cite{l}.
It seems to us important to notice that 
this new proof uses in an essential way higher stacks, as
it is based on the morphism $\pi$ from $\underline{Coh}(X)$
to the stack of perfect complexes, sending a coherent sheaf
$E$ on $X$ to the complex $\mathbb{R}\underline{Hom}(E_{0},E)$, 
for $E_{0}$ a compact generator in $D_{qcoh}(X)$.}
\item
\emph{
When $X$ is no longer proper and smooth, 
the $D^{-}$-stack $\mathcal{M}_{L_{pe}(X)}$ only classifies
pseudo-perfect complexes on $X$. These are the 
complexes of quasi-coherent sheaves $E$ on $X$ such that 
for any perfect complexes $F$ on $X$ the 
complex of $k$-modules $\mathbb{R}\underline{Hom}(F,E)$
is perfect. When $X$ is only assumed to be smooth
then a pseudo-perfect complex if perfect (by \ref{l2}), and
when $X$ is only assumed to be proper and flat then 
perfect complexes are pseudo-perfect (by \ref{l2}). 
In general, the objects classified by 
$\mathcal{M}_{L_{pe}(X)}$ have no clear relations
with perfect complexes on $X$, and this is the reason
why we did not consider this general
situation, even if corollary 
\ref{c5} stays most probably correct under much more
general assumptions than smoothness and properness. }
\end{enumerate}
\end{rmk}

Let $E\in D_{qcoh}(X)$ be a compact generator, and let us consider
the induced morphism
$$\pi : \mathbb{R}\underline{Perf}(X) \longrightarrow \mathbb{R}\underline{Perf},$$
as during the proof
of theorem \ref{t1}. The morphism $\pi$ simply sends a 
perfect complex $F$ on $X$ to the perfect complex
$\mathbb{R}\underline{Hom}(E,F)$.
Let $\nu : \mathbb{Z}\longrightarrow \mathbb{N}$ be a
function with finite support, and let us consider
the sub-$D^{-}$-stack $\mathbb{R}\underline{Perf}(X)^{\nu}$ of
$\mathbb{R}\underline{Perf}(X)$ defined as in \ref{p6}.

\begin{cor}\label{c6}
The $D^{-}$-stack $\mathbb{R}\underline{Perf}(X)^{\nu}$ is
$n$-geometric for some $n$, and strongly of finite presentation.
\end{cor}

Corollaries \ref{c5} and \ref{c6} provide in particular a proof of 
the result \cite[Thm. 21.5]{hs}. \\

Corollary \ref{c6} has an interesting consequences when 
$X$ is a smooth and projective variety over a field $k$, as in this case
there is a natural choice for the generator $E$ of $D_{qcoh}(X)$. 
Indeed, recall from \cite{bv} that there exists an 
integers $d$ such that the vector bundle 
$E:=\oplus_{i\in [0,d]}\mathcal{O}(-i)$ is 
a compact generator of $D_{qcoh}(X)$, where
$\mathcal{O}(1)$ is a very ample line bundle on $X$. 
>From corollary 
\ref{c6} we thus get a boundeness
criterion for perfect complexes on $X$, which 
seems to be a generalization of Kleiman criterion. 
For this, we will say that a family of
objects $\{E_{i}\}_{i\in I}$ in $D_{pe}(X)$ is bounded, if 
the corresponding family of points $E_{i}\in \mathbb{R}\underline{Parf}(X)(k)$
belongs to an open sub-$D^{-}$-stack which is 
strongly of finite presentation.

\begin{cor}\label{c6+}
Let $k$ be a field and $X$ be a smooth and projective 
variety over $k$. Let $\mathcal{O}(1)$ be a very ample
line bundle on $X$. Then, there exists an integer $d$, such
that the following condition is satisfied:

A family of perfect complexes $\{E_{i}\}_{i\in I}$ on $X$ is bounded
if and only if there exists a function $\nu : \mathbb{Z} \longrightarrow \mathbb{N}$
with finite support, such that 
$$Dim_{k}\mathbb{H}^{k}(X,E_{i}(j)) \leq \nu(k) \qquad \forall k\; \forall 
j\in [0,d] \; \forall i\in I.$$
\end{cor}

\textit{Proof:} This follows from the fact that 
any open sub-$D^{-}$-stack strongly of finite presentation 
of $\mathbb{R}\underline{Perf}(X)$ is contained in some
$\mathbb{R}\underline{Perf}(X)^{\nu}$, together with the fact that 
for the generator $E:=\oplus_{i\in [0,d]}\mathcal{O}(-i)$
the morphism 
$$\pi : \mathbb{R}\underline{Perf}(X) \longrightarrow \mathbb{R}\underline{Perf}$$
sends a perfect complex $F$ on $X$ to 
$\oplus_{i\in [0,d]}\mathbb{R}\Gamma(X,F(-i))$. \hfill $\Box$ \\

We now pass to our second example. Let $B$ be an associative $k$-algebra, which is
assumed to be projective and of finite type as a $k$-module. We also assume that
$B$ is perfect as a complex of $B\otimes_{k}B^{op}$-modules. Then, the
dg-category with a unique object $B$ is smooth and proper, and thus
by theorem \ref{t1} the $D^{-}$-stack $\mathcal{M}_{B}$ is locally geometric.
By definition, for any $A\in sk-CAlg$, there exists a natural equivalence
between $\mathcal{M}_{B}(A)$ and the nerve of the category of quasi-isomorphisms
between perfect $B\otimes_{k}N(A)$-dg-modules. Therefore,
$\mathcal{M}_{B}$ is a moduli for perfect complexes of $B$-modules.

Now, let us assume that $Q$ is an oriented finite Quiver.
Let $B=k[Q]^{op}$ be the opposite $k$-algebra associated to $Q$. The $B^{op}$-modules are precisely
the representations of $Q$ in $k$-modules, and thus
the dg-category $\widehat{B}$ can be identified with the
dg-category of complexes of representations of $Q$ with coefficients in
$k$-modules. The pseudo-perfect objects in $\widehat{B}$
are simply the complexes of representations of $B$ whose
underlying complex of $k$-modules is perfect.

\begin{df}\label{d11}
With the notations above, the $D^{-}$-stack $\mathcal{M}_{B}$ is called
the \emph{$D^{-}$-stack of pseudo-perfect complexes of representations of $Q$}. 
It is denoted
by $\mathbb{R}\underline{Parf}(Q)$.
\end{df}

The $k$-algebra $B$ is easily seen to be 
homotopically finitely presented in the model
category of dg-algebras. If the quiver $Q$ has no loops, then
$B$ is furthermore a smooth and proper dg-category. Therefore, theorem
\ref{t1} implies the following corollary.

\begin{cor}\label{c7}
The $D^{-}$-stack $\mathbb{R}\underline{Parf}(Q)$ is locally geometric
and locally of finite type.
For a global point $E\in \mathbb{R}\underline{Parf}(Q)(k)$, corresponding to
a complex of representations of $Q$, we have
$$T_{\mathbb{R}\underline{Parf}(Q),E}\simeq \mathbb{R}\underline{Hom}(E,E)[1],$$
where $\mathbb{R}\underline{Hom}$ denotes the derived
$Hom$'s of complexes of representations of $Q$.
\end{cor}

One has a natural projection, corresponding to the natural
generator of $\widehat{B}$ (i.e. $B$ itself)
$$\mathbb{R}\underline{Parf}(Q) \longrightarrow \mathbb{R}\underline{Perf}.$$
As in proposition \ref{p6}, 
one gets for any $\nu : \mathbb{Z} \longrightarrow \mathbb{N}$
with finite support, an open sub-$D^{-}$-stack
$\mathbb{R}\underline{Parf}(Q)^{\nu}$.

\begin{cor}\label{c8}
The $D^{-}$-stack $\mathbb{R}\underline{Parf}(Q)^{\nu}$
is $n$-geometric for some $n$, and strongly of finite presentation.
\end{cor}

\end{document}